\documentclass[EJP, preprint, NODS, nobera, nofontexpansion]{ejpecp}

\usepackage{amsmath}
\usepackage{amssymb}
\usepackage{mathrsfs}
\usepackage{adjustbox}
\usepackage{array}
\usepackage{longtable}
\usepackage{multirow}
\usepackage{arydshln}
\usepackage{pdflscape}
\usepackage{afterpage}
\usepackage{capt-of}
\usepackage[T1]{fontenc}

\makeatletter
\patchcmd{\AP@savetop}{\par}{\@@par}{}{}
\AddToHook{shipout/after}{%
  \ifnum\value{page}=5\relax
    \afterpage{%
      \let\par\@@par
      \begin{landscape}
\noindent\begin{minipage}[t]{\linewidth}
\vspace{0pt}
\captionof{table}{Comparison rates by profile regime and test geometry.  Each
entry is an $O(\cdot)$ bound. Blank assumption cells mean that no
additional assumption of that type is imposed in the displayed
specialization. A slash marks a cell where no separate bound in the
corresponding quadratic-form modulus is recorded. The number $k$
of classes in the last two rows may grow. All displayed operator
specializations have $|\Sigma|_*^\infty=O(n)$; the general theorem retains
an explicit trace factor and imposes no such restriction.}
\label{tab:comparison_rates}
\centering
{\setlength{\tabcolsep}{3pt}
\setlength{\extrarowheight}{2pt}
\renewcommand{\arraystretch}{1.35}
\begin{adjustbox}{max width=\linewidth,max totalheight=6.3cm,center}
\begin{tabular}{|V|V:V:V:V:V|V:V:V:V:V:V|}
\hline
\multirow{3}{*}{\textbf{Result}}
& \multicolumn{5}{c|}{\textbf{Assumptions}}
& \multicolumn{6}{c|}{\textbf{Rates for }
$\boldsymbol{\|\tr(A(G^z-\tilde G^z))\|_{L^q}}$} \\
\cline{2-12}
& \multirow{2}{*}{\shortstack{\textbf{Profile}\\\textbf{structure}}}
& \multirow{2}{*}{$\boldsymbol{p/n}$}
& \multirow{2}{*}{$\boldsymbol{|\Sigma|_{\op}^\infty}$}
& \multirow{2}{*}{$\boldsymbol{|\Sigma|_{\hs}^\infty}$}
& \multirow{2}{*}{$\boldsymbol{|\Sigma|_*^\infty}$}
& \multicolumn{2}{c:}{$\boldsymbol{|A|_{\op}\leq1}$}
& \multicolumn{2}{c:}{$\boldsymbol{|A|_{\hs}\leq1}$}
& \multicolumn{2}{c|}{\textbf{Stieltjes t. }$\boldsymbol{(A=I_p/p)}$} \\
\cline{7-12}
& & & & &
& $\boldsymbol{q\geq1}$ & $\boldsymbol{q\geq2}$
& $\boldsymbol{q\geq1}$ & $\boldsymbol{q\geq2}$
& $\boldsymbol{q\geq1}$ & $\boldsymbol{q\geq2}$ \\
\hline\hline
\shortstack{Th.~\ref{the:operator_profile_comparison}\\Th.~\ref{the:noncommuting_limited_second_moment_profiles} for $|\cdot|_{\hs}$}
& Arbitrary
& $O(1)$
& $O(1)$
& 
& 
& \multirow{2}[6]{*}[0.5ex]{$1+M_{\op}^{(q)}$}
& \multirow{2}[6]{*}[0.5ex]{$\displaystyle\begin{gathered}\frac{M_{\op}^{(q)}}{\sqrt n}\\[-2pt]{}+\frac{1+(M_{\op}^{(2)})^2}{n}\end{gathered}$}
& $\displaystyle M_{\hs}^{(q)}+\frac{\sqrt p}{n}$
& $\displaystyle\begin{gathered}\frac{M_{\hs}^{(q)}}{\sqrt n}\\[-2pt]{}+\frac{[1+(M_{\hs}^{(2)})^2]\sqrt p}{n}\end{gathered}$
& \multirow{2}[6]{*}[0.5ex]{$\displaystyle\frac{1+M_{\op}^{(q)}}p$}
& \multirow{2}[6]{*}[0.5ex]{$\displaystyle\begin{gathered}\frac{M_{\op}^{(q)}}{p\sqrt n}\\[-2pt]{}+\frac{1+(M_{\op}^{(2)})^2}{pn}\end{gathered}$} \\
\cline{1-6}\cline{9-10}
Th.~\ref{the:operator_profile_comparison}
& Arbitrary
& 
& $O(1)$
& 
& $O(n)$
& 
& 
& $\diagup$
& $\diagup$
& 
&  \\
\hline
Th.~\ref{the:operator_profile_comparison}
& \shortstack{$k$ well-balanced\\classes\\$\left(\frac{1}{\min_a\pi_a}=O(k)\right)$}
& 
& 
& 
& $O(n)$
& $k(1+M_{\op}^{(q)})$
& $\displaystyle\begin{gathered}\frac{kM_{\op}^{(q)}}{\sqrt n}\\[-2pt]{}+\frac{k^2+(M_{\op}^{(2)})^2}{n}\end{gathered}$
& $\diagup$
& $\diagup$
& $\displaystyle\frac{k(1+M_{\op}^{(q)})}p$
& $\displaystyle\begin{gathered}\frac{kM_{\op}^{(q)}}{p\sqrt n}\\[-2pt]{}+\frac{k^2+(M_{\op}^{(2)})^2}{pn}\end{gathered}$ \\
\hline
Th.~\ref{the:noncommuting_limited_second_moment_profiles}
& \shortstack{$k$ commuting\\profiles\\$(\Sigma^{(a)}\Sigma^{(b)}=\Sigma^{(b)}\Sigma^{(a)})$}
& $O(1)$
& 
& $O(\sqrt p)$
& 
& $\diagup$
& $\diagup$
& $\displaystyle k\left(M_{\hs}^{(q)}+\frac{\sqrt p}{n}\right)$
& $\displaystyle\begin{gathered}\frac{M_{\hs}^{(q)}}{\sqrt n}\\[-2pt]{}+\frac{k[1+(M_{\hs}^{(2)})^2]\sqrt p}{n}\end{gathered}$
& $\displaystyle k\left(\frac{M_{\hs}^{(q)}}{\sqrt p}+\frac1n\right)$
& $\displaystyle\begin{gathered}\frac{M_{\hs}^{(q)}}{\sqrt{np}}\\[-2pt]{}+\frac{k[1+(M_{\hs}^{(2)})^2]}{n}\end{gathered}$ \\
\hline
\end{tabular}%
\end{adjustbox}}
\end{minipage}

\vspace{0.9cm}
\noindent
\begin{minipage}[c]{0.565\linewidth}
    \footnotesize
    \setlength{\abovecaptionskip}{0pt}
    \captionof{figure}{Empirical eigenvalue density of $\frac1n XX^\top$ for
    columns $x_i=\sqrt{r_i}\,y_i$, where the $y_i$ are independent and
    uniform on $\{y\in\mathbb R^p:|y|=\sqrt p\}$, the $r_i$ are independent
    with $\mathbb P(r_i=1/5)=\mathbb P(r_i=9/5)=1/2$, and the two families
    are independent.  Thus $|x_i|^2/p=r_i$,
    $\mathbb E[x_ix_i^\top]=I_p$, and
    $\mathbb E[x_ix_i^\top\mid r_i]=r_iI_p$.  The histogram pools $3000$
    independent matrices with $p=50$ and $n=100$.  The dark dashed curve is
    the Marchenko--Pastur prediction based on the unconditional second
    moment, whereas the solid blue curve is the equal-weight two-profile
    prediction with profile values $I_p/5$ and $9I_p/5$.}
    \label{fig:radial_mixture_two_profile}
\end{minipage}\hfill
\begin{minipage}[c]{0.395\linewidth}
    \centering
    \includegraphics[width=\linewidth]
        {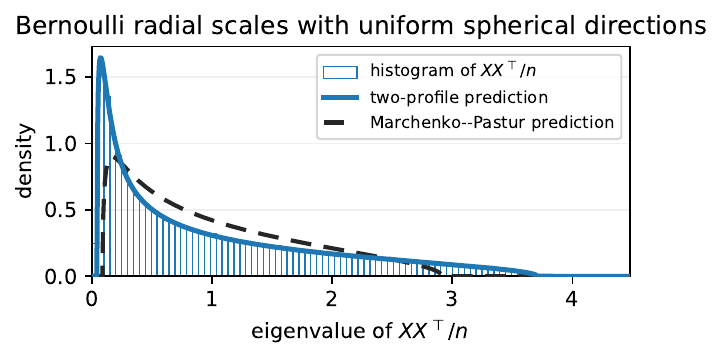}
\end{minipage}
\end{landscape}
    }%
  \fi
}
\makeatother

\newcolumntype{V}{>{\adjustbox{valign=m}\bgroup}c<{\egroup}}
\newcommand{\ii}{\mathbf{i}}
\DeclareMathOperator{\op}{op}
\DeclareMathOperator{\tr}{Tr}
\DeclareMathOperator{\diag}{Diag}

\DeclareMathOperator{\cone}{cone}
\DeclareMathOperator{\hs}{HS}
\DeclareMathOperator{\Var}{\mathbb V}
\SHORTTITLE{Resolvent convergence for Heterogeneous profiles}
\TITLE{Resolvent Convergence under \\ Heterogeneous Second-Moment Profiles \\ and Quadratic-Form Control}
\AUTHORS{Cosme Louart\footnote{School of Data Science, The Chinese University of Hong Kong (Shenzhen), Shenzhen, China}}
\KEYWORDS{sample covariance matrices; deterministic equivalents; resolvent estimates; heterogeneous covariance profiles; dependent coordinates; nonasymptotic bounds}

\AMSSUBJ{60B20, 15B52, 62H10}

\SUBMITTED{}

\ABSTRACT{
Let $X=(x_1,\ldots,x_n)\in\mathbb R^{p\times n}$ have independent real
columns with finite second moments, not necessarily centered or identically
distributed. For $z\in\mathbb C\setminus[0,\infty)$, we compare
\[
G^z=\left(\frac1nXX^\top-zI_p\right)^{-1}
\]
with a deterministic matrix $\tilde G^z$ defined by the second-moment
profile $(\Sigma_i)_{i\in[n]}$, where
$\Sigma_i=\mathbb E[x_ix_i^\top]$.

For deterministic matrices \(A\), we prove finite-dimensional \(L^q\) bounds for
$\tr(A(G^z-\tilde G^z))$ in terms of moments of centered quadratic forms,
without assuming independence among the coordinates of a column.
Operator-test bounds impose no trace-growth or aspect-ratio restriction
and display explicitly the effect of the sizes and repetitions of the
second-moment matrices $(\Sigma_i)_{i\in[n]}$. 
Under
$\sqrt p\sup_{i\in[n]}|\Sigma_i|_{\hs}=O(n)$, the Hilbert–Schmidt bounds involve approximation of the profile by commuting positive-semidefinite matrices, with explicit dependence on the approximation error and the dimension of their span.
We also obtain corresponding bounds for \(\tr(A G^z B G^z)\), with deterministic \(A,B\). All bounds are locally uniform in
$z\in\mathbb C\setminus[0,\infty)$.
}
\begin{document}

\section*{Introduction}
\addcontentsline{toc}{section}{Introduction}

We study deterministic approximations of the resolvent of a sample
second-moment matrix. Let $X=(x_1,\ldots,x_n)\in\mathbb R^{p\times n}$
have independent columns with finite second moments, not necessarily
centered or identically distributed, and set
\[
    \Sigma_i:=\mathbb E[x_ix_i^\top],
    \qquad
    G^z:=\left(\frac1nXX^\top-zI_p\right)^{-1},
    \quad z\in\mathbb C\setminus[0,\infty).
\]
The collection $(\Sigma_i)_{i\in[n]}$ is called the second-moment profile.
We approximate $G^z$ by a deterministic matrix $\tilde G^z$ depending only
on this collection, and bound $\tr(A(G^z-\tilde G^z))$ in $L^q$ for
deterministic $A$. This includes normalized traces, which describe the
empirical spectral distribution, and directional quantities obtained
with $A=uu^\top$.

The Marchenko--Pastur theorem \cite{MarchenkoPastur1967} and its extensions
provide the classical description of sample covariance spectra; see
\cite{BaiSilverstein2010,PasturShcherbina2011}.
Deterministic equivalents provide approximations defined at each matrix
dimension, notably for variance-profile and information-plus-noise models
\cite{DozierSilverstein2007,HachemLoubatonNajim2007}.
Concentration of quadratic forms provides a way to go beyond independence
of coordinates, as in Bai--Zhou \cite{BaiZhou2008} and Yaskov
\cite{Yaskov2016ECP}. Related results treat heterogeneous second-moment
matrices \cite{BenaychGeorgesCouillet2016,MeiWangYao2023} and quantitative
resolvent estimates under concentration assumptions
\cite{Chouard2022,FanMaPaquetteWang2026}.

Here we obtain finite-dimensional bounds that separate the dependence on
the matrices $\Sigma_i$, the test $A$, and the fluctuations of the centered
quadratic forms
\begin{align}
\label{eq:quadratic_form}
    x_i^\top Bx_i-\tr(\Sigma_iB).
\end{align}
The matrices $\Sigma_i$ may all differ and need not commute, and no limiting
assumption is imposed on their collection. Independence is required only
between columns. The estimates give convergence criteria according to the
size of these quadratic-form fluctuations, uniformly for $z$ in a fixed
compact subset of $\mathbb C\setminus[0,\infty)$.

The results are organized along two axes. The first is the norm used to
measure $A$. The Hilbert--Schmidt norm comparison combines $|A|_{\hs}$, bounds
on $|\Sigma_i|_{\hs}$, and moments of quadratic forms~\eqref{eq:quadratic_form} with
$|B|_{\hs}\leq1$. The operator norm comparison instead uses $|A|_{\op}$,
the nuclear norms $|\Sigma_i|_*=\tr(\Sigma_i)$, and the corresponding
moments for $|B|_{\op}\leq1$. The former can exploit smaller quadratic-form
fluctuations, while the latter can benefit from a smaller test operator norm and
allows an unrestricted aspect ratio $p/n$. The rates also reflect
relations within the second-moment profile $(\Sigma_i)_{i\in[n]}$: groups of equal
matrices enter the operator bound, whereas approximation by commuting
matrices enters the Hilbert--Schmidt bound.

The second axis is the moment order $q$ of the centered quadratic forms in \eqref{eq:quadratic_form}.
For each norm, we give bounds for $q\geq1$ and refined bounds for $q\geq2$.
The bounds with $q=1$, do not require the centered
quadratic forms to have finite second moment. 
At $q=2$, uniformly bounded quadratic-form second moments yield a gain of order $\sqrt n$ in the
leading moment-dependent term.
For i.i.d. entries with a fixed centered unit-variance law, the first-moment argument recovers weak convergence in probability to the Marchenko–Pastur law as \(p\to\infty\) and \(p/n\to c\in[0,\infty)\), without higher entry moments (Corollary~\ref{cor:mp_iid_finite_variance}).

The proofs combine leave-one-out
identities with stability of the deterministic fixed-point equations.
Successive class replacements handle heterogeneous collections
$(\Sigma_i)_{i\in[n]}$ in the operator norm comparison. This replacement procedure
produces resolvents with general matrices $Z$ in place of $zI_p$.
We therefore establish our results for matrix parameters $Z$ whose
numerical range is separated from $\mathbb R_+$. Additive
deformation and Cauchy's formula then also give deterministic
approximations for resolvent sandwiches $G^zBG^z$, with arbitrary
deterministic complex $B$. At negative spectral parameters, resolvents
are regularized inverses; these sandwiches occur, in particular, in
ridge-regression prediction errors.

Section~\ref{sec:main_results} states the results, compares them with
related work, and gives examples, with selected rates summarized in
Table~\ref{tab:comparison_rates}.
Sections~\ref{sec:deterministic_study} and~\ref{sec:probabilistic_results}
contain the deterministic and probabilistic proofs, respectively.
Subsection~\ref{sec:linear_response_superoperator} treats resolvent
sandwiches; auxiliary results are collected in the appendices.

\section*{Notation}
\label{sec:notation}
\addcontentsline{toc}{section}{Notation}
Set
\[
    \mathbb R_+:=[0,\infty),
    \qquad
    \mathbb H:=\{z\in\mathbb C:\Im z>0\}.
\]
Let $\delta_0$ denote the Dirac mass at zero.
For $S\subset\mathbb C$, its convex conic hull is denoted by
\[
    \cone(S)
    :=\left\{
        \sum_{\ell=1}^m a_\ell s_\ell:
        m\in\mathbb N,\ a_\ell\in\mathbb R_+,\ s_\ell\in S
      \right\}.
\]
For $q,r>0$, $q\wedge r$ denotes their minimum.  We write
$[n]:=\{1,\ldots,n\}$.  For $\mathbb K\in\{\mathbb R,\mathbb C\}$, let
$\mathcal M_p(\mathbb K):=\mathbb K^{p\times p}$ and abbreviate
$\mathcal M_p:=\mathcal M_p(\mathbb C)$. For $S\subset\mathbb C$, set
\[
    \mathcal D_n(S)
    :=\{\diag(d_1,\ldots,d_n):d_i\in S\}.
\]
Here $\ii=\sqrt{-1}$, $M^\top$ is the transpose, and
$M^*=\overline M^\top$ is the conjugate transpose.
For Hermitian matrices $A,B\in\mathcal M_n$, the Loewner order
$A\preceq B$ means $x^*(B-A)x\geq0$ for every $x\in\mathbb C^n$;
we write $B\succeq A$ for the reverse order. For any square matrix $A$, set
\[
    \Re A:=\frac{A+A^*}{2},
    \qquad
    \Im A:=\frac{A-A^*}{2\ii}.
\]
For $x\in\mathbb C^n$, $D=\diag(x)$ has entries $D_i=x_i$.  The Euclidean,
operator, Hilbert--Schmidt, and nuclear norms are denoted by
$|\cdot|$, $|\cdot|_{\op}$, $|\cdot|_{\hs}$, and $|\cdot|_*$, respectively.
Thus
\begin{align*}
    |M|_{\hs}^2&=\tr(MM^*),
    &
    |M|_{\hs}&=\sup_{|A|_{\hs}\leq1}|\tr(A^*M)|,
    &
    |M|_*&=\sup_{|A|_{\op}\leq1}|\tr(A^*M)|.
\end{align*}
For $Z\in\mathcal M_p$, its numerical range is
\[
    \mathcal W(Z):=\{u^*Zu:u\in\mathbb C^p,\ |u|=1\}.
\]
It is compact and convex by the Toeplitz--Hausdorff theorem.
For $1\leq q<\infty$, $L^q$ denotes the space of real or complex random
variables $\xi$ satisfying $\mathbb E\left[|\xi|^q\right]<\infty$, equipped with the norm
\[
    \|\xi\|_{L^q}:=\bigl(\mathbb E\left[|\xi|^q\right]\bigr)^{1/q}.
\]
For a finite matrix profile $\Sigma=(\Sigma_i)_{i\in[n]}$, set
\begin{align}\label{eq:def_profile_hs_infty}
    |\Sigma|_\nu^\infty:=\max_{i\in[n]}|\Sigma_i|_\nu,
    \qquad \nu\in\{\op,\hs,*\}.
\end{align}
For vectors and diagonal matrices,
$|u|_\infty:=\max_i|u_i|$ and
$|D|_\infty:=\max_i|D_i|=|D|_{\op}$.  An index of the specialized notation appears in Appendix~\ref{sec:specific_notations_appendix}.

\tableofcontents

\section{Main results and examples}
\label{sec:main_results}

\subsection{Deterministic results}
\label{sse:main_deterministic_results}
We first define the deterministic equivalent and record the properties
needed for the probabilistic comparisons: existence, uniqueness,
holomorphic dependence, and Stieltjes-transform representations.
Columnwise self-consistent equations for heterogeneous correlation
matrices appear in \cite[Theorem~1]{WagnerCouilletDebbahSlock2012} and in
\cite{kammoun2016no,Yin2020DifferentCorrelations}. These works formulate
deterministic equivalents in asymptotic random-matrix settings;
\cite{WagnerCouilletDebbahSlock2012} also allows a positive-semidefinite
additive deformation. Their equations are already defined at each matrix
dimension. We use a formulation valid for every
finite $p,n$ and every real symmetric positive-semidefinite profile
$\Sigma_1,\ldots,\Sigma_n$, without an aspect-ratio or profile-size
assumption. The construction is also related
to the theory of self-consistent equations developed in
\cite{AjankiErdosKruger2019}.

\begin{theorem}\label{the:fixed_point_existence_uniqueness}
For every $z\in\mathbb C\setminus\mathbb R_+$, there are unique scalars
$\tilde\Gamma_1^z,\ldots,\tilde\Gamma_n^z\in\{a-b/z:a,b\geq0\}$
and a matrix $\tilde G^z$ satisfying
\begin{align*}
    \tilde G^z
    =-\frac1z\left(
        I_p+\frac1n\sum_{j=1}^n\tilde\Gamma_j^z\Sigma_j
      \right)^{-1},\qquad 
    \tilde\Gamma_i^z
    =-\frac1{z\left(1+\frac1n\tr(\Sigma_i\tilde G^z)\right)},
    \qquad i\in[n].
\end{align*}
The mappings $z\mapsto\tilde\Gamma_i^z$ and $z\mapsto\tilde G^z$
are holomorphic on $\mathbb C\setminus\mathbb R_+$.
\end{theorem}

The deterministic construction and the resolvent comparison estimates extend to complex matrix parameters \(Z\) satisfying \(d(\mathcal W(Z),\mathbb R_+)>0\). Their matrix formulations, including the dependence of the constants on \(Z\), are given in 
Sections~\ref{sec:deterministic_study} and~\ref{sec:probabilistic_results}.
Set
\[
    \tilde g(z):=\frac1p\tr(\tilde G^z),\qquad z\in\mathbb C\setminus\mathbb R_+.
\]

\begin{theorem}\label{the:g_tilde_stieltjes}
  For every $i\in[n]$, the mapping
  $z\mapsto\tilde\Gamma_i^z$ is the Stieltjes transform of a probability
  measure of compact support $\tilde\mu_i$ on $\mathbb R_+$.  The mapping
  $z\mapsto\tilde g(z)$ is the Stieltjes transform of a probability measure
  $\tilde\mu$ on $\mathbb R_+$, with
  \begin{align}\label{eq:relation_mu_mu_i}
      \tilde\mu
      =\left(1-\frac np\right)\delta_0
       +\frac1p\sum_{i=1}^n\tilde\mu_i.    
  \end{align}
  In addition, $z\mapsto\tilde G^z$ is the Stieltjes transform of a positive
  matrix-valued measure on $\mathbb R_+$ with total mass $I_p$.
\end{theorem}
  The right-hand side of \eqref{eq:relation_mu_mu_i} is understood as a finite signed combination; the
  identity asserts that this combination is in fact the positive probability
  measure whose Stieltjes transform is $\tilde g$.

Section~\ref{sec:deterministic_study} establishes existence and
uniqueness, the scalar and matrix-valued Stieltjes-transform representations,
an explicit support bound, and the stability estimate used in the
probabilistic comparisons.

\subsection{Probabilistic results}
\label{sse:main_probabilistic_results}
Fix an index set $\Theta\subset\mathbb N^2$.  For every
$(p,n)\in\Theta$, let $X=(x_1,\ldots,x_n)\in \mathcal M_{p,n}(\mathbb R)$ have independent real columns
with finite second moments, and retain the notation $\Sigma_i$ and $G^z$ from the introduction.  The dependence on $(p,n)$ is suppressed.

Fix a compact set $K\subset\mathbb C\setminus\mathbb R_+$. Throughout this section, $O(\cdot)$
means uniformly for $(p,n)\in\Theta$ and $z\in K$. The implicit constants
may depend on $K$, on the fixed moment order, and on the constants appearing
in the assumptions, but not on the deterministic test matrices or on the
profile quantities displayed in the bounds. The
assumption $p=O(n)$ in particular means that there exists a constant $C>0$ such that $(p,n)\in \Theta \implies p\leq Cn$. No uniformity is
claimed as $K$ approaches $\mathbb R_+$.

\subsubsection{Quadratic-form moduli and directional estimates}
\label{sssec:main_quadratic_form_moduli} 
For $q\geq1$, define
\begin{align}\label{eq:def_M_hs_op_main}
    M_{\hs}^{(q)}
    :=
    \sup_{\genfrac{}{}{0pt}{2}{i\in[n]}{|B|_{\hs}\leq1}}
    \left\|x_i^\top Bx_i-\tr(\Sigma_iB)\right\|_{L^q},
    \qquad
    M_{\op}^{(q)}
    :=
    \sup_{\genfrac{}{}{0pt}{2}{i\in[n]}{|B|_{\op}\leq1}}
    \left\|x_i^\top Bx_i-\tr(\Sigma_iB)\right\|_{L^q}.
\end{align}
The suprema are over deterministic complex matrices $B\in\mathcal M_p$.
Both moduli may depend on $(p,n)$ and may vanish. The reduction to
real symmetric tests, together with its conditional version, is recorded
in Lemma~\ref{lem:complex_conditional_test_reduction}.

\paragraph{Size of the two moduli.}
The standing second-moment assumption on the columns ensures that
$M_{\hs}^{(1)}$ and $M_{\op}^{(1)}$ are finite, by
Lemma~\ref{lem:uncentered_quadratic_form_bound}; higher-order moduli
may be infinite. They are nondecreasing in the moment order:
for $1\leq r\leq q$, $M_{\hs}^{(r)}\leq M_{\hs}^{(q)}$ and $M_{\op}^{(r)}\leq M_{\op}^{(q)}$.

The inclusions between the two matrix unit balls give
\[
    M_{\hs}^{(q)}
    \leq M_{\op}^{(q)}
    \leq\sqrt p\,M_{\hs}^{(q)}.
\]
In particular, $M_{\hs}^{(q)}=o(\sqrt p)$ implies
$M_{\op}^{(q)}=o(p)$, but the two moduli need not have the same order.

Their dependence on the second-moment matrices $\Sigma_i$ is explicit
in familiar models. Suppose that $x_i=\Sigma_i^{1/2}y_i$, where
the coordinates of each $y_i$ are independent, centered, and
standardized, and let $K_4:=\sup_{i,j}\|(y_i)_j\|_{L^4}<\infty$. 
Lemma~\ref{lem:profile_quadratic_form_moduli} gives
\[
    M_{\hs}^{(2)}
    =O\!\left(K_4^2|\Sigma|_{\op}^\infty\right),
    \qquad
    M_{\op}^{(2)}
    =O\!\left(K_4^2|\Sigma|_{\hs}^\infty\right),
\]
with the corresponding first-moment bounds following by monotonicity.
For centered Gaussian columns, these spectral scales are attained
for every fixed $q\geq2$:
\[
    M_{\hs}^{(q)}\asymp_q|\Sigma|_{\op}^\infty,
    \qquad
    M_{\op}^{(q)}\asymp_q|\Sigma|_{\hs}^\infty,
\]
by Lemma~\ref{lem:gaussian_quadratic_form_moduli}. Thus the
Hilbert--Schmidt normalization reflects the largest eigenvalues of
the matrices $\Sigma_i$, whereas the operator normalization reflects
their Hilbert--Schmidt sizes.

A large Hilbert--Schmidt size can also force a large
operator-normalized modulus without coordinate independence.
Indeed, Lemma~\ref{lem:quadratic_form_spectral_spread} gives
\[
    |\Sigma|_{\hs}^\infty
    \leq |\Sigma|_{\op}^\infty+M_{\op}^{(2)}.
\]

\paragraph{Consequences for resolvent tests.}
The comparison theorems bound
$\|\tr(A(G^z-\tilde G^z))\|_{L^q}$ for deterministic $A$, measured
in Hilbert--Schmidt or operator norm. Their constants are uniform
in the test matrix; any supremum over tests is taken outside the
$L^q$ norm.
Table~\ref{tab:comparison_rates} records selected specializations of Theorems~\ref{the:operator_profile_comparison} and~\ref{the:noncommuting_limited_second_moment_profiles}. When both comparison theorems apply, one may take the smaller bound.

The ordering of the moduli does not order the resulting resolvent
bounds: the Hilbert--Schmidt norm comparison uses the smaller modulus
but the larger test norm $|A|_{\hs}$, and the two comparisons have
different assumptions and coefficients involving the collection
$(\Sigma_i)_{i\in[n]}$. For high-rank tests, the smaller value of $|A|_{\op}$
can favor the operator norm comparison. For directional tests, however,
\[
    |uu^\top|_{\op}=|uu^\top|_{\hs}=|u|^2,
    \qquad u\in\mathbb R^p,
\]
so a smaller Hilbert--Schmidt modulus can improve the estimate
without increasing the test norm.

For example, for centered Gaussian columns with $\Sigma_i=I_p$
and $p\asymp n$, the two moduli have orders $1$ and $\sqrt p$,
respectively. For every fixed $q\geq2$,
Theorem~\ref{the:noncommuting_limited_second_moment_profiles}
therefore gives
\[
    \left\|u^\top(G^z-\tilde G^z)u\right\|_{L^q}
    =O\!\left(\frac{|u|^2}{\sqrt n}\right),
\]
whereas the operator-test estimate in
Theorem~\ref{the:operator_profile_comparison} does not yield
a vanishing bound for a fixed unit direction in this regime.

\subsubsection{Operator-norm route}
\label{sssec:main_operator_route}
Let $\Sigma^{(1)},\ldots,\Sigma^{(k)}$, with $k\leq n$, be the distinct
second-moment matrices. For $a\in[k]$, set
\[
 I_a:=\{i:\Sigma_i=\Sigma^{(a)}\},\qquad
 n_a:=|I_a|,\qquad \pi_a:=\frac{n_a}{n}>0.
\]
The number and sizes of these classes may depend on the dimensions. Note that columns that follow different laws but admit the same second-moment matrix belong to the same class. 
Define
\begin{align}\label{eq:def_tau_profile}
 \tau_\Sigma:=\max_{a\in[k]}
 \frac{1+|\Sigma^{(a)}|_{\op}}{1+\pi_a|\Sigma^{(a)}|_{\op}}.
\end{align} 
The coefficient depends
only on these matrices and their proportions, independently of $z$ and the
test. Unlike $\sigma_\Sigma$ below, it does not measure approximation:
close but unequal matrices with large operator norm belong to different classes.

\begin{remark}[Operator size and class proportions]
\label{rem:tau_profile_specializations}
Since $0<\pi_a\leq1$, each ratio in \eqref{eq:def_tau_profile}
is at least \(1\) and at most $\min (1+|\Sigma^{(a)}|_{\op},\pi_a^{-1})$. Thus
\[
 1\leq\tau_\Sigma\leq
 \min\{1+|\Sigma|_{\op}^\infty,(\min_a\pi_a)^{-1}\}.
\]
A common profile has $\tau_\Sigma=1$, including the zero profile.
Bounded operator profiles have bounded $\tau_\Sigma$ without a balance
condition. Equally sized $k$ classes give $\tau_\Sigma\leq k$; more
generally, $(\min_a\pi_a)^{-1}=O(k)$ gives $\tau_\Sigma=O(k)$, with
$k$ allowed to grow. 
\end{remark}

\begin{theorem}[Operator norm comparison with explicit trace dependence]
\label{the:operator_profile_comparison}
For every deterministic $A\in\mathcal M_p$ and every $q\geq1$,
\begin{align}
 &\|\tr(A(G^z-\tilde G^z))\|_{L^q}
 =O\left(|A|_{\op}
 \left(1+\frac{|\Sigma|_*^\infty}n\right)^{3/2}
 \tau_\Sigma\bigl(1+M_{\op}^{(q)}\bigr)\right).
\end{align}
For every $q\geq2$, one also has
\begin{align}
 &\|\tr(A(G^z-\tilde G^z))\|_{L^q}
 =O\left(|A|_{\op}
 \left(1+\frac{|\Sigma|_*^\infty}n\right)^{3/2}\left[
 \frac{\tau_\Sigma M_{\op}^{(q)}}{\sqrt n}
 +\frac{\tau_\Sigma^2+(M_{\op}^{(2)})^2}{n}
 \right]\right).
\end{align}
If $K\subset(-\infty,0)$, the factor
$\left(1+\frac{|\Sigma|_*^\infty}{n}\right)^{3/2}$ in both bounds may be
replaced by $1+\frac{|\Sigma|_*^\infty}{n}$.
\end{theorem}

\subsubsection{Hilbert--Schmidt-norm route}
\label{sssec:main_hs_route}
Define
\begin{align}\label{eq:def_sigma_profile}
    \sigma_\Sigma
    &:={}
    \inf_{\substack{
        T_1,\ldots,T_n\in\mathcal M_p(\mathbb R)\\
        T_i=T_i^\top\succeq0\ (i\in[n])\\
        T_iT_j=T_jT_i\ (i,j\in[n])
    }}
    \left(
        1+\dim_{\mathbb R}
        \operatorname{span}\{T_1,\ldots,T_n\}
    \right)
    \left(
        1+\frac1n\sum_{i=1}^n|\Sigma_i-T_i|_{\op}
    \right).
\end{align}
The coefficient combines the dimension of the span of a commuting approximating family 
with its average operator-norm error.  Any such family gives an
upper bound.
The matrices $T_i$ enter only the error coefficient; the
deterministic equivalent $\tilde G^z$ is always computed
from the original profile $(\Sigma_i)_{i\in[n]}$.

\begin{remark}[Simple bounds for the profile coefficient]
\label{rem:profile_interaction_simplifications}
The zero family gives
\[
    1\leq\sigma_\Sigma
    \leq1+\frac1n\sum_{i=1}^n|\Sigma_i|_{\op}
    \leq1+|\Sigma|_{\op}^\infty.
\]
Hence bounded operator profiles have $\sigma_\Sigma=O(1)$.  If the
$\Sigma_i$ commute pairwise, choosing $T_i=\Sigma_i$ gives
\[
    \sigma_\Sigma
    \leq
    1+\dim_{\mathbb R}
    \operatorname{span}\{\Sigma_1,\ldots,\Sigma_n\}.
\]
Thus a profile with $k\geq1$ distinct commuting values satisfies
$\sigma_\Sigma=O(k)$.  Every common profile satisfies
$\sigma_\Sigma\leq2$, while the zero profile has $\sigma_\Sigma=1$.
These bounds require no class-balance condition and no bound on the norms of
the commuting matrices.  Examples~\ref{exa:profile_commuting_low_dimension}--
\ref{exa:profile_commuting_matrix_perturbations} below illustrate these
bounds and the
additional flexibility obtained by using commuting approximants in
\eqref{eq:def_sigma_profile}.
\end{remark}

\begin{theorem}[Hilbert--Schmidt norm comparison under a profile bound]
\label{the:noncommuting_limited_second_moment_profiles}
Assume that
\begin{align}\label{eq:hs_profile_size_condition}
    \sqrt p\,|\Sigma|_{\hs}^\infty=O(n).
\end{align}
For every deterministic $A\in\mathcal M_p$ and every $q\geq1$,
\begin{align}\label{eq:noncommuting_limited_second_moment_uniform}
    \left\|
        \tr\left(A(G^z-\tilde G^z)\right)
    \right\|_{L^q}
    &={}
    O\left(
        \sigma_\Sigma|A|_{\hs}\left[
            M_{\hs}^{(q)}+\frac{|\Sigma|_{\hs}^\infty}n
        \right]
    \right).
\end{align}
For every $q\geq2$,
\begin{align}\label{eq:noncommuting_limited_second_moment_refined}
    \left\|
        \tr\left(A(G^z-\tilde G^z)\right)
    \right\|_{L^q}
    &={}
    O\left(
        |A|_{\hs}
        \left[
            \frac{M_{\hs}^{(q)}}{\sqrt n}
            +\sigma_\Sigma
             \frac{|\Sigma|_{\hs}^\infty+(M_{\hs}^{(2)})^2\sqrt p}{n}
        \right]
    \right).
\end{align}
\end{theorem}


\subsubsection{Resolvent sandwiches}
\label{sssec:main_resolvent_sandwiches}
For $z\in\mathbb C\setminus\mathbb R_+$ and a deterministic matrix $A\in\mathcal M_p$, set
\[
    \mathcal S^z[A]:=G^zAG^z,
\]
and define $\widetilde{\mathcal S}^{\,z}[A]$ by
\begin{align}
\label{eq:deterministic_linear_response_equation}
    \widetilde{\mathcal S}^{\,z}[A]
    ={}&
    \tilde G^z A\tilde G^z
    +
    \frac1{n^2}\sum_{i=1}^n
    \frac{
        \tr\!\left(
            \Sigma_i\widetilde{\mathcal S}^{\,z}[A]
        \right)
    }{
        \left(
            1+\frac1n\tr(\Sigma_i\tilde G^z)
        \right)^2
    }
    \tilde G^z\Sigma_i\tilde G^z.
\end{align}
Proposition~\ref{pro:deterministic_linear_response_equation} shows that
this equation has a unique solution, that
$A\mapsto\widetilde{\mathcal S}^{\,z}[A]$ is linear, and that, for the matrix-parameter equivalent defined in Section~\ref{sec:deterministic_study},
\begin{align*}
    \widetilde{\mathcal S}^{\,z}[A]
=-\left.\frac{d}{dt}\widetilde G^{\,zI_p-tA}\right|_{t=0}.
\end{align*}

\begin{corollary}[Random-to-deterministic resolvent sandwiches]
\label{cor:resolvent_sandwich_comparison}
Let $A,B\in\mathcal M_p$ be deterministic, and let $q\geq1$. In each of the
minima below, the second entry is available only when $q\geq2$.

\begin{enumerate}
\item[\textnormal{(i)}]
$\begin{aligned}[t]
&\left\|
    \tr\!\left(
        B\left[
            \mathcal S^z[A]
            -\widetilde{\mathcal S}^{\,z}[A]
        \right]
    \right)
\right\|_{L^q}
\nonumber\\
&\quad=
O\!\Biggl(
    |A|_{\op}|B|_{\op}
    \left(1+\frac{|\Sigma|_*^\infty}{n}\right)^{3/2}
    \min\left\{
        \tau_\Sigma\bigl(1+M_{\op}^{(q)}\bigr),
        \frac{\tau_\Sigma M_{\op}^{(q)}}{\sqrt n}
        +\frac{\tau_\Sigma^2+(M_{\op}^{(2)})^2}{n}
    \right\}
\Biggr).
\end{aligned}$

\item[\textnormal{(ii)}]
If $\sqrt p\,|\Sigma|_{\hs}^\infty=O(n)$,
then
\begin{align*}
&\left\|
    \tr\!\left(
        B\left[
            \mathcal S^z[A]
            -\widetilde{\mathcal S}^{\,z}[A]
        \right]
    \right)
\right\|_{L^q}
\nonumber\\
&\quad=
O\!\Biggl(
    |A|_{\op}|B|_{\hs}
    \min\left\{
        \sigma_\Sigma
        \left(
            M_{\hs}^{(q)}
            +\frac{|\Sigma|_{\hs}^\infty}{n}
        \right),
        \frac{M_{\hs}^{(q)}}{\sqrt n}
        +\sigma_\Sigma
        \frac{
            |\Sigma|_{\hs}^\infty
            +(M_{\hs}^{(2)})^2\sqrt p
        }{n}
    \right\}
\Biggr).
\end{align*}
\end{enumerate}
\end{corollary}

The proof, given in
Subsection~\ref{sec:linear_response_superoperator}, follows directly from
Theorems~\ref{the:operator_profile_comparison}
and~\ref{the:noncommuting_limited_second_moment_profiles}. The same argument
extends, for each fixed $k$, to symmetrized sandwiches with $k+1$ resolvent
factors and deterministic insertions $A_1,\ldots,A_k$ by taking successive
mixed derivatives of $\widetilde G^{\,zI_p-t_1A_1-\cdots-t_kA_k}$
at $t_1=\cdots=t_k=0$.

\paragraph{Ridge regression.}
Let $y=X^\top\beta+\varepsilon$, where $\beta$ is deterministic,
$\mathbb E[\varepsilon\mid X]=0$, and
$\mathbb E[\varepsilon\varepsilon^\top\mid X]=\sigma_\varepsilon^2I_n$.
For fixed $\lambda>0$, the ridge estimator
$\hat\beta_{\lambda}=G^{-\lambda}Xy/n$ has conditional bias
$-\lambda G^{-\lambda}\beta$.  If a test vector $x_0$, independent of
$(X,\varepsilon)$, has second moment $\Sigma$, then
\[
\mathbb E\!\left[
    (x_0^\top\hat\beta_{\lambda}-x_0^\top\beta)^2\mid X
  \right]
    =
    \lambda^2\beta^\top G^{-\lambda}\Sigma G^{-\lambda}\beta
    +\frac{\sigma_\varepsilon^2}{n}
       \tr\!\left(\Sigma\left[G^{-\lambda}
                    -\lambda(G^{-\lambda})^2\right]\right).
\]
This follows from
$\hat\beta_{\lambda}-\beta=-\lambda G^{-\lambda}\beta
+G^{-\lambda}X\varepsilon/n$ and
$\frac1n G^{-\lambda}XX^\top G^{-\lambda}=G^{-\lambda}-\lambda(G^{-\lambda})^2$.
Assume the classical conditions $p=O(n)$ and $M_{\hs}^{(2)}=O(1)$, which give
\[
    M_{\op}^{(2)}
    \leq\sqrt p\,M_{\hs}^{(2)}
    =O(\sqrt p).
\]
Suppose also that $\sup_{i\in[n]}|\Sigma_i|_{\op}=O(1)$ and
$|\Sigma|_{\op}=O(1)$ (the former assumption, together with $p=O(n)$,
implies $\sqrt p\,|\Sigma|_{\hs}^{\infty}=O(n)$).
We use the two bounds of
Corollary~\ref{cor:resolvent_sandwich_comparison} with $q=2$:
part~\textnormal{(i)} controls the squared-resolvent term in the variance
with insertion $I_p$ and test $\Sigma$, keeping the test factor
$|\Sigma|_{\op}$ bounded even when $|\Sigma|_{\hs}$ is of order $\sqrt p$;
part~\textnormal{(ii)} controls the squared-bias term with insertion
$\Sigma$ and test $\beta\beta^\top$, using
$|\beta\beta^\top|_{\hs}=|\beta|^2$ and the bounded modulus
$M_{\hs}^{(2)}$.
The one-resolvent theorems control the remaining variance term.  No
commutation between $\Sigma$ and the training profiles is required.

\subsection{Rates, comparisons, and examples}

\subsubsection{Comparison with earlier results}

The first-moment bounds give quantitative sufficient conditions within the
quadratic-form approach of Bai--Zhou \cite{BaiZhou2008} and Yaskov
\cite{Yaskov2016ECP,YaskovUniversality2014}. In the isotropic regime
$\Sigma_i=I_p$ with $p\asymp n$,
Theorem~\ref{the:operator_profile_comparison} gives
\[
    \left\|\frac1p\tr(G^z-\tilde G^z)\right\|_{L^1}
    =O\!\left(\frac{1+M_{\op}^{(1)}}p\right),
\]
so $M_{\op}^{(1)}=o(p)$ yields convergence with an explicit error bound.
Yaskov's universality principle already permits independent,
non-identically distributed columns; our comparisons add explicit rates
and $L^q$ control of deterministic matrix tests.

For heterogeneous second-moment matrices, the constructions in
\cite{WagnerCouilletDebbahSlock2012,Yin2020DifferentCorrelations}
use linear transforms of independent coordinates; our assumptions instead
concern quadratic forms of the columns themselves. 
Under the size and moment assumptions of our comparison theorems, we do not require the Gaussian-column and fixed-class assumptions
of \cite{BenaychGeorgesCouillet2016}, the asymptotic simultaneous
diagonalization assumed in \cite{MeiWangYao2023}, or the average-isotropy
assumptions of \cite{DembczakLytova2022}. 
The coefficients \(\tau_\Sigma\) and \(\sigma_\Sigma\) describe how relations among the second-moment matrices enter the bounds.
In particular, both coefficients
remain bounded when $\max_i|\Sigma_i|_{\op}=O(1)$, even if all $n$
second-moment matrices differ.

For directional tests, Hachem--Loubaton--Najim--Vallet
\cite{HachemLoubatonNajimVallet2013} obtain quantitative bilinear-form
estimates in a noncentered independent-entry model with a separable
variance profile. In the real, centered version of that model with bounded
variance profile and $p\asymp n$,
Theorem~\ref{the:noncommuting_limited_second_moment_profiles} gives the
same $n^{-1/2}$ scale for the $L^2$ error using uniformly bounded fourth
entry moments, whereas their Theorem~1.1 assumes eighth moments for this
bound. More generally, our theorem uses $M_{\hs}^{(2)}$ directly and does
not require coordinate independence.

Unlike the proportional-dimensional setting of
\cite{HachemLoubatonNajimVallet2013}, the operator estimate is formulated
without a restriction on $p/n$; whether it gives convergence is
specified by its trace and quadratic-form factors. Existing results for
$p\gg n$ include the eigenvalue-statistic CLTs of Ding--Wang
\cite{DingWang2025}, for independent entries and a common diagonal covariance,
and the dimension-free ridge-regression equivalents of Cheng--Montanari
\cite{ChengMontanari2024} and Misiakiewicz--Saeed
\cite{MisiakiewiczSaeed2024}, under i.i.d. sampling and spectral and
concentration assumptions. In the common-covariance case at $z<0$, our
fixed-point system reduces to the same effective-regularization equation
as in the ridge-regression results.
Theorem~\ref{the:operator_profile_comparison}
instead allows different matrices $\Sigma_i$ and measures the error through
their traces and quadratic-form moments. Example~\ref{exa:bounded_hs_profiles}
takes $p=n^4$ and centered Gaussian columns with covariance eigenvalues
$j^{-3/4}$ and growing traces $4n+O(1)$, and yields
\[
    \|\tr(A(G^z-\tilde G^z))\|_{L^q}
    =O\!\left(\frac{|A|_{\op}}{\sqrt n}\right),
    \qquad q\geq2,
\]
uniformly over the choices of covariance eigenbases.
This is an additive fixed-test comparison, not a CLT or the relative
risk approximation, potentially at vanishing regularization, obtained
in the cited works.

Finally, our fixed-$K$ estimates do not provide spectral confinement
as in \cite{BenaychGeorgesCouillet2016,ZhuangZhangXuSong2026}, or local
laws as in \cite{KnowlesYin2017,AltErdosKruger2017,FanMaPaquetteWang2026}.
For a common covariance in proportional dimensions,
Ma--Misiakiewicz~\cite{MaMisiakiewicz2026} obtain anisotropic local laws
under quadratic-form concentration without an additional cumulant
assumption, but with stronger tail and regularity hypotheses than
those used in our fixed-moment comparisons.
Approaching-axis estimates under Lipschitz concentration are also given
in \cite{Chouard2022}.

\subsubsection{Selected illustrations}
For the normalized-trace estimates below, write $\lambda_1,\ldots,\lambda_p$
for the eigenvalues of $\frac1n XX^\top$ and define the empirical spectral
measure and its Stieltjes transform by
\[
    \mu:=\frac1p\sum_{j=1}^p\delta_{\lambda_j},
    \qquad
    g(z):=\frac1p\tr(G^z)
         =\int_{\mathbb R_+}\frac{\mu(d\lambda)}{\lambda-z},
    \qquad z\in\mathbb C\setminus\mathbb R_+.
\]

\paragraph{Operator-norm route.}
\begin{example}[Balanced strong signals with Gaussian noise]
\label{exa:balanced_strong_signal_noise}
Let $x_i=v_a+y_i$ for $i\in I_a$, where $v_a\in\mathbb R^p$ are
deterministic and the $y_i\sim\mathcal N(0,\Sigma)$ are independent.
Assume that the class proportions $\pi_a=|I_a|/n$ satisfy
$(\min_a\pi_a)^{-1}=O(k)$. Then
\[
    \Sigma^{(a)}=\Sigma+v_av_a^\top,
    \qquad
    \tau_\Sigma=O(k).
\]
Neither orthogonality of the signals nor commutation of the
matrices $\Sigma^{(a)}$ is required.

For real symmetric $B$,
\[
    x_i^\top Bx_i-\tr(\Sigma^{(a)}B)
    =y_i^\top By_i-\tr(\Sigma B)+2v_a^\top By_i.
\]
Gaussian moment bounds and
Lemma~\ref{lem:gaussian_quadratic_form_moduli} therefore give,
for every fixed $q\geq2$,
\[
    M_{\op}^{(q)}
    =O\!\left(
        |\Sigma|_{\hs}
        +|\Sigma|_{\op}^{1/2}\max_a|v_a|
    \right).
\]
Under the trace condition
$\tr\Sigma+\max_a|v_a|^2=O(n)$,
Theorem~\ref{the:operator_profile_comparison} yields, for every
deterministic $A\in\mathcal M_p$,
\begin{align*}
    &\bigl\|\tr\bigl(A(G^z-\tilde G^z)\bigr)\bigr\|_{L^q}
    \\
    &\qquad=O\!\left(
        |A|_{\op}\left[
            \frac{k}{\sqrt n}
            \left(
                |\Sigma|_{\hs}
                +|\Sigma|_{\op}^{1/2}\max_a|v_a|
            \right)
            +\frac{
                k^2+|\Sigma|_{\hs}^2
                +|\Sigma|_{\op}\max_a|v_a|^2
            }{n}
        \right]
    \right).
\end{align*}
For unit-variance noise $\Sigma=I_p$, suppose that
$k=O(1)$, $\max_a|v_a|=O(\sqrt p)$, and $p=O(n)$.
The trace condition holds and $M_{\op}^{(q)}=O(\sqrt p)$.
Since $p/n=O(1)$, the estimates simplify to
\begin{align*}
    \bigl\|\tr\bigl(A(G^z-\tilde G^z)\bigr)\bigr\|_{L^q}
    =O\!\left(|A|_{\op}\sqrt{\frac pn}\right),
    \qquad \text{and}\qquad 
    \|g(z)-\tilde g(z)\|_{L^q}
    =O\!\left(\frac1{\sqrt{np}}\right).
\end{align*}
These estimates allow full-rank noise and second-moment matrices
with $|\Sigma^{(a)}|_{\op}=1+|v_a|^2$ of order $p$.
\end{example}
\begin{example}[Growing trace beyond proportional dimensions]
\label{exa:bounded_hs_profiles}
Let $p=n^4$, and let the columns be independent centered Gaussian vectors
with covariance matrices
\[
    \Sigma_i
    =U_i\diag\bigl(j^{-3/4}:1\leq j\leq p\bigr)U_i^\top,
    \qquad i\in[n],
\]
where the $U_i$ are arbitrary deterministic real orthogonal matrices.
Then
\begin{align*}
    |\Sigma|_{\op}^\infty&=1,
    \qquad
    \bigl(|\Sigma|_{\hs}^\infty\bigr)^2
    =\sum_{j=1}^{n^4}j^{-3/2}\asymp1,\qquad
    |\Sigma|_*^\infty
    =\sum_{j=1}^{n^4}j^{-3/4}
      =4n+O(1).
\end{align*}
Thus $\tau_\Sigma\leq2$ and
$1+|\Sigma|_*^\infty/n=5+O(n^{-1})$.
Lemma~\ref{lem:gaussian_quadratic_form_moduli} also gives
$M_{\op}^{(q)}=O(1)$ for every fixed $q\geq2$.
Consequently, Theorem~\ref{the:operator_profile_comparison} yields
\begin{equation}
\label{eq:square_summable_profile_operator_rate}
    \left\|\tr\!\left(A(G^z-\tilde G^z)\right)\right\|_{L^q}
    =O\!\left(\frac{|A|_{\op}}{\sqrt n}\right),
    \qquad q\geq2,
\end{equation}
for every deterministic $A\in\mathcal M_p$, uniformly for $z$ in any
fixed compact set $K\subset\mathbb C\setminus\mathbb R_+$. The constants are independent of
the choices of $U_i$; no common eigenbasis or bound on the number of
distinct covariance matrices is required. In particular, taking $A=I_p$
gives an $O(n^{-1/2})$ error for the unnormalized trace, although $p/n=n^3$.
The Hilbert--Schmidt norm comparison does not apply, since
\[
    \frac{\sqrt p\,|\Sigma|_{\hs}^\infty}{n}\asymp n,
\]
so its size condition \eqref{eq:hs_profile_size_condition} fails.

The population resolvent $(I_p+\frac1n \sum_{i=1}^n\Sigma_i)^{-1}$ does not provide the same approximation, even
in the common-covariance case $U_i=I_p$ for every $i$.
Write $\Sigma=\diag(j^{-3/4}:1\leq j\leq p)$ for this common covariance.
At $z=-1$, write $\tilde G:=\left.\tilde G^z\right|_{z=-1}$.
The fixed-point equations in
Theorem~\ref{the:fixed_point_existence_uniqueness} can be written as
\begin{align*}
    \tilde G
    =(I_p+\tilde\Gamma\Sigma)^{-1},
    \qquad 0<\tilde\Gamma\leq1,
    \qquad 
    \frac1{\tilde\Gamma}
    =1+\frac1n\sum_{j=1}^{n^4}
       \frac{j^{-3/4}}{1+\tilde\Gamma j^{-3/4}}.
\end{align*}
Since
\[
    0\leq
    \sum_{j=1}^{n^4}j^{-3/4}
    -\sum_{j=1}^{n^4}\frac{j^{-3/4}}{1+\tilde\Gamma j^{-3/4}}
    =\tilde\Gamma\sum_{j=1}^{n^4}\frac{j^{-3/2}}{1+\tilde\Gamma j^{-3/4}}
    =O(1),
\]
we obtain $\tilde\Gamma^{-1}=5+O(n^{-1})$, and hence
$\tilde\Gamma=1/5+O(n^{-1})$.
Moreover, uniformly for $0\leq\Gamma\leq1$,
\[
    \sum_{j=1}^{n^4}
    \frac{j^{-3/4}}{(1+\Gamma j^{-3/4})(1+j^{-3/4})}
    =\sum_{j=1}^{n^4}j^{-3/4}
      +O\!\left(\sum_{j=1}^{n^4}j^{-3/2}\right)
    =4n+O(1).
\]
It follows that
\begin{align*}
    \tr\!\left(
        \tilde G-(I_p+\Sigma)^{-1}
    \right)
    &=(1-\tilde\Gamma)\sum_{j=1}^{n^4}
      \frac{j^{-3/4}}{(1+\tilde\Gamma j^{-3/4})(1+j^{-3/4})}
    =\frac{16}{5}n+O(1).
\end{align*}
Thus the fixed-point correction is of order $n$ in unnormalized trace,
whereas \eqref{eq:square_summable_profile_operator_rate} controls the
random-to-deterministic trace error at order $n^{-1/2}$.
\end{example}

\paragraph{Examples for the commuting-approximation coefficient.}
The next two examples bound $\sigma_\Sigma$; the profile-size and
moment assumptions of Theorem~\ref{the:noncommuting_limited_second_moment_profiles} must still be checked.

\begin{example}[A commuting family of low linear dimension]
\label{exa:profile_commuting_low_dimension}
The span-dimension bound in
Remark~\ref{rem:profile_interaction_simplifications} does not require a
bound on the number of distinct profiles.  For example, if
\[
    \Sigma_i=a_iA+b_iB,
    \qquad
    A,B\succeq0,\quad AB=BA,\quad a_i,b_i\geq0,
\]
then $\sigma_\Sigma\leq3$.  This remains true when all the profiles are
distinct or their norms grow. More generally, if the matrices $\Sigma_i$
are positive semidefinite, pairwise commuting, and span an $r$-dimensional
real vector space, then $\sigma_\Sigma\leq1+r$; no representation by
nonnegative coordinates is required.
\end{example}

\begin{example}[Perturbations of commuting matrices]
\label{exa:profile_commuting_matrix_perturbations}
Let $T^{(1)},\ldots,T^{(k)}$ be pairwise commuting real positive
semidefinite matrices, and suppose that
\[
    \Sigma_i=T^{(a(i))}+E_i\succeq0,
    \qquad a(i)\in[k].
\]
Choosing $T_i=T^{(a(i))}$ in \eqref{eq:def_sigma_profile} gives
\[
    \sigma_\Sigma
    \leq
    \left(
        1+\dim_{\mathbb R}
        \operatorname{span}\{T^{(1)},\ldots,T^{(k)}\}
    \right)
    \left(1+\frac1n\sum_{i=1}^n|E_i|_{\op}\right)
    =O\left(k\left(1+\frac1n\sum_{i=1}^n|E_i|_{\op}\right)\right).
\]
The matrices $E_i$ may be indefinite and need not commute. No bound on
the operator norms $|T^{(a)}|_{\op}$ or balance among the classes is
required.

For $k=1$, write $T=T^{(1)}$, so that $\Sigma_i=T+E_i$. The bound becomes
\[
    \sigma_\Sigma
    \leq2\left(1+\frac1n\sum_{i=1}^n|E_i|_{\op}\right).
\]
Choosing $T=\overline\Sigma:=\frac1n\sum_j\Sigma_j$ gives the explicit
average-deviation estimate
\[
    \sigma_\Sigma
    \leq2\left(1+\frac1n\sum_{i=1}^n
        |\Sigma_i-\overline\Sigma|_{\op}\right).
\]
If at most $m$ perturbations $E_i$ are nonzero and $|E_i|_{\op}=O(1)$, the
average-error factor in the general bound is $1+O(m/n)$.
\end{example}

\paragraph{Why second moments alone do not suffice.}
The following example shows that a common unconditional second moment does
not determine the limiting empirical spectral distribution when quadratic
forms fail to concentrate.  It also shows how conditioning on latent radii
recovers the appropriate heterogeneous-profile deterministic equivalent.
For a broader treatment of spectral limits for elliptical distributions,
see El Karoui~\cite{ELK09}.

\begin{example}[An isotropic radial mixture]
\label{exa:isotropic_radial_mixture}
Let $x_i=\sqrt{r_i}\,y_i$, where the $y_i$ are independent and uniform
on the sphere of radius $\sqrt p$, and the $r_i$ are independent of
the directions and of each other, with
$\mathbb P(r_i=1/5)=\mathbb P(r_i=9/5)=1/2$.
Since $\mathbb E[r_i]=1$ and $\mathbb E[y_i y_i^\top]=I_p$,
independence gives $\Sigma_i=\mathbb E[x_i x_i^\top]=I_p$.
However, $|x_i|^2=pr_i$ and $|r_i-1|=4/5$, so the tests
$B=I_p$ and $B=I_p/\sqrt p$ give
\[
    \frac{M_{\op}^{(q)}}p\geq\frac45,
    \qquad
    \frac{M_{\hs}^{(q)}}{\sqrt p}\geq\frac45,
    \qquad q\geq1.
\]

Suppose that $p,n\to\infty$ with $p/n\to\gamma\in(0,\infty)$,
and condition on the radii $ r=(r_1,\ldots,r_n)$.
The columns remain independent, with conditional second moments
$r_iI_p$. For every real symmetric $B$, the spherical variance bound
\[
    \Var(y_i^\top By_i)\leq2|B|_{\hs}^2
\]
and $r_i\leq9/5$ give a conditional bound $M_{\hs}^{(2)}=O(1)$,
uniformly in $ r$, including complex tests by
transpose symmetrization and real--imaginary decomposition.
Theorem~\ref{the:noncommuting_limited_second_moment_profiles}
therefore yields
\[
    \left\|g(z)-\tilde g_{ r}(z)\right\|_{
        L^2(\,\cdot\,\mid r)}
    =O(n^{-1}),
\]
uniformly in $ r$, where $\tilde g_{ r}$ is the
deterministic equivalent for the conditional profile. Since this
profile consists of scalar matrices,
\[
    \frac1{\tilde g_{ r}(z)}
    =-z+\frac1n\sum_{i=1}^n
    \frac{r_i}{1+(p/n)r_i\tilde g_{ r}(z)}.
\]

The proportions of the two radii converge in probability to $1/2$.
The uniform support bound gives subsequential weak compactness for
the conditional deterministic measures. Passing to the scalar equation
and using the sector uniqueness argument of
Proposition~\ref{pro:fixed_point_secant_contraction} then gives
$\tilde g_{ r}(z)\to m(z)$ in probability, where $m$ is
the unique Stieltjes transform of a probability measure on
$\mathbb R_+$ satisfying
\[
    \frac1{m(z)}
    =-z+\frac1{2(5+\gamma m(z))}
       +\frac9{2(5+9\gamma m(z))}.
\]
The conditional comparison therefore gives weak convergence in
probability of the empirical spectral distribution to this law.

Expansion at infinity shows that this law has first moment one and second moment
\[
    1+\gamma\,\mathbb E[r_i^2]
    =1+\frac{41}{25}\gamma.
\]
This differs from the second moment $1+\gamma$ of the Marchenko--Pastur
law predicted by the unconditional profile $I_p$.
The discrepancy is illustrated in
Figure~\ref{fig:radial_mixture_two_profile}.
\end{example}
\subsubsection{Independent entries}
For independent entries, the first-moment Hilbert--Schmidt norm
comparison gives qualitative convergence when the squared
entries are uniformly integrable, including every fixed
finite-variance entry law. Higher moments yield quantitative
bounds.
For $c\in[0,\infty)$, write $\mu_{\mathrm{MP},c}$ for the
Marchenko--Pastur law with parameter $c$, with
$\mu_{\mathrm{MP},0}:=\delta_1$.

\begin{corollary}[Marchenko--Pastur convergence under finite variance]
\label{cor:mp_iid_finite_variance}
For every $(p,n)\in\Theta$, suppose that the entries of $X$ are
independent copies of a centered unit-variance real random variable
$\xi_{p,n}$. Assume that the family $(\xi_{p,n}^2)_{(p,n)\in\Theta}$
is uniformly integrable and that $p=O(n)$.
Then, along every sequence in $\Theta$ with
$p\to\infty$, for every compact $K\subset\mathbb C\setminus\mathbb R_+$,
\[
    \sup_{z\in K}\|g(z)-\tilde g(z)\|_{L^1}\longrightarrow0.
\]
If, along the same sequence, $p/n\to c\in[0,\infty)$, the empirical
spectral distribution of $\frac1n XX^\top$ converges weakly in probability
to $\mu_{\mathrm{MP},c}$.
\end{corollary}

Its proof appears immediately after the proof of
Theorem~\ref{the:noncommuting_limited_second_moment_profiles}.

\begin{remark}[Rates under higher moments]
\label{rem:iid_finite_moment_rates}
Suppose that $p=O(n)$ and that, within each matrix, the entries are
centered and i.i.d., with a law that may depend on $(p,n)$.
For any fixed $r\in[1, 2]$ such that $\|X_{11}\|_{L^{2r}}=O(1)$, for every deterministic
$A\in\mathcal M_p$,
\begin{align*}
    \left\|\tr\left(A(G^z-\tilde G^z)\right)\right\|_{L^1}
    &=O\left(|A|_{\hs}p^{\frac{1}{r}-\frac{1}{2}}\right),
    &
    \|g(z)-\tilde g(z)\|_{L^1}
    &=O\left(p^{\frac{1}{r}-1}\right).
\end{align*}
Indeed, since
$\mathbb E[X_{11}^2]=O(1)$ and
$\Sigma_i=\mathbb E[X_{11}^2]I_p$, 
$|\Sigma|_{\hs}^\infty=O(\sqrt p)$.
The common scalar profile satisfies $\sigma_\Sigma\leq2$, and
Proposition~\ref{pro:quadratic_form_independent_coordinates_L1} gives
$M_{\hs}^{(1)}=O(p^{\frac{1}{r}-\frac{1}{2}})$. Since
$|\Sigma|_{\hs}^\infty/n=O(\sqrt p/n)=O(1)$,
the first-moment estimate in
Theorem~\ref{the:noncommuting_limited_second_moment_profiles} gives the
first bound, and $|I_p/p|_{\hs}=\frac1{\sqrt p}$ gives the second.

At \(r=1\), the displayed bound does not imply convergence. Uniform integrability, as assumed in Corollary~\ref{cor:mp_iid_finite_variance}, provides a sufficient additional condition.

If $\|X_{11}\|_{L^4}=O(1)$, the variance calculation in
Lemma~\ref{lem:profile_quadratic_form_moduli} gives directly
\[
    M_{\hs}^{(2)}=O(\|X_{11}\|_{L^4}^2)=O(1).
\]
The same comparison theorem with $q=2$
then yields
\[
    \left\|\tr\left(A(G^z-\tilde G^z)\right)\right\|_{L^1}
    \leq \left\|\tr\left(A(G^z-\tilde G^z)\right)\right\|_{L^2}
    =O\left(
        |A|_{\hs}\left[\frac1{\sqrt n}+\frac{\sqrt p}{n}\right]
      \right)
    =O\left(\frac{|A|_{\hs}}{\sqrt n}\right),
\]
and consequently
\[
    \|g(z)-\tilde g(z)\|_{L^1}
    \leq\|g(z)-\tilde g(z)\|_{L^2}
    =O\left(\frac1{\sqrt{np}}\right).
\]
\end{remark}

\clearpage

\section{Deterministic equivalents for sectorial matrix parameters}
\label{sec:deterministic_study}

We prove the deterministic results for general complex matrix parameters.
Throughout the proofs, we write
$\Omega:=\mathbb C\setminus\mathbb R_+$, recalling that $\mathbb R_+=[0,\infty)$.
For $Z\in\mathcal M_p$, let
$\mathcal W(Z):=\{u^*Zu:u\in\mathbb C^p,\ |u|=1\}$ be its numerical
range, and set
\[
    \eta_Z:=d\bigl(\mathcal W(Z),\mathbb R_+\bigr),
    \qquad
    \eta_z:=\eta_{zI_p}=d(z,\mathbb R_+),\qquad z\in\Omega.
\]
We call $Z$ admissible when $\eta_Z>0$, and associate with it the closed
convex sector
\begin{align}\label{eq:def_matrix_resolvent_sector}
    \mathcal C_Z:=\cone\bigl(\mathbb R_+\cup-\mathcal W(Z)\bigr).
\end{align}
The sector $\mathcal C_Z$ specifies the admissible values of the
fixed-point coefficients. It allows a common treatment of nonreal
scalar parameters and negative real parameters. General matrix parameters are needed in the classwise replacement
argument of the operator norm route, and in the additive
deformations used to obtain resolvent sandwiches; see
Subsections~\ref{sse:operator_profile_route}
and~\ref{sec:linear_response_superoperator}.
The supporting phase introduced below gives the real-part bounds
used to control the resolvents and the fixed-point denominators.

\subsection{Matrix sectors and two-sided resolvent calculus}

\begin{figure}[htbp]
    \centering
    \includegraphics[width=0.485\linewidth,
        trim=43bp 33.4bp 43bp 33.4bp,clip]{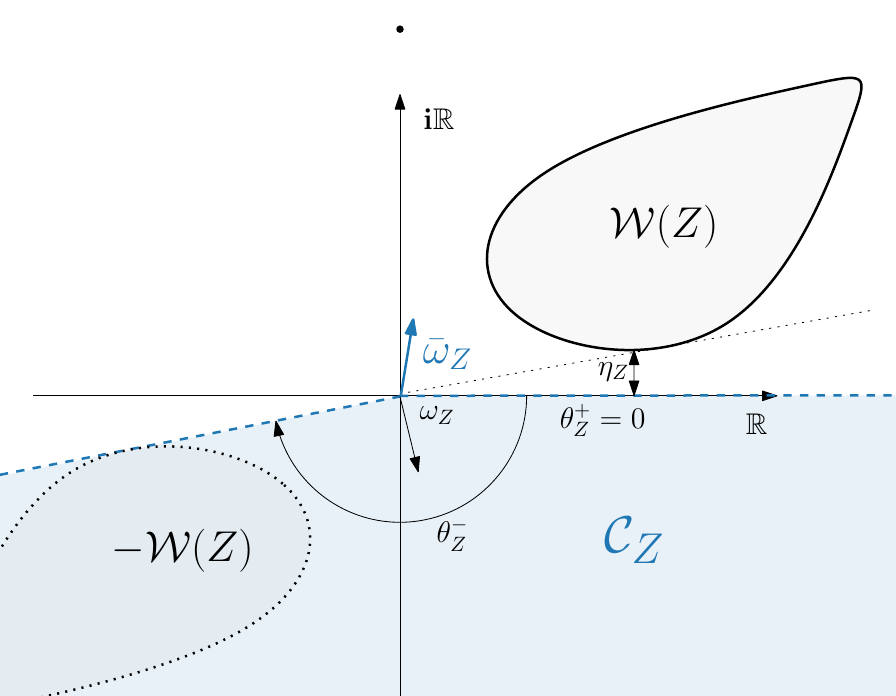}\hfill
    \includegraphics[width=0.485\linewidth,
        trim=43bp 33.4bp 43bp 33.4bp,clip]{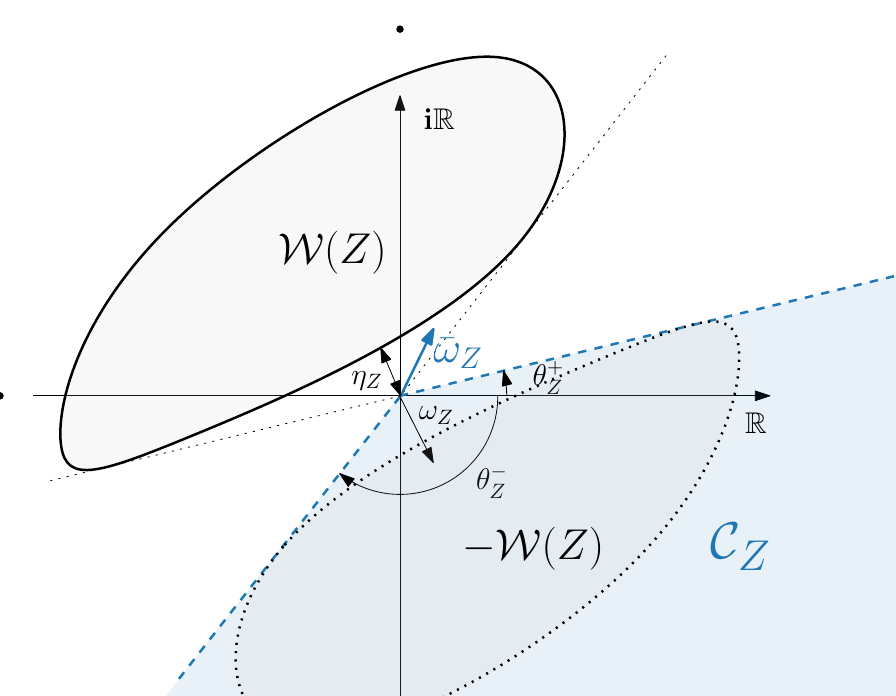}
    \caption{Geometry of the matrix sector $\mathcal C_Z$ for two numerical-range configurations.}
    \label{fig:matrix_parameter_sector_geometry}
\end{figure}

For an admissible matrix $Z$, $\mathcal C_Z$ is convex and contains
$\mathbb R_+$. Its closedness and aperture bound are proved below.
Figure~\ref{fig:matrix_parameter_sector_geometry} illustrates the geometry
in the complex plane.
For an admissible parameter \(Z\), define
\[
    \theta_Z^+
    :=\sup\left\{
        \theta\in(-\pi,\pi):e^{\ii\theta}\in\mathcal C_Z
    \right\},
    \qquad
    \theta_Z^-
    :=\inf\left\{
        \theta\in(-\pi,\pi):e^{\ii\theta}\in\mathcal C_Z
    \right\}.
\]

\begin{lemma}[Aperture of an admissible matrix sector]
\label{lem:matrix_parameter_sector_aperture}
For every admissible matrix \(Z\), the cone \(\mathcal C_Z\) is closed,
and the principal arguments of its nonzero elements form the closed interval
\([\theta_Z^-,\theta_Z^+]\), with \(\theta_Z^-\leq0\leq\theta_Z^+\) and
\begin{align}\label{eq:matrix_parameter_sector_aperture}
    0\leq\theta_Z^+-\theta_Z^-
    \leq2\arccos\left(\frac{\eta_Z}{2|Z|_{\op}}\right)<\pi.
\end{align}
\end{lemma}

\begin{proof}
The set \(B:=\operatorname{conv}(\{1\}\cup-\mathcal W(Z))\) is compact
and convex and does not contain zero. Indeed, an identity
\(0=t-(1-t)w\), with \(t\in[0,1]\) and \(w\in\mathcal W(Z)\), would
force \(t<1\) and \(w=t/(1-t)\in\mathbb R_+\), contrary to
admissibility. Choose \(q\in B\) of minimum modulus.
Every element of \(\mathcal C_Z\) is a nonnegative multiple of a point
of \(B\). If \(r_mv_m\) converges, with \(r_m\geq0\) and \(v_m\in B\),
then \(r_m\leq|r_mv_m|/|q|\) is bounded. Compactness of \(B\) allows
us to pass to a subsequence with \(r_m\to r\geq0\) and \(v_m\to v\in B\),
so the limit \(rv\) belongs to \(\mathcal C_Z\). This proves closedness.

For \(v\in B\) and \(0<t\leq1\), convexity and minimality give
\[
    0\leq\frac{|q+t(v-q)|^2-|q|^2}{t}
    =2\Re\bigl(\overline q(v-q)\bigr)+t|v-q|^2.
\]
Letting \(t\downarrow0\) yields
\(\Re(\overline qv)\geq|q|^2>0\). With
\(M:=\max_{v\in B}|v|\), it follows that
\[
    \Re\left(\frac{\overline q}{|q|}v\right)
    \geq\frac{|q|}{M}|v|,\qquad v\in B.
\]
By addition and homogeneity, the same inequality holds on
\(\mathcal C_Z\). Thus this cone has aperture strictly smaller than
\(\pi\). Since it contains \(1\), it excludes the negative real axis;
its principal arguments therefore form a closed interval
\([\theta_Z^-,\theta_Z^+]\) containing zero.

Put \(\delta:=\theta_Z^+-\theta_Z^-\). Compactness of \(\mathcal W(Z)\)
and \(|w|\geq\eta_Z\) for \(w\in\mathcal W(Z)\) ensure that each
boundary ray is either \(\mathbb R_+\) or contains a point of
\(-\mathcal W(Z)\). If both boundary rays contain points
\(-w_\pm=r_\pm e^{\ii\theta_Z^\pm}\), then convexity gives
\[
    w_0:=\frac{r_-w_++r_+w_-}{r_++r_-}\in\mathcal W(Z),
    \qquad
    \eta_Z\leq|w_0|
    =\frac{2r_+r_-}{r_++r_-}\cos(\delta/2)
    \leq|Z|_{\op}\cos(\delta/2).
\]
If one boundary ray is $\mathbb R_+$, the other contains a point $-w=re^{\pm\ii\delta}$ for $w\in\mathcal W(Z)$.
Taking the point $r\in\mathbb R_+$ as a competitor in the definition of
$\eta_Z$, we obtain
\[
    \eta_Z
    \leq |w-r|
    =
    \left|-re^{\pm\ii\delta}-r\right|
    =
    r\left|1+e^{\pm\ii\delta}\right|
    =
    2r\cos(\delta/2)
    \leq
    2|Z|_{\op}\cos(\delta/2).
\]
For \(\delta=0\), the same conclusion follows from
\(\eta_Z\leq|Z|_{\op}\). These bounds prove
\eqref{eq:matrix_parameter_sector_aperture}.
\end{proof}

Define the bisecting phase and its angular constant by
\begin{align}\label{eq:matrix_parameter_angular_phase}
    \theta_Z:=\frac{\theta_Z^++\theta_Z^-}{2},
    \qquad
    \omega_Z:=e^{\ii\theta_Z},
    \qquad
    c_Z:=\cos\left(\frac{\theta_Z^+-\theta_Z^-}{2}\right)>0.
\end{align}

\begin{lemma}[Supporting phase for an admissible parameter]
\label{lem:matrix_parameter_supporting_phase}
For every admissible matrix \(Z\), the phase and constant in
\eqref{eq:matrix_parameter_angular_phase} satisfy
\begin{align}
    \Re(\overline{\omega_Z}u)\geq c_Z|u|,
    \qquad u\in\mathcal C_Z,
\end{align}
and
\begin{align}\label{eq:matrix_parameter_coercivity}
    \Re(\overline{\omega_Z}(-Z))
    \succeq c_Z\eta_Z I_p\succ0.
\end{align}
\end{lemma}

\begin{proof}
For \(u\neq0\) in \(\mathcal C_Z\), write \(u=|u|e^{\ii\theta}\).
Lemma~\ref{lem:matrix_parameter_sector_aperture} gives
\[
    |\theta-\theta_Z|\leq\frac{\theta_Z^+-\theta_Z^-}{2},
\]
so
\(\Re(\overline{\omega_Z}u)=|u|\cos(\theta-\theta_Z)\geq c_Z|u|\).
The assertion also holds for \(u=0\).
For every unit vector \(v\), one has
\(-v^*Zv\in\mathcal C_Z\) and \(|v^*Zv|\geq\eta_Z\). Hence
\[
    v^*\Re(\overline{\omega_Z}(-Z))v
    =\Re\bigl(\overline{\omega_Z}(-v^*Zv)\bigr)
    \geq c_Z|v^*Zv|\geq c_Z\eta_Z,
\]
which proves \eqref{eq:matrix_parameter_coercivity}.
\end{proof}

We write
\[
    \overline{\mathcal C_Z}:=\{\bar w:w\in\mathcal C_Z\}
\]
for the conjugate cone. The next lemma applies to any supporting phase
satisfying its hypotheses.

\begin{lemma}[Sectorial resolvent calculus]
\label{lem:distance_resolvent_and_ward}
Let \(Z\in\mathcal M_p\), let \(\mathcal C_Z\) be given by
\eqref{eq:def_matrix_resolvent_sector}, and suppose that a unit phase
\(\omega\) and constants \(c,\kappa>0\) satisfy
\begin{equation}\label{eq:sectorial_resolvent_hypotheses}
    \forall v\in \mathcal C_Z:\quad \Re(\overline{\omega}v)\geq c|v|
    \qquad
    \text{and}
    \qquad
    \Re(\overline{\omega}(-Z))\succeq\kappa I_p.
\end{equation}
Given \(D\in\mathcal D_n(\mathcal C_Z)\) and positive semidefinite
matrices \(S_1,\ldots,S_n\in\mathcal M_p\), set
\[
    S:=\frac1n\sum_{j=1}^nD_jS_j
    \qquad 
    \text{and}
    \qquad 
    H=(S-Z)^{-1}.
\]
Then \(H\) is well defined and satisfies
\[
    |H|_{\op}\leq\frac1\kappa.
\]
Moreover, $\Re(\omega H)=H^*\Re\!\left(\overline{\omega}(S-Z)\right)H
      =H\Re\!\left(\overline{\omega}(S-Z)\right)H^*$, and the following
\emph{Ward inequalities} hold:
\begin{align*}
    H^*H,\ HH^*
    &\preceq\frac1\kappa\Re(\omega H).
\end{align*}
For \(x\in\mathbb C^p\) and \(\Sigma\succeq0\),
\[
    x^*Hx\in\overline{\mathcal C_Z},
    \qquad
    \tr(\Sigma H)\in\overline{\mathcal C_Z}.
\]
Finally, for every \(w\in\overline{\mathcal C_Z}\),
\[
    |1+w|\geq c(1+|w|),
    \qquad (1+w)^{-1}\in\mathcal C_Z.
\]
\end{lemma}

\begin{proof}
For every unit vector \(u\), the definition of \(\mathcal C_Z\) gives
\(-u^*Zu\in\mathcal C_Z\). Since \(u^*S_ju\geq0\) and
\(D_j\in\mathcal C_Z\), one also has
\[
    u^*Su=\frac1n\sum_jD_j(u^*S_ju)\in\mathcal C_Z,
    \qquad u^*(S-Z)u\in\mathcal C_Z.
\]
These inclusions extend to every \(u\) by homogeneity.
Since \(1\in\mathcal C_Z\), the supporting phase satisfies
\(\Re(\omega)\geq c\). In particular,
\(\Re(\overline{\omega}(s-w))\geq cs+\kappa\) for
\(s\geq0\) and \(w\in\mathcal W(Z)\), so these hypotheses imply
\(\eta_Z\geq\kappa>0\). The supporting-phase and coercivity bounds give
\[
    \Re\!\left(\overline{\omega}(S-Z)\right)
    \succeq\frac cn\sum_j|D_j|S_j+\kappa I_p
    \succeq\kappa I_p.
\]
Thus \(|(S-Z)u|\geq\kappa|u|\), proving invertibility and the inverse
norm bound. Multiplying the identity \(H(S-Z)=(S-Z)H=I_p\) and its
adjoint in the two orders gives the displayed inverse identities.
Coercivity then yields both Ward inequalities.

For the sector inclusions, set \(y=Hx\). Then
\[
    x^*Hx=((S-Z)y)^*y
    =\overline{y^*(S-Z)y}\in\overline{\mathcal C_Z}.
\]
Diagonalizing \(\Sigma\succeq0\) and summing its eigenvalue-weighted
quadratic forms proves the trace assertion.



Since \(1\in\mathcal C_Z\), the supporting phase gives
\(\Re(\omega)\geq c\). For \(w\in\overline{\mathcal C_Z}\), it also
gives \(\Re(\omega w)\geq c|w|\), so
\[
    |1+w|\geq\Re(\omega(1+w))\geq c(1+|w|).
\]
Moreover,
\[
    (1+w)^{-1}=\frac{\overline{1+w}}{|1+w|^2}\in\mathcal C_Z.
\]
\end{proof}

Expectations of integrable random variables taking values in either
sector, and line segments joining points in either sector, remain in that
sector. All assertions of the lemma remain valid when $\mathcal C_Z$
is replaced throughout by a larger closed convex sector containing
$\mathcal C_Z$ and satisfying the same supporting-phase bound. The expectation and line-segment statements
follow from closedness and convexity. The extension to larger sectors
follows by the same proof, retaining the coercivity bound in
\eqref{eq:sectorial_resolvent_hypotheses} and allowing the coefficients
\(D_j\) in the larger sector.

For every admissible \(Z\), Lemma~\ref{lem:matrix_parameter_supporting_phase}
lets us apply the lemma with \(\omega=\omega_Z\), \(c=c_Z\), and
\(\kappa=c_Z\eta_Z\).

\subsection{The matrix-parameter fixed point}
Let $\Sigma_1,\ldots,\Sigma_n\in\mathcal M_p(\mathbb R)$ be symmetric
positive-semidefinite matrices.
For \(D\in\mathcal D_n(\mathcal C_Z)\), define
\begin{align}\label{eq:def_matrix_parameter_resolvent}
    G_D^Z
    &:={}
    \left(
        \frac1n\sum_{j=1}^nD_j\Sigma_j-Z
    \right)^{-1},
    \qquad
    \Phi^Z(D)_i
    :={}
    \left(
        1+\frac1n\tr(\Sigma_iG_D^Z)
    \right)^{-1},
    \qquad i\in[n].
\end{align}

By Lemma~\ref{lem:distance_resolvent_and_ward},
\[
    \frac1n\tr(\Sigma_iG_D^Z)
    \in\overline{\mathcal C_Z}.
\]
The reciprocal bounds in the same lemma give
\[
    \Phi^Z(D)_i\in\mathcal C_Z,
    \qquad
    |\Phi^Z(D)_i|\leq c_Z^{-1}.
\]
Thus \(\Phi^Z\) maps the nonempty compact convex set
\[
    \left\{
        D\in\mathcal D_n(\mathcal C_Z):
        |D|_\infty\leq c_Z^{-1}
    \right\}
\]
continuously into itself.  Brouwer's fixed-point theorem therefore gives the following result.

\begin{lemma}[Sector invariance and existence]
\label{lem:fixed_point_compact_invariant_set}
For every admissible \(Z\), the mapping \(\Phi^Z\) has a fixed point
\(\widetilde\Delta^Z\in\mathcal D_n(\mathcal C_Z)\), and
\[
    |\widetilde\Delta^Z|_\infty\leq c_Z^{-1}.
\]
\end{lemma}

For a fixed point supplied by the lemma, define
\[
    \widetilde G^Z:=G_{\widetilde\Delta^Z}^Z.
\]
This matrix satisfies
\begin{equation}\label{eq:matrix_parameter_deterministic_equation}
    \widetilde G^Z
    =
    \left(
        -Z+\frac1n\sum_{i=1}^n
        \frac{\Sigma_i}
        {1+\frac1n\tr(\Sigma_i\widetilde G^Z)}
    \right)^{-1}.
\end{equation}

We now prove uniqueness using an estimate that will also be used for
stability. With the phase \(\omega_Z\) defined in
\eqref{eq:matrix_parameter_angular_phase}, define the \emph{sector part of
\(w\)} relative to \(\mathcal C_Z\), for a scalar or square matrix \(w\), by
\begin{align}\label{eq:def_weighted_real_part}
    \mathfrak C_Z(w)
    :=\Re(\overline{\omega_Z}w),
    \qquad
    \mathfrak C_{\bar Z}(w)
    =\Re(\omega_Z w).
\end{align}
Here \(\Re\) denotes the real part for scalars and the Hermitian part for
matrices. The second expression is the sector part of \(w\) relative to
the conjugate sector. The second identity follows from
\(\omega_{\bar Z}=\overline{\omega_Z}\), since complex conjugation reflects
the sector and its bisecting phase. In particular,
\(\mathfrak C_Z(1)=\mathfrak C_{\bar Z}(1)\geq c_Z>0\).
For Hermitian \(\Sigma\) and a square matrix \(A\), taking the trace
commutes with taking the sector part:
\[
    \tr\!\left(\Sigma\mathfrak C_{\bar Z}(A)\right)
    =\mathfrak C_{\bar Z}\!\left(\tr(\Sigma A)\right).
\]
For profiles \(D,D'\) for which the two resolvents are defined, set
\[
    \Psi^Z(D,D')_{ij}
    :=
    \frac{D_iD_i'}{n^2}
    \tr\!\left(
        \Sigma_iG_D^Z\Sigma_jG_{D'}^Z
    \right).
\]
For \(D,D'\in\mathcal D_n(\mathcal C_Z\setminus\{0\})\), define the
inverse defect
\begin{equation}\label{eq:def_inverse_defect}
    r_i(D):=\frac1{D_i}-1-\frac1n\tr(\Sigma_iG_D^Z).
\end{equation}
Thus \(r(D)=0\) precisely when \(D\) is a fixed point of \(\Phi^Z\).
The resolvent identity gives
\begin{align}\label{eq:two_profile_resolvent_identity}
    G_D^Z-G_{D'}^Z
    =
    -\frac1nG_D^Z
    \left(\sum_{j=1}^n(D_j-D_j')\Sigma_j\right)
    G_{D'}^Z.
\end{align}
Subtracting \eqref{eq:def_inverse_defect} evaluated at \(D\) and \(D'\), using
\(D_i^{-1}-(D_i')^{-1}=-(D_i-D_i')/(D_iD_i')\), and substituting the
resolvent identity~\eqref{eq:two_profile_resolvent_identity} gives, for each
\(i\in[n]\),
\[
    r_i(D)-r_i(D')
    =
    -\frac{D_i-D_i'}{D_iD_i'}
    +\frac1{n^2}\sum_{j=1}^n(D_j-D_j')
    \tr\!\left(\Sigma_iG_D^Z\Sigma_jG_{D'}^Z\right).
\]
Multiplying the \(i\)-th identity by \(-D_iD_i'\) and using the definition
of \(\Psi^Z(D,D')\) yields the exact secant equation
\begin{align}\label{eq:two_profile_inverse_defect_secant}
    \bigl(I_n-\Psi^Z(D,D')\bigr)(D-D')
    =
    -DD'\bigl(r(D)-r(D')\bigr).
\end{align}
Here \(D-D'\) and \(r(D)-r(D')\) are viewed as vectors, and the
diagonal matrix \(DD'\) acts coordinatewise. For two fixed
points, this reduces to
\(D-D'=\Psi^Z(D,D')(D-D')\); a strict contraction of \(\Psi^Z(D,D')\)
therefore proves uniqueness. Taking \(D'=\widetilde\Delta^Z\) instead
leaves only \(r(D)\) on the right, giving
\eqref{eq:exact_weighted_inverse_defect_secant} below. Bounds on this
defect then control \(D-\widetilde\Delta^Z\), and
\eqref{eq:two_profile_resolvent_identity} transfers them to
\(G_D^Z-\widetilde G^Z\).

For such profiles, equip
\(\mathbb C^n\) with
\[
    |u|_{D,D'}
    :=
    \max_{i\in[n]}
    \frac{|u_i|}
    {\sqrt{\mathfrak C_Z(D_i)\mathfrak C_Z(D_i')}}.
\]
For a linear map \(T:\mathbb C^n\to\mathbb C^n\), denote its induced
operator norm by
\[
    |T|_{D,D'}
    :=
    \sup_{u\neq0}
    \frac{|Tu|_{D,D'}}{|u|_{D,D'}}.
\]

\begin{lemma}[Weighted secant estimate]
\label{lem:weighted_secant_estimate}
For \(D,D'\in\mathcal D_n(\mathcal C_Z\setminus\{0\})\) and every
\(u\in\mathbb C^n\),
\begin{multline}\label{eq:weighted_secant_estimate}
    |\Psi^Z(D,D')u|_{D,D'}
    \leq
    \max_{i\in[n]}
    \sqrt{
        \frac{|D_i|^2\frac1n\mathfrak C_{\bar Z}\!\left(\tr(\Sigma_iG_D^Z)\right)}{\mathfrak C_Z(D_i)}
        \frac{|D_i'|^2\frac1n\mathfrak C_{\bar Z}\!\left(\tr(\Sigma_iG_{D'}^Z)\right)}{\mathfrak C_Z(D_i')}
    }
    |u|_{D,D'}.
\end{multline}
The quantities
\(\frac1n\mathfrak C_{\bar Z}\!\left(\tr(\Sigma_iG_D^Z)\right)\)
are nonnegative.
\end{lemma}

\begin{proof}
For each \(i,j\), the Hilbert--Schmidt Cauchy--Schwarz inequality gives
\begin{align*}
|\Psi^Z(D,D')_{ij}|
&=\frac{|D_iD'_i|}{n^2}\left|
    \tr(\Sigma_iG_D^Z\Sigma_jG_{D'}^Z)
\right|\\
&\leq
\frac{|D_iD'_i|}{n^2}\tr\!\left(
    \Sigma_iG_D^Z\Sigma_j(G_D^Z)^*
\right)^{1/2}
\tr\!\left(
    \Sigma_i(G_{D'}^Z)^*\Sigma_jG_{D'}^Z
\right)^{1/2}.
\end{align*}
Since \(-\mathcal W(Z)\subset\mathcal C_Z\),
Lemma~\ref{lem:matrix_parameter_supporting_phase} gives $\mathfrak C_Z(-Z)\succeq0$. 
Then the two inverse identities in
Lemma~\ref{lem:distance_resolvent_and_ward} give
\begin{align*}
    \frac1n\sum_j\mathfrak C_Z(D_j)
        G_D^Z\Sigma_j(G_D^Z)^*
    &=\mathfrak C_{\bar Z}(G_D^Z)
        -G_D^Z\mathfrak C_Z(-Z)(G_D^Z)^*
    \preceq \mathfrak C_{\bar Z}(G_D^Z).
\end{align*}
Interchanging the two outer resolvent factors on the left gives the same
upper bound \(\mathfrak C_{\bar Z}(G_D^Z)\).
Multiplying each absolute kernel entry by
$\sqrt{\mathfrak C_Z(D_j)\mathfrak C_Z(D_j')}$ and applying Cauchy--Schwarz in $j$
therefore gives
\[
    \sum_{j=1}^n|\Psi^Z(D,D')_{ij}|
       \sqrt{\mathfrak C_Z(D_j)\mathfrak C_Z(D_j')}
    \leq
    \frac{|D_iD_i'|}{n}
    \sqrt{
        \tr(\Sigma_i\mathfrak C_{\bar Z}(G_D^Z))
        \tr(\Sigma_i\mathfrak C_{\bar Z}(G_{D'}^Z))
    }.
\]
Bounding
$|(\Psi^Z(D,D')u)_i|$ by this row sum times $|u|_{D,D'}$ and dividing by
\(\sqrt{\mathfrak C_Z(D_i)\mathfrak C_Z(D_i')}\) proves
\eqref{eq:weighted_secant_estimate}.
\end{proof}

The same estimate holds if the canonical phase is replaced by any
unit phase satisfying the hypotheses of
Lemma~\ref{lem:distance_resolvent_and_ward}, provided both sector parts
and the weights in $|\cdot|_{D,D'}$ are formed with
that phase.

\begin{proposition}[Fixed-point secant contraction and uniqueness]
\label{pro:fixed_point_secant_contraction}
If \(D,D'\in\mathcal D_n(\mathcal C_Z)\) are fixed points of \(\Phi^Z\),
then
\begin{equation}\label{eq:fixed_point_secant_contraction}
    |\Psi^Z(D,D')|_{D,D'}
    \leq
    \max_{i\in[n]}
    \sqrt{
        \frac{\frac1n\mathfrak C_{\bar Z}\!\left(\tr(\Sigma_iG_D^Z)\right)}{\mathfrak C_Z(1)+\frac1n\mathfrak C_{\bar Z}\!\left(\tr(\Sigma_iG_D^Z)\right)}
        \frac{\frac1n\mathfrak C_{\bar Z}\!\left(\tr(\Sigma_iG_{D'}^Z)\right)}{\mathfrak C_Z(1)+\frac1n\mathfrak C_{\bar Z}\!\left(\tr(\Sigma_iG_{D'}^Z)\right)}
    }
    <1.
\end{equation}
Consequently, the fixed point is unique in
\(\mathcal D_n(\mathcal C_Z)\), and
\(I_n-\Psi^Z(D,D)\) is invertible.

\end{proposition}

\begin{proof}
By \eqref{eq:def_matrix_parameter_resolvent}, each component of a fixed
point is the inverse of a finite nonzero scalar. Consequently,
\(\mathfrak C_Z(D_i)\geq c_Z|D_i|>0\), and similarly for \(D'\),
so the weighted norms are well defined. At a fixed point,
\[
    \frac1{D_i}
    =
    1+\frac1n\tr(\Sigma_iG_D^Z).
\]
Taking the sector part relative to \(\overline{\mathcal C_Z}\) gives the exact identity
\[
    \frac{\mathfrak C_Z(D_i)}{|D_i|^2}
    =
    \mathfrak C_Z(1)+\frac1n\mathfrak C_{\bar Z}\!\left(\tr(\Sigma_iG_D^Z)\right).
\]
Substitution in Lemma~\ref{lem:weighted_secant_estimate} proves
\eqref{eq:fixed_point_secant_contraction}; it is strict because
\(\mathfrak C_Z(1)\geq c_Z>0\).

If \(U:=D-D'\), equation
\eqref{eq:two_profile_inverse_defect_secant} with \(r(D)=r(D')=0\) gives
\[
    U=\Psi^Z(D,D')U.
\]
The strict contraction gives \(U=0\), hence \(D=D'\).
Taking \(D'=D\), the same norm bound gives the invertibility of
\(I_n-\Psi^Z(D,D)\) by the Neumann series.
\end{proof}

\begin{proposition}[Sectorial matrix-parameter deterministic equivalent]
\label{pro:matrix_parameter_fixed_point}
For every admissible matrix \(Z\), there is a unique
\(\widetilde\Delta^Z\in\mathcal D_n(\mathcal C_Z)\) satisfying
\[
    \widetilde\Delta^Z=\Phi^Z(\widetilde\Delta^Z).
\]
The associated matrix \(\widetilde G^Z=G_{\widetilde\Delta^Z}^Z\)
satisfies \eqref{eq:matrix_parameter_deterministic_equation}, with
\[
    |\widetilde\Delta^Z|_\infty\leq\frac1{c_Z},
    \qquad
    |\widetilde G^Z|_{\op}\leq\frac1{c_Z\eta_Z}.
\]
\end{proposition}

\begin{proof}
Existence follows from
Lemma~\ref{lem:fixed_point_compact_invariant_set}, and uniqueness from
Proposition~\ref{pro:fixed_point_secant_contraction}.
The bounds follow from the existence lemma and
Lemma~\ref{lem:distance_resolvent_and_ward}.
\end{proof}

\begin{lemma}[Trace identity for the deterministic equivalent]
\label{lem:deterministic_equivalent_trace_identity}
At the matrix-parameter fixed point,
\[
    \tr(Z\widetilde G^Z)
    =
    n-p-\sum_{i=1}^n\widetilde\Delta_i^Z.
\]
\end{lemma}

\begin{proof}
The identity
\[
    \left(
        -Z+\frac1n\sum_i\widetilde\Delta_i^Z\Sigma_i
    \right)\widetilde G^Z=I_p
\]
and the fixed-point equations give
\[
    p
    =
    -\tr(Z\widetilde G^Z)
    +\frac1n\sum_i\widetilde\Delta_i^Z
        \tr(\Sigma_i\widetilde G^Z)
    =
    -\tr(Z\widetilde G^Z)
    +\sum_i(1-\widetilde\Delta_i^Z).
\]
This proves the identity.
\end{proof}

\subsection{Holomorphic dependence}

\begin{proposition}[Holomorphic dependence on the matrix parameter]
\label{pro:matrix_parameter_analyticity}
Let \(U\subset\mathbb C\) be open, and let
\(t\mapsto Z_t\in\mathcal M_p\) have holomorphic entries.
Assume that
\[
    d\bigl(\mathcal W(Z_t),\mathbb R_+\bigr)>0,
    \qquad t\in U.
\]
Then the mappings
\[
    t\longmapsto\widetilde\Delta^{Z_t},
    \qquad
    t\longmapsto\widetilde G^{Z_t}
\]
are holomorphic on \(U\).
\end{proposition}

\begin{proof}
Fix \(t_0\in U\), and choose \(r>0\) so that the closed disk
\(\overline B(t_0,r)\) is contained in \(U\). For \(s,t\) in this disk,
\[
    |\eta_{Z_s}-\eta_{Z_t}|\leq |Z_s-Z_t|_{\op},
\]
because \(|v^*(Z_s-Z_t)v|\leq|Z_s-Z_t|_{\op}\) for every unit vector
\(v\). Consequently, continuity and compactness give
\[
    \eta_*:=\inf_{t\in\overline B(t_0,r)}\eta_{Z_t}>0.
\]
Lemma~\ref{lem:matrix_parameter_sector_aperture} also gives
\[
    c_*:=\inf_{t\in\overline B(t_0,r)}c_{Z_t}
    \geq
    \frac{\eta_*}
    {2\sup_{t\in\overline B(t_0,r)}|Z_t|_{\op}}
    >0.
\]

To identify the solution given by the implicit-function theorem,
we first prove continuity at \(t_0\). Let \(t_m\to t_0\) in this disk,
and write \(Z_m:=Z_{t_m}\) and \(Z_0:=Z_{t_0}\).

We need closedness when the matrix parameter varies.
Let \(u_m\in\mathcal C_{Z_m}\) satisfy \(u_m\to u\).
By the definition of \(\mathcal C_{Z_m}\) and convexity of \(\mathcal W(Z_m)\),
choose \(a_m,b_m\geq0\) and \(w_m\in\mathcal W(Z_m)\) such that
\[
    u_m=a_m-b_mw_m.
\]
The uniform supporting and separation
bounds give \(b_m\leq |u_m|/(c_*\eta_*)=O(1)\).
Since \(d(w_m,\mathcal W(Z_0))\leq|Z_m-Z_0|_{\op}\), we have
\[
    d(u_m,\mathcal C_{Z_0})
    \leq b_m|Z_m-Z_0|_{\op}\longrightarrow0,
\]
and closedness of \(\mathcal C_{Z_0}\) yields
\(u\in\mathcal C_{Z_0}\).

The bound \(|\widetilde\Delta^{Z_m}|_\infty\leq c_*^{-1}\),
supplied by Lemma~\ref{lem:fixed_point_compact_invariant_set}, makes
\((\widetilde\Delta^{Z_m})_m\) bounded.
Every convergent subsequence therefore has a limit in
\(\mathcal D_n(\mathcal C_{Z_0})\), which satisfies the fixed-point
equation by the resolvent bounds and continuity of matrix inversion.
Uniqueness identifies this limit with \(\widetilde\Delta^{Z_0}\).
Consequently,
\[
    \widetilde\Delta^{Z_m}\longrightarrow
    \widetilde\Delta^{Z_0}.
\]

Now \(\Phi^{Z_t}(D)\) denotes the rational
expression in \eqref{eq:def_matrix_parameter_resolvent}, defined on the
open subset of \(U\times\mathcal D_n(\mathbb C)\) where the matrix defining
\(G_D^{Z_t}\) is invertible and all scalar denominators are nonzero.
This open set contains \((t_0,\widetilde\Delta^{Z_{t_0}})\). The mapping
\[
    F(t,D):=D-\Phi^{Z_t}(D)
\]
is therefore holomorphic in a neighborhood of this point.
Its derivative with respect to \(D\) is
\[
    \partial_DF(t_0,\widetilde\Delta^{Z_{t_0}})
    =
    I_n-\Psi^{Z_{t_0}}\!\left(
        \widetilde\Delta^{Z_{t_0}},
        \widetilde\Delta^{Z_{t_0}}
    \right),
\]
which is invertible by
Proposition~\ref{pro:fixed_point_secant_contraction}.
The holomorphic implicit-function theorem gives a holomorphic solution
near \(t_0\), unique among solutions near
\(\widetilde\Delta^{Z_{t_0}}\).
The continuity just proved ensures that this solution equals
\(\widetilde\Delta^{Z_t}\).

Since \(t_0\) is arbitrary,
\(t\mapsto\widetilde\Delta^{Z_t}\) is holomorphic on \(U\).
Holomorphy of \(t\mapsto\widetilde G^{Z_t}\) follows from
\eqref{eq:def_matrix_parameter_resolvent}.
\end{proof}

\subsection{Scalar spectral parameters and Stieltjes transforms}

The scalar quantities of Theorem~\ref{the:fixed_point_existence_uniqueness}
are related to the matrix fixed point by
\[
    \tilde\Gamma_i^z=-\frac1z\widetilde\Delta_i^{zI_p},
    \qquad
    \tilde G^z=\widetilde G^{zI_p},\qquad z\in\Omega.
\]
Write $\tilde\Gamma^z:=\diag(\tilde\Gamma_1^z,\ldots,\tilde\Gamma_n^z)$
and set
\[
    \mathcal C_z:=\mathcal C_{zI_p}
    =\cone\bigl(\mathbb R_+\cup\{-z\}\bigr).
\]
Since $-z^{-1}\mathcal C_z=\overline{\mathcal C_z}$, the following is a
direct corollary of Propositions~\ref{pro:matrix_parameter_fixed_point}
and~\ref{pro:matrix_parameter_analyticity}.

\begin{corollary}[Scalar fixed point and holomorphy]
\label{cor:tilde_Gamma_analytic}
For every $z\in\Omega$, the matrix $\tilde\Gamma^z$ is the unique
$D\in\mathcal D_n(\overline{\mathcal C_z})$ satisfying
\[
 D_i=-\frac1{
  z\left(1+\frac1n\tr(\Sigma_iG_{-zD}^{zI_p})\right)
 },\qquad i\in[n].
\]
Moreover, the mappings
\[
    z\longmapsto\tilde\Gamma^z,
    \qquad
    z\longmapsto\tilde G^z
\]
are holomorphic on \(\Omega\).
\end{corollary}

\begin{proof}[Proof of Theorem~\ref{the:fixed_point_existence_uniqueness}]
Corollary~\ref{cor:tilde_Gamma_analytic} gives existence, uniqueness, and
holomorphy, since
$\overline{\mathcal C_z}=\{a-b/z:a,b\geq0\}$.
Substituting $\widetilde\Delta_i^{zI_p}=-z\tilde\Gamma_i^z$ into
the definition $\tilde G^z=G_{\widetilde\Delta^{zI_p}}^{zI_p}$ gives
the stated formula for $\tilde G^z$.
\end{proof}

Lemma~\ref{lem:deterministic_equivalent_trace_identity}, applied with
$Z=zI_p$, gives
\[
    z\tr(\tilde G^z)
    =n-p+z\sum_{i=1}^n\tilde\Gamma_i^z,
\]
and hence
\begin{align}\label{eq:scalar_companion_transform_identity}
    \tilde g(z)
    =\frac1z\left(\frac np-1\right)
     +\frac1p\sum_{i=1}^n\tilde\Gamma_i^z.
\end{align}
Uniqueness gives the following conjugation identities; the last
equality also uses \((\tilde G^z)^\top=\tilde G^z\).
\begin{equation}\label{eq:scalar_conjugation_identities}
    \tilde\Gamma_i^{\bar z}
    =
    \overline{\tilde\Gamma_i^z},
    \qquad
    \tilde G^{\bar z}
    =
    \overline{\tilde G^z}
    =(\tilde G^z)^*.
\end{equation}

We next prove that each mapping \(z\mapsto\tilde\Gamma_i^z\)
is a Stieltjes transform. We use the characterization
and inversion formula in
Theorem~\ref{the:condition_stieltjes_transform} of
Appendix~\ref{app:auxiliary_results}.

\begin{proposition}\label{pro:tilde_Gamma_stieltjes_transform}
For every \(i\in[n]\), there is a probability measure
\(\tilde\mu_i\) supported on \(\mathbb R_+\) such that
\[
    \tilde\Gamma_i^z
    =
    \int_{\mathbb R_+}\frac{\tilde\mu_i(dt)}{t-z},
    \qquad z\in\Omega.
\]
\end{proposition}

\begin{proof}
By Corollary~\ref{cor:tilde_Gamma_analytic}, the mappings
\(z\mapsto\tilde\Gamma_i^z\) are holomorphic on \(\Omega\).
For \(\Im z>0\), their values lie in
\(\overline{\mathcal C_z}\), so
\(\Im(\tilde\Gamma_i^z)\geq0\) and
\(\Im(z\tilde\Gamma_i^z)\geq0\). The inverse identity gives
\begin{align}\label{eq:scalar_resolvent_imaginary_part}
    \Im(z\tilde G^z)
    =\frac{|z|^2}{n}\sum_{j=1}^n
      \Im(\tilde\Gamma_j^z)
      \tilde G^z\Sigma_j\tilde G^{\bar z}
    \succeq0.
\end{align}
The fixed-point equation therefore yields
\[
    \Im\!\left(-\frac1{\tilde\Gamma_i^z}\right)
    =\Im z+\frac1n\tr\!\left(\Sigma_i\Im(z\tilde G^z)\right)
    >0,
\]
hence \(\Im(\tilde\Gamma_i^z)>0\). Since
\(\Im(-1/\tilde\Gamma_i^z)\leq|\tilde\Gamma_i^z|^{-1}\), the same
identity gives \(|\tilde\Gamma_i^z|\leq(\Im z)^{-1}\), and in particular
\[
    |\tilde\Gamma_i^{\ii y}|\leq y^{-1}.
\]
For fixed \(p,n\), letting \(y\to\infty\) gives
\[
    \frac1n\sum_j\tilde\Gamma_j^{\ii y}\Sigma_j
    \longrightarrow0,
    \qquad
    \ii y\,\tilde G^{\ii y}\longrightarrow-I_p.
\]
The fixed-point equation therefore gives
\[
    \ii y\,\tilde\Gamma_i^{\ii y}\longrightarrow-1.
\]
Theorem~\ref{the:condition_stieltjes_transform} produces a probability
measure on \(\mathbb R\).  Since
\[
    \Im(z\tilde\Gamma_i^z)\geq0,
    \qquad \Im z>0,
\]
the same theorem places its support in \(\mathbb R_+\).  Its Stieltjes transform agrees with \(z\mapsto\tilde\Gamma_i^z\)
on the upper half-plane. Since both functions are holomorphic on
\(\Omega\), they agree there by the identity theorem.
\end{proof}

\begin{proposition}\label{pro:matrix_stieltjes_representation}
The mapping \(z\mapsto\tilde G^z\) is the Stieltjes transform of a positive
matrix-valued measure \(\tilde M\) on \(\mathbb R_+\) with total mass
\(I_p\).  Consequently,
\[
    |\tilde G^z|_{\op}\leq\frac1{\eta_z},
    \qquad z\in\Omega.
\]
\end{proposition}

\begin{proof}
For \(\Im z>0\),
\[
    \tilde G^z
    =
    -\left(
        zI_p+\frac zn\sum_j\tilde\Gamma_j^z\Sigma_j
    \right)^{-1}.
\]
Since \(\Im(z\tilde\Gamma_j^z)\geq0\), the inverse identity gives
\[
    \Im(\tilde G^z)
    =
    \tilde G^{\bar z}
    \left(
        \Im(z)I_p+\frac1n\sum_j\Im(z\tilde\Gamma_j^z)\Sigma_j
    \right)
    \tilde G^z
    \succ0.
\]
Furthermore, \eqref{eq:scalar_resolvent_imaginary_part} and the preceding
proof give
\[
    \Im(z\tilde G^z)
    \succeq0,
    \qquad
    \ii y\,\tilde G^{\ii y}\longrightarrow-I_p.
\]
Fix \(u\in\mathbb C^p\).  For \(u=0\), the corresponding zero measure
gives the assertion trivially.
If \(u\neq0\), apply
Theorem~\ref{the:condition_stieltjes_transform} to the normalized scalar
mapping
\[
    z\longmapsto |u|^{-2}u^*\tilde G^zu.
\]
It has positive imaginary part when \(\Im z>0\); its product with \(z\) has
nonnegative imaginary part there, and
\(\ii y|u|^{-2}u^*\tilde G^{\ii y}u\to-1\).  Consequently,
\(z\mapsto u^*\tilde G^zu\) is the Stieltjes transform of a positive
measure on \(\mathbb R_+\) of total mass \(|u|^2\).
Polarization and uniqueness of scalar Stieltjes transforms give a
positive matrix-valued measure \(\tilde M\) with total mass \(I_p\).
The identity theorem, applied entrywise, extends this representation
to \(\Omega\):
\[
    \tilde G^z
    =
    \int_{\mathbb R_+}\frac{\tilde M(dt)}{t-z},
    \qquad z\in\Omega.
\]
Positivity of the matrix-valued measure and Cauchy--Schwarz for its scalar
polarizations give, for \(u,v\in\mathbb C^p\),
\[
    |u^*\tilde G^zv|
    \leq\frac{|u|\,|v|}{\eta_z}.
\]
Taking the supremum over unit vectors \(u,v\) gives the stated
operator-norm bound.
\end{proof}

\subsection{Support of the deterministic equivalent measure}
\label{sse:deterministic_equivalent_support}

To complete the support assertion in
Theorem~\ref{the:g_tilde_stieltjes}, set
\[
    x_\Sigma
    :=
    \max\left\{\frac{8p}{n},4\right\}|\Sigma|_{\op}^\infty.
\]
\begin{lemma}\label{lem:large_real_part_resolvent_bounds}
Assume that \(|\Sigma|_{\op}^\infty>0\). For every \(z\in\mathbb C\)
with \(\Im z>0\) and \(\Re z\geq x_\Sigma\),
\[
    \Re\left(-\frac1{\tilde\Gamma_i^z}\right)
    \geq\frac{\Re z}{2},
    \qquad i\in[n],
\]
and
\[
    |z\tilde G^z|_{\op}\leq2.
\]
\end{lemma}

\begin{proof}
Put \(x:=\Re z\geq x_\Sigma\).  Consider the compact convex set
\[
    \left\{
        D\in\mathcal D_n(\mathcal C_z):
        |D|_\infty\leq\frac{2|z|}{x}
    \right\}.
\]
It contains $I_n$.  If \(D\) belongs to this set, then
\[
    \left|
        \frac1n\sum_j\left(-\frac{D_j}{z}\right)\Sigma_j
    \right|_{\op}
    \leq\frac{2|\Sigma|_{\op}^\infty}x\leq\frac12,
\]
so the Neumann series gives \(|zG_D^{zI_p}|_{\op}\leq2\).  Hence
\[
    \Re\left(
        z+\frac zn\tr(\Sigma_iG_D^{zI_p})
    \right)
    \geq x-\frac{2p}{n}|\Sigma|_{\op}^\infty
    \geq\frac x2.
\]
The map \(\Phi^{zI_p}\) maps \(\mathcal D_n(\mathcal C_z)\) into itself.  Its image satisfies
\[
    -\frac1z\Phi^{zI_p}(D)_i
    =-\frac1{
        z+\frac zn\tr(\Sigma_iG_D^{zI_p})
      },
\]
and therefore
\[
    |\Phi^{zI_p}(D)_i|
    \leq\frac{2|z|}{x}.
\]
Thus $\Phi^{zI_p}$ maps the displayed set into itself.  Brouwer's fixed-point theorem gives a fixed point in this set. Proposition~\ref{pro:matrix_parameter_fixed_point}
identifies this fixed point with
$\widetilde\Delta^{zI_p}$.  Since
$\tilde\Gamma^z=-z^{-1}\widetilde\Delta^{zI_p}$, the two asserted bounds
follow.
\end{proof}

\begin{proof}[Proof of Theorem~\ref{the:g_tilde_stieltjes}]
Proposition~\ref{pro:tilde_Gamma_stieltjes_transform} provides the
probability measures \(\tilde\mu_i\). By
Proposition~\ref{pro:matrix_stieltjes_representation}, \(\tilde G^z\) is
the Stieltjes transform of a positive matrix-valued measure
\(\tilde M\) on \(\mathbb R_+\) with total mass \(I_p\).
Taking its normalized trace defines the probability measure
\(\tilde\mu:=\frac1p\tr\tilde M\), whose Stieltjes transform is
\(\tilde g\). Equation~\eqref{eq:scalar_companion_transform_identity}
shows that this transform also equals that of the finite signed measure
\((1-n/p)\delta_0+\frac1p\sum_i\tilde\mu_i\).
Uniqueness of Stieltjes transforms therefore proves the stated measure
identity.

If \(|\Sigma|_{\op}^\infty=0\), every \(\Sigma_i\) vanishes, so
\(\tilde\Gamma_i^z=-1/z\) and \(\tilde G^z=-I_p/z\). Thus
\(\tilde\mu_i=\tilde\mu=\delta_0\), and the support assertion follows.
Assume henceforth that \(|\Sigma|_{\op}^\infty>0\), and define
\[
    M(z):=\max_{j\in[n]}\Im(\tilde\Gamma_j^z).
\]
For \(z=x+\ii y\) with \(x\in\mathbb R\) and \(y>0\),
\eqref{eq:scalar_resolvent_imaginary_part} gives
\begin{align}
    \Im\left(-\frac1{\tilde\Gamma_i^z}\right)
    &=
    y+\frac{|z|^2}{n^2}\sum_{j=1}^n
    \Im(\tilde\Gamma_j^z)
    \tr\!\left(
        \Sigma_i\tilde G^z
        \Sigma_j\tilde G^{\bar z}
    \right).
\end{align}
For \(x\geq x_\Sigma\),
Lemma~\ref{lem:large_real_part_resolvent_bounds} gives
\(|z\tilde G^z|_{\op}\leq2\), and hence
\[
    |z|^2\tr\!\left(
        \Sigma_i\tilde G^z
        \Sigma_j\tilde G^{\bar z}
    \right)
    =
    \left|
        \Sigma_i^{1/2}z\tilde G^z\Sigma_j^{1/2}
    \right|_{\hs}^{2}
    \leq4p(|\Sigma|_{\op}^\infty)^2.
\]
Thus
\[
    \Im\left(-\frac1{\tilde\Gamma_i^z}\right)
    \leq
    y+\frac{4p(|\Sigma|_{\op}^\infty)^2}{n}M(z).
\]
The real-part estimate in
Lemma~\ref{lem:large_real_part_resolvent_bounds} yields
\[
    \Im(\tilde\Gamma_i^z)
    \leq
    \frac4{x^2}
    \Im\left(-\frac1{\tilde\Gamma_i^z}\right).
\]
Taking the maximum over \(i\) gives
\[
    M(z)
    \leq
    \frac{4y}{x^2}
    +
    \frac{16p(|\Sigma|_{\op}^\infty)^2}{nx^2}M(z).
\]
The last coefficient is at most \(1/2\) for
\(x\geq x_\Sigma\): if \(p\leq n/2\), use
\(x_\Sigma\geq4|\Sigma|_{\op}^\infty\), while if \(p\geq n/2\), use
\(x_\Sigma\geq8p|\Sigma|_{\op}^\infty/n\).  Therefore
\[
    M(x+\ii y)\leq\frac{8y}{x^2},
    \qquad x\geq x_\Sigma.
\]

Fix \(i\in[n]\), and let \(x_\Sigma<a<b\) satisfy
\(\tilde\mu_i(\{a\})=\tilde\mu_i(\{b\})=0\).  Uniformly for \(x\in[a,b]\),
\[
    0\leq\Im(\tilde\Gamma_i^{x+\ii y})
    \leq\frac{8y}{x^2}.
\]
Letting \(y\downarrow0\) in the Stieltjes inversion formula gives
\(\tilde\mu_i([a,b])=0\).  Exhausting
\((x_\Sigma,\infty)\) by such intervals proves
\[
    \operatorname{supp}(\tilde\mu_i)\subset[0,x_\Sigma].
\]
The measure identity proved above now gives
\[
    \operatorname{supp}(\tilde\mu)
    \subset
    \{0\}\cup\bigcup_{i=1}^n\operatorname{supp}(\tilde\mu_i)
    \subset[0,x_\Sigma]
    =  \left[
          0,\,
          \max\left\{\frac{8p}{n},4\right\}
          |\Sigma|_{\op}^\infty
      \right].
\]
\end{proof}

\section{Probabilistic estimates}
\label{sec:probabilistic_results}

We separate the resolvent error into its fluctuation and its bias.
The common estimates first control $G^Z-\mathbb E[G^Z]$ by independent-copy
replacement. For the bias, each route introduces a deterministic diagonal matrix
\(D\) and estimates both \(\mathbb E[G^Z]-G_D^Z\) and the inverse defect
\(r(D)\) defined in \eqref{eq:def_inverse_defect}.
The Hilbert--Schmidt route uses coordinatewise defect bounds,
whereas the operator route uses a weighted sum of absolute defects.
The two transfer estimates below then control
\(G_D^Z-\widetilde G^Z\).
The final subsection uses additive deformation to obtain the
corresponding sandwich estimates.

\subsection{Common probabilistic estimates and auxiliary pivots}

Fix a nonempty compact set $K\subset\mathbb C$ such that\footnote{The inequality $\eta_K\leq 1$ will be used in particular in \eqref{eq:use_of_eta_leq_1}.}
\begin{equation}\label{eq:def_family_parameter_bounds}
    \eta_K
    :=1\wedge d(\operatorname{conv}K,\mathbb R_+)>0.
\end{equation}
The convex-hull restriction entails no loss of generality for the
scalar results: every compact subset of $\Omega$ admits a finite cover
by compact sets satisfying it. Uniform estimates on the original
compact set follow by taking the maximum of the resulting constants.

Define
\begin{equation}\label{eq:def_common_probabilistic_sector}
    \mathcal C_K
    :=\operatorname{cone}\bigl(\mathbb R_+\cup(-K)\bigr),
\end{equation}
where $\operatorname{cone}$ denotes the convex conic hull.
The geometric proof of
Lemma~\ref{lem:matrix_parameter_sector_aperture}, applied to
$\operatorname{conv}K$, shows that $\mathcal C_K$ is a closed sector
of aperture less than $\pi$.
Let $\omega_K$ be its bisecting unit phase and set
\begin{equation}\label{eq:def_uniform_sector_constant}
    c_K
    :=\inf_{w\in\mathcal C_K\setminus\{0\}}
       \frac{\Re(\overline{\omega_K}w)}{|w|}>0.
\end{equation}
For a scalar or square matrix $w$, define its \emph{sector part} relative
to $\mathcal C_K$ by
\begin{align}\label{eq:def_common_probabilistic_real_parts}
    \mathfrak C_K(w)&:=\Re(\overline{\omega_K}w),
    &
    \mathfrak C_{\bar K}(w)&:=\Re(\omega_Kw).
\end{align}
For matrices, $\Re$ denotes the Hermitian part. The second expression
is the sector part of $w$ relative to the conjugate sector
$\overline{\mathcal C_K}$.
The sector and phase remain fixed throughout the section, including
conditional arguments.

For every $p\geq1$, define
\begin{equation}\label{eq:universal_sectorial_base_class}
    \mathcal K_p
    :=\left\{
        Z\in\mathcal M_p:
        \mathcal C_Z\subseteq\mathcal C_K,
        \quad
        \mathfrak C_K(-Z)\succeq c_K\eta_K I_p
    \right\}.
\end{equation}
Every $z\in K$ gives $zI_p\in\mathcal K_p$.
For any $Z\in\mathcal K_p$, the defining conditions imply
\begin{equation}\label{eq:common_probabilistic_sector_bounds}
    \eta_Z\geq c_K\eta_K>0,
    \qquad c_Z\geq c_K.
\end{equation}


For a parameter family $Z=(Z_p)_{(p,n)\in\Theta}$, write
$Z\in\mathcal K$ when $Z_p\in\mathcal K_p$ at every index; in matrix
formulas, $Z$ denotes the current component, with its indices suppressed.
We use the $O(\cdot)$ convention of
Section~\ref{sse:main_probabilistic_results}, with
$Z\in\mathcal K$ in place of $z\in K$ and with $K$ fixed.
In particular, $|Z|_{\op}=O(1)$ means $\sup_{(p,n)\in\Theta}|Z_p|_{\op}<\infty$.


\subsubsection{Resolvent and leave-one-out estimates}
\label{sse:universal_probabilistic_toolbox}

Given $Z\in\mathcal K$, put
\[
    G:=G^Z=\left(\frac1n XX^\top-Z\right)^{-1}.
\]
For $i\in[n]$, let
\[
    X_{-i}:=(x_1,\ldots,x_{i-1},0,x_{i+1},\ldots,x_n),
    \qquad
    G_{-i}:=G_{-i}^Z:=
    \left(\frac1nX_{-i}X_{-i}^{\top}-Z\right)^{-1}.
\]

\paragraph{Resolvent, sector, and pivot identities.}

Define the scaled Schur pivots by
\begin{align}\label{eq:def_matrix_schur_pivot}
    \Delta_i^Z
    &:=
    \left(
        1+\frac1n x_i^{\top}G_{-i}^Zx_i
    \right)^{-1},
    \qquad i\in[n].
\end{align}
We suppress the superscript $Z$ below and put
$\Delta:=\diag(\Delta_1,\ldots,\Delta_n)$.

Apply Lemma~\ref{lem:distance_resolvent_and_ward} with
$\omega=\omega_K$, $c=c_K$, and
$\kappa=c_K\eta_K$.
Every $T\in\{G,G_{-i}\}$ is well defined and satisfies for all $u\in\mathbb C^p$
\begin{align}\label{eq:universal_sectorial_resolvent_bounds}
 |T|_{\op}&\leq\frac1{c_K\eta_K},\qquad
 u^*Tu\in\overline{\mathcal C_Z},
 \qquad 
 T^*T,\ TT^*\preceq
 \frac1{c_K\eta_K}\mathfrak C_{\bar K}(T).
\end{align}
The sector conclusions and reciprocal bounds of
Lemma~\ref{lem:distance_resolvent_and_ward} also give
\begin{align*}
    x_i^{\top}G_{-i}x_i,\ \tr(\Sigma_iG_{-i})
    &\in\overline{\mathcal C_Z},
    \qquad 
    \Delta_i\in\mathcal C_Z\setminus\{0\},
    \qquad |\Delta_i|\leq\frac1{c_K}.
\end{align*}
Expectations and convex combinations remain in the same sectors. The lemma's reciprocal estimates therefore
also apply to them.

Sherman--Morrison gives
\begin{align}\label{eq:leave_one_out_resolvent_identities}
    G-G_{-i}
    &=
    -\frac{\Delta_i}{n}
      G_{-i}x_ix_i^{\top}G_{-i},
    &
    Gx_i&=\Delta_iG_{-i}x_i,
    &
    x_i^{\top}G&=\Delta_i x_i^{\top}G_{-i}.
\end{align}
Consequently,
\begin{equation}\label{eq:schur_pivot_trace_identity}
    \frac1n x_i^{\top}Gx_i=1-\Delta_i.
\end{equation}

\paragraph{Concentration and leave-one-out estimates.}

Set
\begin{equation}
\label{eq:def_c_nq}
    c_{n,q}:=n^{1/(q\wedge2)-1},\qquad q\geq1.
\end{equation}

\begin{proposition}[Random resolvent fluctuation]
\label{pro:resolvent_concentration}
For every deterministic $A\in\mathcal M_p$ and every $q\geq1$,
\[
    \left\|
        \tr\bigl(A(G-\mathbb E[G])\bigr)
    \right\|_{L^q}
    =
    O\left(
        c_{n,q}|A|_{\op}M_{\op}^{(q)}
    \right).
\]
If $\sqrt p\,|\Sigma|_{\hs}^\infty=O(n)$, then also
\[
    \left\|
        \tr\bigl(A(G-\mathbb E[G])\bigr)
    \right\|_{L^q}
    =O\left(c_{n,q}|A|_{\hs}M_{\hs}^{(q)}\right).
\]
The same estimates hold with $G$ replaced by any $G_{-i}$.
\end{proposition}

\begin{proof}
    Let $x_i'$ be an independent copy of $x_i$, independent of all columns,
and write $G^{(i)}$ for the resolvent obtained by replacing $x_i$ with
$x_i'$. It satisfies the same bounds
\eqref{eq:universal_sectorial_resolvent_bounds} by
Lemma~\ref{lem:distance_resolvent_and_ward}.
Conditionally on $X_{-i}$, write
\begin{align*}
    t&:=x_i^{\top}G_{-i}x_i,
    &s&:=x_i^{\top}G_{-i}AG_{-i}x_i,\\
    t'&:=(x_i')^\top G_{-i}x_i',
    &s'&:=(x_i')^\top G_{-i}AG_{-i}x_i',
\end{align*}
and
\[
    t^\Sigma:=\tr(\Sigma_iG_{-i}),
    \qquad
    s^\Sigma:=\tr(\Sigma_iG_{-i}AG_{-i}).
\]
The rank-one formula gives
\begin{align}\label{eq:replacement_identity_profile}
    \tr\left(A(G-G^{(i)})\right)
    =
    -\frac1n\left(
        \frac{s}{1+t/n}
        -
        \frac{s'}{1+t'/n}
    \right).
\end{align}
Moreover,
\begin{align}\label{eq:self_normalized_decomposition_profile}
    \frac{s}{1+t/n}
    -
    \frac{s^\Sigma}{1+t^\Sigma/n}
    =
    \frac{s-s^\Sigma}{1+t/n}
    -
    \frac{s^\Sigma}{n(1+t^\Sigma/n)}
    \frac{t-t^\Sigma}{1+t/n}.
\end{align}
The term $s^\Sigma/(1+t^\Sigma/n)$ depends only on $X_{-i}$ and
cancels between the two copies in \eqref{eq:replacement_identity_profile}.
Thus each term in the replacement difference contains a quadratic
form centered conditionally on \(X_{-i}\).

We first use operator tests. Hilbert--Schmidt Cauchy--Schwarz gives
\begin{align*}
    |s^\Sigma|
    \leq
    |A|_{\op}
    \sqrt{
        \tr(\Sigma_iG_{-i}^*G_{-i})
        \tr(\Sigma_iG_{-i}G_{-i}^*)
    }
    =O\!\left(|A|_{\op}\mathfrak C_{\bar K}(t^\Sigma)\right),
\end{align*}
thanks to the two sectorial Ward inequalities
\eqref{eq:universal_sectorial_resolvent_bounds}:
\[
    G_{-i}^*G_{-i},\ G_{-i}G_{-i}^*
    \preceq
    \frac{\mathfrak C_{\bar K}(G_{-i})}
         {c_K\eta_K}.
\]
Since $t^\Sigma\in\overline{\mathcal C_Z}$,
Lemma~\ref{lem:distance_resolvent_and_ward} gives
\begin{align}\label{eq:self_normalized_mean_operator_profile}
    \frac{|s^\Sigma|}
         {n|1+t^\Sigma/n|}
    =O(|A|_{\op}).
\end{align}
Conditional quadratic-form control and
\eqref{eq:universal_sectorial_resolvent_bounds} give
\begin{align*}
    \|s-s^\Sigma\|_{L^q(x_i\mid X_{-i})}
    &=
    O\left(M_{\op}^{(q)}|A|_{\op}\right)
    \qquad \text{and}\qquad 
    \|t-t^\Sigma\|_{L^q(x_i\mid X_{-i})}
    &=
    O\left(M_{\op}^{(q)}\right).
\end{align*}
Since \(G_{-i}\) depends only on \(X_{-i}\),
Lemma~\ref{lem:complex_conditional_test_reduction}
applies to the complex test matrices \(G_{-i}AG_{-i}\) and \(G_{-i}\).  Combining these
estimates with
\eqref{eq:self_normalized_decomposition_profile} and
\eqref{eq:self_normalized_mean_operator_profile} yields
\[
    \left\|
        \frac{s}{1+t/n}
        -
        \frac{s^\Sigma}{1+t^\Sigma/n}
    \right\|_{L^q(x_i\mid X_{-i})}
    =
    O\left(M_{\op}^{(q)}|A|_{\op}\right).
\]

For Hilbert--Schmidt tests, conditional quadratic-form control and the
resolvent bounds give
\begin{align*}
    \|s-s^\Sigma\|_{L^q(x_i\mid X_{-i})}
    &=O\left(M_{\hs}^{(q)}|A|_{\hs}\right),\\
    \|t-t^\Sigma\|_{L^q(x_i\mid X_{-i})}
    &=O\left(M_{\hs}^{(q)}\sqrt p\right),
    &|s^\Sigma|&=O\left(|\Sigma_i|_{\hs}|A|_{\hs}\right).
\end{align*}
Together with the uniform reciprocal bounds of
Lemma~\ref{lem:distance_resolvent_and_ward},
\eqref{eq:self_normalized_decomposition_profile} therefore yields
\[
    \left\|\frac{s}{1+t/n}
       -\frac{s^\Sigma}{1+t^\Sigma/n}\right\|_{L^q(x_i\mid X_{-i})}
    =O\left(M_{\hs}^{(q)}|A|_{\hs}
       \left(1+\frac{\sqrt p\,|\Sigma|_{\hs}^\infty}n\right)\right)
    =O\left(M_{\hs}^{(q)}|A|_{\hs}\right).
\]
The bounds for the operator and Hilbert--Schmidt tests also hold
with \(x_i'\) in place of \(x_i\).
In particular, \eqref{eq:replacement_identity_profile} gives
\[
 \|\tr(A(G-G^{(i)}))\|_{L^q(x_i,x_i'\mid X_{-i})}
 =O\left(\frac{M_{\op}^{(q)}|A|_{\op}}n\right).
\]
Under \eqref{eq:hs_profile_size_condition}, it also gives
\[
 \|\tr(A(G-G^{(i)}))\|_{L^q(x_i,x_i'\mid X_{-i})}
 =O\left(\frac{M_{\hs}^{(q)}|A|_{\hs}}n\right).
\]
After integrating these conditional moment bounds, the
replacement-moment inequality in
Theorem~\ref{the:efron_stein_Lq_simple}
gives both fluctuation estimates.
For $G_{-i}$, apply the same argument as a function of
$(x_j)_{j\ne i}$, replacing a column $x_j$ with $j\ne i$; the
normalization remains $1/n$. The case $n=1$ is trivial.
\end{proof}

\begin{lemma}[Mean Hilbert--Schmidt leave-one-out bound]
\label{lem:leave_one_out_profile_bounds}
For every $i\in[n]$,
\begin{align}\label{eq:leave_one_out_hs_profile_bound}
    |\mathbb E[G-G_{-i}]|_{\hs}
    = O\left(\frac{|\Sigma_i|_{\hs}}n\right).
\end{align}
\end{lemma}

\begin{proof}
Fix a deterministic matrix $A$ with $|A|_{\hs}\leq1$.
The rank-one identity \eqref{eq:leave_one_out_resolvent_identities} and
the pivot bound in Lemma~\ref{lem:distance_resolvent_and_ward} imply
\[
    \left|\tr\!\left(A\mathbb E[G-G_{-i}]\right)\right|
    =O\!\left(
      \frac1n\mathbb E\left[\left|x_i^\top G_{-i}AG_{-i}x_i\right|\right]
    \right).
\]
Condition on $X_{-i}$ and apply
Lemma~\ref{lem:uncentered_quadratic_form_bound} to
$B=G_{-i}AG_{-i}$.  Since
$|G_{-i}AG_{-i}|_{\hs}=O(|A|_{\hs})$,
\[
    \left|\tr\!\left(A\mathbb E[G-G_{-i}]\right)\right|
    = O\left(\frac{|\Sigma_i|_{\hs}}n\right).
\]
Hilbert--Schmidt duality proves
\eqref{eq:leave_one_out_hs_profile_bound}.
\end{proof}

\subsubsection{Auxiliary resolvents and inverse defects}

Fix \(Z\in\mathcal K\), and write
\[
    \widetilde\Delta:=\widetilde\Delta^Z,
    \qquad
    \widetilde G:=\widetilde G^Z=G_{\widetilde\Delta}^Z.
\]
For a deterministic \(D\in\mathcal D_n(\mathcal C_Z)\), set
\[
    G_D:=G_D^Z
    :=
    \left(
        \frac1n\sum_{i=1}^nD_i\Sigma_i-Z
    \right)^{-1}.
\]
Applying Lemma~\ref{lem:distance_resolvent_and_ward} with the same sectorial
constants shows that $G_D$ is well defined and satisfies
\eqref{eq:universal_sectorial_resolvent_bounds}.
The two probabilistic routes use the same decomposition
\[
    G-\widetilde G
    =
    (G-\mathbb E[G])
    +
    (\mathbb E[G]-G_D)
    +
    (G_D-\widetilde G).
\]
Proposition~\ref{pro:resolvent_concentration} controls the first term.
The two routes estimate the middle term using their respective pivots.

To estimate \(G_D-\widetilde G\), we use the inverse defect \(r(D)\)
defined in \eqref{eq:def_inverse_defect}. For
\(D\in\mathcal D_n(\mathcal C_Z\setminus\{0\})\), uniqueness of the
deterministic fixed point gives \(r(D)=0\) exactly when
\(D=\widetilde\Delta\). Applying
\eqref{eq:two_profile_resolvent_identity} with
\(D'=\widetilde\Delta\) gives
\begin{align}
    G_D-\widetilde G
    =
    -\frac1n
    G_D
    \left(
        \sum_{j=1}^n
        (D_j-\widetilde\Delta_j)\Sigma_j
    \right)
    \widetilde G.
\end{align}
The same specialization in \eqref{eq:two_profile_inverse_defect_secant},
using \(r(\widetilde\Delta)=0\), gives
\begin{align}
\label{eq:exact_weighted_inverse_defect_secant}
    D-\widetilde\Delta
    =
    -D\widetilde\Delta\,r(D)
    +
    \Psi^Z(D,\widetilde\Delta)
    (D-\widetilde\Delta).
\end{align}

The secant equation \eqref{eq:exact_weighted_inverse_defect_secant}
transfers a bound on \(r(D)\) to \(G_D-\widetilde G\).
The appropriate measure of \(r(D)\) depends
on the test geometry: the Hilbert--Schmidt route uses
\(|r(D)|_\infty\), whereas the operator-norm route uses
\(\sum_i|D_ir_i(D)|\).

\subsubsection{Choice of the auxiliary pivot}
\label{sse:universal_profile_pivots}

We choose \(D\) so that the estimate of \(\mathbb E[G]-G_D\) also
gives the defect bound needed in
Proposition~\ref{pro:two_inverse_defect_transfers}.

\paragraph{Hilbert--Schmidt-norm route.}
Define
\[
    \check\Delta_i
    :=
    \bigl(\mathbb E[\Delta_i^{-1}]\bigr)^{-1}
    =
    \left(
        1+\frac1n
        \tr(\Sigma_i\mathbb E[G_{-i}])
    \right)^{-1},
\]
and set
\[
    \check\Delta
    :=
    \diag(\check\Delta_1,\ldots,\check\Delta_n),
    \qquad
    \check G:=\check G^Z:=G_{\check\Delta}^Z.
\]
The sector and reciprocal bounds in
Lemma~\ref{lem:distance_resolvent_and_ward} give
\[
    \check\Delta_i\in\mathcal C_Z\setminus\{0\},
    \qquad
    |\check\Delta|_\infty=O(1),
\]
so \(\check G\) is well defined.  By the definitions of \(\check\Delta_i\) and \(r_i\),
\begin{align}
\label{eq:check_pivot_exact_defect}
    r_i(\check\Delta)
    =
    \frac1n
    \tr\!\left(
        \Sigma_i(
            \mathbb E[G_{-i}]-\check G
        )
    \right).
\end{align}
Consequently,
\[
    |r_i(\check\Delta)|
    \leq
    \frac{|\Sigma_i|_{\hs}}n
    \left(
        |\mathbb E[G]-\check G|_{\hs}
        +
        |\mathbb E[G_{-i}-G]|_{\hs}
    \right).
\]
Thus the mean-resolvent bias and the mean leave-one-out error control
\(|r(\check\Delta)|_\infty\) directly. The check pivot also satisfies
\[
 \check\Delta_i-\Delta_i
 =\frac{\check\Delta_i\Delta_i}{n}
 \left(x_i^\top G_{-i}x_i
       -\mathbb E[x_i^\top G_{-i}x_i]\right).
\]
The uniform bounds on the pivots reduce the estimate of
\(\check\Delta_i-\Delta_i\) to the centered quantity on the right-hand side. The factor \(|\Sigma_i|_{\hs}/n\) in the bound for
\(r_i(\check\Delta)\) comes from the Hilbert--Schmidt trace inequality.

\paragraph{Operator-norm route.}
For tests \(A\) measured in operator norm, we bound
\(G_D-\widetilde G\) in nuclear norm. The transfer uses the weighted defect
\[
    \sum_{i=1}^n|D_ir_i(D)|.
\]
For this bound, we use expectations of the Schur pivots.
First suppose that \(\Sigma_i=\Sigma\) for every \(i\), and set
\[
    \widehat\delta
    :=
    \frac1n\sum_{i=1}^n\mathbb E[\Delta_i],
    \qquad
    \widehat G
    :=
    (\widehat\delta\Sigma-Z)^{-1}.
\]
The two resolvent equations are
\begin{align}
    \frac1nXX^\top G
    &=
    I_p+ZG,
    &
    \widehat\delta\Sigma\widehat G
    &=
    I_p+Z\widehat G.
\end{align}
Since
\[
    \mathbb E\!\left[
        \tr\!\left(
            \frac1nXX^\top G
        \right)
    \right]
    =
    n(1-\widehat\delta),
\]
the common inverse defect satisfies
\[
    n\widehat\delta\,
    r_i(\widehat\delta I_n)
    =
    \tr\!\left(
        Z(\mathbb E[G]-\widehat G)
    \right).
\]
This identity controls the weighted defect by testing
\(\mathbb E[G]-\widehat G\) against \(Z\).

For several classes, the same identity is applied at each replacement
stage using the expected pivots defined in
\eqref{eq:class_stage_pivots}.  If the final deterministic pivot is constant on \(I_a\), with
\(D_i=\widehat\delta^{(a)}\) for \(i\in I_a\), then
\[
    n_a(1-D_i)
    -
    \pi_aD_i
    \tr(\Sigma^{(a)}G_D)
    =
    n_aD_ir_i(D),
    \qquad i\in I_a.
\]
The successive replacement argument estimates the change from the
expected random contribution \(n_a(1-D_i)\) to the final deterministic
contribution
\(\pi_aD_i\tr(\Sigma^{(a)}G_D)\).  Summing the absolute error bounds over the classes controls
\(\sum_i|D_i r_i(D)|\).  Replacing the classwise
expected pivots by the check pivots would introduce
\[
    \bigl(\mathbb E[\Delta_i^{-1}]\bigr)^{-1}
    -
    \mathbb E[\Delta_i],
\]
which would require a separate estimate in this argument.

Accordingly, the two routes use different pivots for a precise
algebraic reason:
\[
    \check\Delta
    \quad\longrightarrow\quad
    |r(\check\Delta)|_\infty
    \quad\longrightarrow\quad
    |\check G-\widetilde G|_{\hs},
\]
whereas
\[
    D_i=\widehat\delta^{(a)}
    \ \text{on }I_a
    \quad\longrightarrow\quad
    \sum_i|D_ir_i(D)|
    \quad\longrightarrow\quad
    |G_D-\widetilde G|_*.
\]

\subsubsection{Transfer from the inverse defect}
\label{sse:weighted_pivot_transfer}

The following proposition bounds \(G_D^Z-\widetilde G^Z\) in
Hilbert--Schmidt and nuclear norm using the secant equation
\eqref{eq:exact_weighted_inverse_defect_secant}.
Part~\textnormal{(i)} uses a weighted supremum norm, while
part~\textnormal{(ii)} uses weighted column sums.

For the weighted secant estimates below, use
\[
    |u|_{D,D'}
    :=
    \max_{i\in[n]}
    \frac{|u_i|}
    {\sqrt{\mathfrak C_K(D_i)
           \mathfrak C_K(D_i')}},
\]
and the corresponding induced operator norm.

\begin{proposition}[Two transfers from the inverse defect]
\label{pro:two_inverse_defect_transfers}
Let \(Z\in\mathcal K\) and

\[
    D\in\mathcal D_n(\mathcal C_Z\setminus\{0\}).
\]

Write \(\widetilde\Delta=\widetilde\Delta^Z\) and
\(\widetilde G=\widetilde G^Z\).

\begin{enumerate}

\item[\textnormal{(i)}]
Assume that

\[
    |\Sigma|_*^\infty=O(n),
    \qquad
    |D|_\infty=O(1).
\]

Then
\begin{align}
\label{eq:inverse_defect_stability_hs}
|G_D^Z-\widetilde G|_{\hs}
=
O\!\left(
\min\{|\Sigma|_{\hs}^\infty,\sqrt p\}
|r(D)|_\infty
\right).
\end{align}

\item[\textnormal{(ii)}]
For every \(D\in\mathcal D_n(\mathcal C_Z\setminus\{0\})\) such that $|D|_\infty=O(1)$, assume there exists $\rho>1$ satisfying $|D_i|^{-1}=O(\rho)$ and $|\widetilde\Delta_i|^{-1}=O(\rho)$. Then
\begin{align}
\label{eq:sectorial_weighted_defects_bounded_pivots}
|G_D^Z-\widetilde G|_*
=
O\!\left(
\sqrt\rho
\sum_{i=1}^n |D_i r_i(D)|
\right).
\end{align}

\end{enumerate}
\end{proposition}

\begin{proof}
Since \(r(\widetilde\Delta)=0\), the secant equation reads
\begin{equation}
\label{eq:inverse_defect_secant_used_here}
D-\widetilde\Delta
=
-D\widetilde\Delta r(D)
+
\Psi^Z(D,\widetilde\Delta)(D-\widetilde\Delta) .
\end{equation}

\medskip
\noindent
\emph{Proof of \textnormal{(i)}.}
The argument is the same weighted contraction as in the uniqueness proof,
with the defect \(r(D)\) acting as a forcing term.

The assumption \( |\Sigma|_*^\infty=O(n)\), together with the sectorial
resolvent bound, gives
\[
    \frac1n\mathfrak C_{\bar K}(\tr(\Sigma_iG_D^Z))+\frac1n\mathfrak C_{\bar K}(\tr(\Sigma_i\widetilde G^Z))=O(1)
\]
uniformly in \(i\).
Suppose first that
\[
    |r(D)|_\infty\leq \frac{c_K}2.
\]
Taking the sector part relative to \(\overline{\mathcal C_K}\) in the
definition of \(r_i(D)\), and using the
fixed-point equation for \(\widetilde\Delta\), gives
\begin{align*}
\frac{\mathfrak C_K(D_i)}{|D_i|^2}
&=
\mathfrak C_K(1)+\frac1n\mathfrak C_{\bar K}(\tr(\Sigma_iG_D^Z))
+\mathfrak C_{\bar K}(r_i(D))
    \geq
\frac{c_K}2+\frac1n\mathfrak C_{\bar K}(\tr(\Sigma_iG_D^Z)),
\\
\frac{\mathfrak C_K(\widetilde\Delta_i)}
{|\widetilde\Delta_i|^2}
&=
\mathfrak C_K(1)+\frac1n\mathfrak C_{\bar K}(\tr(\Sigma_i\widetilde G^Z))
    \geq
c_K+\frac1n\mathfrak C_{\bar K}(\tr(\Sigma_i\widetilde G^Z)).
\end{align*}
Hence Lemma~\ref{lem:weighted_secant_estimate} gives, uniformly in
\(Z\in\mathcal K\),
\[
    |\Psi^Z(D,\widetilde\Delta)|_{D,\widetilde\Delta}
    \leq 1-c
\]
for some \(c>0\). Applying this to
\eqref{eq:inverse_defect_secant_used_here} therefore yields
\[
 |D-\widetilde\Delta|_{D,\widetilde\Delta}
 =O\!\left(
  \max_{i\in[n]}
  \frac{|D_i\widetilde\Delta_i r_i(D)|}
       {\sqrt{\mathfrak C_K(D_i)
                    \mathfrak C_K(\widetilde\Delta_i)}}
 \right)
 =O(|r(D)|_\infty).
\]
Since \(D\) is bounded and
\(|\widetilde\Delta|_\infty=O(1)\) by
Proposition~\ref{pro:matrix_parameter_fixed_point}, this also gives
\begin{equation}
|D-\widetilde\Delta|_\infty
=
O(|r(D)|_\infty).
\end{equation}
The resolvent identity now gives directly
\[
    |G_D^Z-\widetilde G|_{\hs}
    \leq
    \frac1n\sum_i
        |D_i-\widetilde\Delta_i|\,
        |G_D^Z\Sigma_i\widetilde G|_{\hs}
    =
    O\!\left(
        |\Sigma|_{\hs}^\infty |r(D)|_\infty
    \right).
\]
We also obtain a bound with \(\sqrt p\) in place of
\(|\Sigma|_{\hs}^\infty\). Using
\[
    |D_i-\widetilde\Delta_i|
    \leq
    |D-\widetilde\Delta|_{D,\widetilde\Delta}
    \sqrt{\mathfrak C_K(D_i)\mathfrak C_K(\widetilde\Delta_i)}
\]
in the same resolvent identity, followed by Cauchy--Schwarz in \(i\) and
the two Ward inequalities from
Lemma~\ref{lem:distance_resolvent_and_ward}, gives
\[
    |G_D^Z-\widetilde G|_{\op}
    =
    O\!\left(
        |D-\widetilde\Delta|_{D,\widetilde\Delta}
    \right)
    =
    O(|r(D)|_\infty).
\]
Consequently,
\[
    |G_D^Z-\widetilde G|_{\hs}
    =
    O(\sqrt p\,|r(D)|_\infty).
\]
Taking the better of the two estimates proves
\eqref{eq:inverse_defect_stability_hs} when
\(|r(D)|_\infty\leq c_K/2\).
If instead
\(|r(D)|_\infty>c_K/2\), the resolvent identity and the uniform bounds on
\(D,\widetilde\Delta,G_D^Z,\widetilde G\) give
\[
    |G_D^Z-\widetilde G|_{\hs}
    =
    O\!\left(
        \min\{|\Sigma|_{\hs}^\infty,\sqrt p\}
    \right),
\]
and the factor \(|r(D)|_\infty\) may be inserted because it is bounded
below by a positive constant. This proves \textnormal{(i)}.
\medskip

\noindent
\emph{Proof of \textnormal{(ii)}.}
Here we retain the term involving the sector part
\(\mathfrak C_K(-Z)\) in the Ward identities.
Since \(Z\in\mathcal K\), we have
\[
    \mathfrak C_K(-Z)
    \succeq c_K\eta_KI_p.
\]
This term will play two roles: it controls the nuclear norm of the
resolvent difference, and it appears as the deficit in the weighted column sums of
the secant kernel.
For \(i\in[n]\), set
\begin{align*}
    h_i
    :=
    \frac1n
    \tr\!\left(
        \Sigma_i(G_D^Z)^*\mathfrak C_K(-Z)G_D^Z
    \right),
    \qquad 
    \widetilde h_i
    :=
    \frac1n
    \tr\!\left(
        \Sigma_i\widetilde G \mathfrak C_K(-Z)\widetilde G^*
    \right).
\end{align*}
By coercivity of \(\mathfrak C_K(-Z)\) and the resolvent identity,
\begin{align}
\label{eq:nuclear_transfer_hi}
|G_D^Z-\widetilde G|_*
&\leq
\frac1n\sum_i
|D_i-\widetilde\Delta_i|\,
|G_D^Z\Sigma_i\widetilde G|_*
=
O\left(
\sum_i |D_i-\widetilde\Delta_i|\sqrt{h_i\widetilde h_i}
\right).
\end{align}
Thus it remains to control the sum on the right.
Use the column weights
\[
    w_i
    :=
    \frac{
        \sqrt{\mathfrak C_K(D_i)\mathfrak C_K(\widetilde\Delta_i)}
    }{
        |D_i\widetilde\Delta_i|
    }.
\]
These weights cancel \(|D_i\widetilde\Delta_i|\) in the \(i\)-th row
of the secant matrix. Cauchy--Schwarz and the Ward identities then give
\begin{align}
\label{eq:column_secant_bound_simple}
\sum_i
w_i
|\Psi^Z(D,\widetilde\Delta)_{ij}|
&\leq \frac1{n^2} \sqrt{
\sum_i
\mathfrak C_K(D_i)
\tr\left(
\Sigma_iG_D^Z
\Sigma_j(G_D^Z)^*
\right)}\sqrt{\sum_i
\mathfrak C_K(\widetilde\Delta_i)
\tr\left(
\Sigma_i\widetilde G^*
\Sigma_j\widetilde G
\right)} \nonumber\\
&=
\sqrt{
\frac1n\mathfrak C_{\bar K}(\tr(\Sigma_jG_D^Z))-h_j
}
\sqrt{
\frac1n\mathfrak C_{\bar K}(\tr(\Sigma_j\widetilde G^Z))-\widetilde h_j
}.
\end{align}
Applying the two-term Cauchy--Schwarz inequality to
\eqref{eq:column_secant_bound_simple} gives
\begin{equation}
\sum_i
w_i|\Psi^Z(D,\widetilde\Delta)_{ij}|
+
\sqrt{h_j\widetilde h_j}
    \leq
\sqrt{
\frac1n\mathfrak C_{\bar K}(\tr(\Sigma_jG_D^Z))
}
\sqrt{
\frac1n\mathfrak C_{\bar K}(\tr(\Sigma_j\widetilde G^Z))
}.
\end{equation}
The inverse-defect equation \eqref{eq:def_inverse_defect} and the
fixed-point equation in Proposition~\ref{pro:matrix_parameter_fixed_point}
imply
\begin{align*}
 \frac1n\mathfrak C_{\bar K}(\tr(\Sigma_jG_D^Z))
 &\leq\frac{\mathfrak C_K(D_j)}{|D_j|^2}+|r_j(D)|,
 \qquad \frac1n\mathfrak C_{\bar K}(\tr(\Sigma_j\widetilde G^Z))
 \leq\frac{\mathfrak C_K(\widetilde\Delta_j)}
              {|\widetilde\Delta_j|^2}.
\end{align*}
All quantities under the square roots are nonnegative, and
$\mathfrak C_K(D_j)>0$ by the sector assumptions. Using
$\sqrt{1+t}\leq1+t/2$ for $t\geq0$ therefore yields
\begin{equation}
\label{eq:column_deficit_with_defect}
\sum_i
w_i|\Psi^Z(D,\widetilde\Delta)_{ij}|
+
\sqrt{h_j\widetilde h_j}
    \leq
w_j
+
\frac{|r_j(D)|}{2}
\sqrt{
\frac{
\mathfrak C_K(\widetilde\Delta_j)/
|\widetilde\Delta_j|^2
}{
\mathfrak C_K(D_j)/|D_j|^2
}
}.
\end{equation}
We can now return to the secant equation
\eqref{eq:inverse_defect_secant_used_here}. Taking absolute values,
multiplying the \(i\)-th coordinate by \(w_i\), and summing over \(i\), gives
\begin{align*}
\sum_iw_i|D_i-\widetilde\Delta_i|
&\leq
\sum_i
|D_i\widetilde\Delta_i r_i(D)|w_i
+
\sum_j|D_j-\widetilde\Delta_j|
\sum_iw_i
|\Psi^Z(D,\widetilde\Delta)_{ij}|.
\end{align*}
Insert \eqref{eq:column_deficit_with_defect} and cancel
$\sum_i w_i|D_i-\widetilde\Delta_i|$ from the two sides to obtain
\begin{align}
\label{eq:key_nuclear_pivot_transfer}
\sum_i|D_i-\widetilde\Delta_i|\sqrt{h_i\widetilde h_i}
\leq
\sum_i |r_i(D)|
\left[
\sqrt{
\mathfrak C_K(D_i)\mathfrak C_K(\widetilde\Delta_i)
}
+
\frac{|D_i-\widetilde\Delta_i|}{2}
\sqrt{
\frac{
\mathfrak C_K(\widetilde\Delta_i)/
|\widetilde\Delta_i|^2
}{
\mathfrak C_K(D_i)/|D_i|^2
}
}
\right].
\end{align}

Since $c_K|u|\leq\mathfrak C_K(u)\leq|u|$ for every
$u\in\mathcal C_Z$, we have 
\[
    \sqrt{
        \mathfrak C_K(D_i)\mathfrak C_K(\widetilde\Delta_i)
    }
    =
    O\!\left(
        \sqrt{|D_i\widetilde\Delta_i|}
    \right)
\qquad \text{and}\qquad  
    \sqrt{
        \frac{
            \mathfrak C_K(\widetilde\Delta_i)/
            |\widetilde\Delta_i|^2
        }{
            \mathfrak C_K(D_i)/|D_i|^2
        }
    }
    =
    O\!\left(
        \sqrt{\frac{|D_i|}{|\widetilde\Delta_i|}}
    \right).
\] 
Using also 
$|D_i-\widetilde\Delta_i|
    \leq
    |D_i|+|\widetilde\Delta_i|$,
we can deduce from \eqref{eq:nuclear_transfer_hi} and~\eqref{eq:key_nuclear_pivot_transfer}: 
\begin{align}
    \label{eq:sectorial_weighted_defects_transfer}
|G_D^Z-\widetilde G|_*\leq
\frac1{c_K\eta_K}
\sum_i|D_i-\widetilde\Delta_i|\sqrt{h_i\widetilde h_i}
\leq
    O\!\left(
        \sum_i
        |D_i r_i(D)|
        \left[
            \sqrt{\frac{|D_i|}{|\widetilde\Delta_i|}}
            +
            \sqrt{\frac{|\widetilde\Delta_i|}{|D_i|}}
        \right]
    \right).
\end{align}
Finally, under the stated bounded-pivot assumptions, 
\[
    \frac{|D_i|}{|\widetilde\Delta_i|}
    +
    \frac{|\widetilde\Delta_i|}{|D_i|}
    =
    O(\rho),
\] 
uniformly in \(i\), and hence each square-root ratio in
\eqref{eq:sectorial_weighted_defects_transfer} is
\(O(\sqrt\rho)\). This proves
\eqref{eq:sectorial_weighted_defects_bounded_pivots}.
\end{proof} 

\begin{corollary}[Inverse-defect transfer for negative Hermitian parameters]
\label{cor:negative_hermitian_inverse_defect_transfer}
Let \(Z\in\mathcal K\) satisfy \(Z=Z^*\prec0\), and let
\(D\in\mathcal D_n((0,\infty))\). Then
\begin{align}
    |G_D^Z-\widetilde G^Z|_*
    \leq
    \frac1{\eta_Z}
    \sum_{i=1}^n|D_i r_i(D)|.
\end{align}
The coefficient $\eta_Z^{-1}$ is uniformly bounded on $\mathcal K$
by \eqref{eq:common_probabilistic_sector_bounds}.
\end{corollary}

\begin{proof}
Set
\[
    \alpha_i
    :=\frac{D_i}{\Phi^Z(D)_i}
    =1-D_i r_i(D)>0.
\]
For nonnegative diagonal coefficients, increasing one entry decreases
the resolvent in the positive-semidefinite order.
Hence \(\Phi^Z\) is coordinatewise increasing. Consequently, the continuous map
\[
    E\longmapsto
    \diag\bigl(\min\{\alpha_i,1\}\Phi^Z(E)_i:i\in[n]\bigr)
\]
maps the box \(0\leq E_i\leq\min\{D_i,\widetilde\Delta_i^Z\}\) into
itself: its \(i\)-th coordinate is at most both
\(\alpha_i\Phi^Z(D)_i=D_i\) and
\(\Phi^Z(\widetilde\Delta^Z)_i=\widetilde\Delta_i^Z\).
Brouwer's fixed-point theorem, as in
Lemma~\ref{lem:fixed_point_compact_invariant_set},
gives a fixed point \(E\) in this box.
Its coordinates are positive because every coordinate of the map
is positive. Since \(E_i\leq\min\{D_i,\widetilde\Delta_i^Z\}\),
\[
    G_E^Z\succeq G_D^Z,
    \qquad
    G_E^Z\succeq\widetilde G^Z.
\]
Taking traces in the resolvent equations
\eqref{eq:def_matrix_parameter_resolvent} and using
\(E_i/\Phi^Z(E)_i=\min\{\alpha_i,1\}\) gives
\begin{align*}
    \tr\!\left((-Z)(G_E^Z-G_D^Z)\right)
    &=
    \sum_i\bigl[\alpha_i-\min\{\alpha_i,1\}-(D_i-E_i)\bigr]
    \leq\sum_i\bigl(\alpha_i-\min\{\alpha_i,1\}\bigr),\\
    \tr\!\left((-Z)(G_E^Z-\widetilde G^Z)\right)
    &=
    \sum_i\bigl[1-\min\{\alpha_i,1\}
        -(\widetilde\Delta_i^Z-E_i)\bigr]
    \leq\sum_i\bigl(1-\min\{\alpha_i,1\}\bigr).
\end{align*}
Both differences are positive semidefinite, so their nuclear norms
equal their traces. Using \(-Z\succeq\eta_ZI_p\) and
the triangle inequality gives
\[
    |G_D^Z-\widetilde G^Z|_*
    \leq
    |G_E^Z-G_D^Z|_*+|G_E^Z-\widetilde G^Z|_*
    \leq
    \frac1{\eta_Z}\sum_i|\alpha_i-1|
    =
    \frac1{\eta_Z}\sum_i|D_i r_i(D)|.
\]
\end{proof}

\subsection{Hilbert--Schmidt-norm route}
\label{sse:hs_profile_route}
\label{sse:limited_second_moment_profiles}

We use \(\check G^Z\) and \(\check\Delta\) as defined in
Section~\ref{sse:universal_profile_pivots}.
We first bound \(\mathbb E[G^Z]-\check G^Z\) and
\(r(\check\Delta)\) in terms of \(J(\check\Delta)\), defined below.
Under the additional assumption \(|Z|_{\op}=O(1)\), we then show that
\(J(\check\Delta)=O(\sigma_\Sigma)\) and apply the Hilbert--Schmidt transfer.

For a fixed $Z\in\mathcal K$,
define
\begin{align}\label{eq:def_matrix_check_pivot_interaction}
    J(D)
    &:={}
    1+\frac1n\sum_{i=1}^n
    |D_i|\,|\Sigma_iG_D^Z|_{\op},
    \qquad D\in\mathcal D_n(\mathcal C_Z).
\end{align}

\begin{proposition}[Bias and inverse defect of the check pivot]
\label{pro:interaction_controlled_check_pivot}
Assume that $\sqrt p\,|\Sigma|_{\hs}^\infty=O(n)$. Then
\begin{align}
 |\mathbb E[G^Z]-\check G^Z|_{\hs}
 &=O\left(J(\check\Delta)\left[
 \frac{|\Sigma|_{\hs}^\infty}n
 +\min\left\{
   M_{\hs}^{(1)},
   \frac{(M_{\hs}^{(2)})^2\sqrt p}{n}
 \right\}\right]\right),
 \label{eq:interaction_check_bias}\\
 |r(\check\Delta)|_\infty
 &=O\left(\frac{|\Sigma|_{\hs}^\infty}nJ(\check\Delta)\left[
 \frac{|\Sigma|_{\hs}^\infty}n
 +\min\left\{
   M_{\hs}^{(1)},
   \frac{(M_{\hs}^{(2)})^2\sqrt p}{n}
 \right\}\right]\right).
 \label{eq:interaction_check_defect}
\end{align}
\end{proposition}

\begin{proof}
Put
\[
 \varepsilon:=\frac{|\Sigma|_{\hs}^\infty}n
 +\min\left\{
  M_{\hs}^{(1)},
  \frac{(M_{\hs}^{(2)})^2\sqrt p}{n}
 \right\}.
\]
We first record two consequences of the profile-size condition:
\[
 \frac{|\Sigma|_*^\infty}n=O(1),
 \qquad \frac{|\Sigma|_{\hs}^\infty}n=O\!\left(\frac1{\sqrt p}\right)=O(1).
\]
The check-pivot definition in Section~\ref{sse:universal_profile_pivots}
and the uniform resolvent bound
\eqref{eq:universal_sectorial_resolvent_bounds} give
\[
 |\check\Delta_i|^{-1}
 \leq 1+C\frac{\tr(\Sigma_i)}n=O(1),
\]
so the definition of \(J\) gives $\frac1n\sum_i|\Sigma_i\check G^Z|_{\op}=O\left(J(\check\Delta)\right)$.

To estimate the bias by Hilbert--Schmidt duality, fix a deterministic
matrix $A$ with $|A|_{\hs}\leq1$, and set
\[
 t_i=x_i^\top G_{-i}^Z x_i,
 \qquad s_i=x_i^\top\check G^Z A G_{-i}^Z x_i.
\]
Independence gives
$\mathbb Es_i=\tr(\Sigma_i\check G^Z A\,\mathbb E[G_{-i}^Z])$.
The resolvent identity and $G^Zx_i=\Delta_iG_{-i}^Zx_i$ therefore give
\begin{align*}
 \tr\bigl(A(\mathbb E[G^Z]-\check G^Z)\bigr)
 &=\frac1n\sum_i\left[
     \check\Delta_i\tr\!\left(\Sigma_i\check G^Z A\,\mathbb E[G^Z]\right)
     -\mathbb E[\Delta_i s_i]
   \right]
   =\frac1n\sum_i(a_i+b_i),
\end{align*}
where
\begin{align*}
 a_i&=\check\Delta_i\tr\left(
       \Sigma_i\check G^Z A\,\mathbb E[G^Z-G_{-i}^Z]\right)
       \qquad \text{and}\qquad
       b_i=\mathbb E\left[s_i(\check\Delta_i-\Delta_i)\right].
\end{align*}
For \(a_i\), the mean leave-one-out estimate
\eqref{eq:leave_one_out_hs_profile_bound} gives
\[
 \max_i|\mathbb E[G^Z-G_{-i}^Z]|_{\hs}
 =O\left(\frac{|\Sigma|_{\hs}^\infty}n\right)=O(\varepsilon),
\]
since $|\Sigma|_{\hs}^\infty/n\leq\varepsilon$ by definition.
The definition of $J$ then yields
\[
 \frac1n\sum_i|a_i|
 =O\left(J(\check\Delta)\frac{|\Sigma|_{\hs}^\infty}n\right).
\]

For $b_i$, the check-pivot definition gives directly
\[
 b_i
 =\mathbb E\left[s_i\left(
   \frac1{1+\mathbb Et_i/n}-\frac1{1+t_i/n}
 \right)\right].
\]
We use Lemma~\ref{lem:reciprocal_difference_estimate} from
Appendix~\ref{app:auxiliary_results}. Let $\mathbb E_i$ denote
conditional expectation given $X_{-i}$, and write
\[
 s_i^\Sigma:=\mathbb E_i s_i=\tr(\Sigma_i\check G^Z A G_{-i}^Z),
 \qquad
 t_i^\Sigma:=\mathbb E_i t_i=\tr(\Sigma_iG_{-i}^Z).
\]
Inserting the reciprocal at $t_i^\Sigma$ and replacing $s_i$ by its
conditional mean in the second term gives
\begin{align}\label{eq:hs_two_level_schur_decomposition}
 b_i
 ={}&\mathbb E\left[s_i\left(
       \frac1{1+t_i^\Sigma/n}-\frac1{1+t_i/n}
     \right)\right]
     +\mathbb E\left[s_i^\Sigma\left(
       \frac1{1+\mathbb Et_i/n}-\frac1{1+t_i^\Sigma/n}
     \right)\right].
\end{align}
We apply the lemma conditionally on \(X_{-i}\) to \((s_i,t_i)\)
in the first term, and without conditioning to
\((s_i^\Sigma,t_i^\Sigma)\) in the second. Its denominator hypotheses follow
from Lemma~\ref{lem:distance_resolvent_and_ward}, since $t_i$, $t_i^\Sigma$,
and $\mathbb Et_i$ lie in $\overline{\mathcal C_Z}$.

Conditional quadratic-form control from
Lemma~\ref{lem:complex_conditional_test_reduction},
Proposition~\ref{pro:resolvent_concentration}, and Cauchy--Schwarz in the
Hilbert--Schmidt norm give, for $q\in\{1,2\}$ with $M_{\hs}^{(q)}<\infty$,
\begin{align*}
 \|s_i-s_i^\Sigma\|_{L^q}
 &=O(M_{\hs}^{(q)}),\qquad
 &\|s_i^\Sigma-\mathbb Es_i\|_{L^q}
 &=O\left(c_{n,q}|\Sigma_i\check G^Z|_{\op}M_{\hs}^{(q)}\right),\\
 \|t_i-t_i^\Sigma\|_{L^q}
 &=O(\sqrt p\,M_{\hs}^{(q)}),
 \qquad&\|t_i^\Sigma-\mathbb Et_i\|_{L^q}
 &=O\left(c_{n,q}|\Sigma|_{\hs}^\infty M_{\hs}^{(q)}\right),\\
 |s_i^\Sigma|&=O\left(|\Sigma|_{\hs}^\infty\right),
 \qquad&|\mathbb Es_i|&=O\left(\sqrt p\,|\Sigma_i\check G^Z|_{\op}\right).
\end{align*}
The bounds for \(s_i-s_i^\Sigma\) and \(t_i-t_i^\Sigma\)
also hold for their conditional \(L^q\) norms given \(X_{-i}\),
uniformly almost surely.
Using $c_{n,1}=1$ and $\sqrt p\,|\Sigma|_{\hs}^\infty/n=O(1)$, the first-moment estimate
of Lemma~\ref{lem:reciprocal_difference_estimate} bounds the two terms in
\eqref{eq:hs_two_level_schur_decomposition} by $O(M_{\hs}^{(1)})$
and $O\left(|\Sigma_i\check G^Z|_{\op}M_{\hs}^{(1)}\right)$, respectively.
When $M_{\hs}^{(2)}<\infty$,
its second-moment estimate and $c_{n,2}=n^{-1/2}$ give
\begin{align*}
 |b_i|
 &=O\Biggl((M_{\hs}^{(2)})^2\biggl[
   \frac{\sqrt p}{n}\left(1+\frac{\sqrt p\,|\Sigma|_{\hs}^\infty}{n}\right)
   +|\Sigma_i\check G^Z|_{\op}\left(
     \frac{|\Sigma|_{\hs}^\infty}{n^2}
     +\frac{\sqrt p\,(|\Sigma|_{\hs}^\infty)^2}{n^3}
   \right)\biggr]\Biggr)\\
   &=O\left(\left(1+|\Sigma_i\check G^Z|_{\op}\right)
   \frac{(M_{\hs}^{(2)})^2\sqrt p}{n}\right),
\end{align*}
since $|\Sigma|_{\hs}^\infty/n=O(p^{-1/2})$ and $p\geq1$.
Taking the better estimate and using
$\frac1n\sum_i|\Sigma_i\check G^Z|_{\op}=O(J(\check\Delta))$ yields
\[
 \frac1n\sum_i|b_i|
 =O\left(J(\check\Delta)
 \min\left\{M_{\hs}^{(1)},
             \frac{(M_{\hs}^{(2)})^2\sqrt p}{n}\right\}\right).
\]
Combining the bounds for $a_i$ and $b_i$
in the bias decomposition above gives
$|\tr(A(\mathbb E[G^Z]-\check G^Z))|=O(J(\check\Delta)\varepsilon)$.
Hilbert--Schmidt duality proves
\eqref{eq:interaction_check_bias}.

Finally, the exact check-pivot defect identity
\eqref{eq:check_pivot_exact_defect} gives
\begin{align*}
 |r_i(\check\Delta)|
 &\leq\frac{|\Sigma_i|_{\hs}}n
 \left(
 |\mathbb E[G^Z]-\check G^Z|_{\hs}
 +|\mathbb E[G_{-i}^Z-G^Z]|_{\hs}\right)
 =O\left(\frac{|\Sigma_i|_{\hs}}nJ(\check\Delta)\varepsilon\right).
\end{align*}
The last estimate uses \eqref{eq:interaction_check_bias} and
\eqref{eq:leave_one_out_hs_profile_bound} in
Lemma~\ref{lem:leave_one_out_profile_bounds}, together with
$J(\check\Delta)\geq1$.
Taking the maximum over $i$ proves
\eqref{eq:interaction_check_defect}.
\end{proof}

\begin{lemma}[Interaction bound]
\label{lem:sigma_profile_interaction}
Assume that $|Z|_{\op}=O(1)$. If $D\in\mathcal D_n(\mathcal C_Z)$
satisfies $|D|_\infty=O(1)$, then
\[
    J(D)=O(\sigma_\Sigma).
\]
\end{lemma}

To bound the contribution of a commuting approximation to $J(D)$,
we first prove the following span estimate.

\begin{lemma}[Positive matrices and span dimension]
\label{lem:positive_profile_span_bound}
Let $C_1,\ldots,C_n$ be Hermitian positive semidefinite matrices such that
$\sum_iC_i\preceq I_p$.  Then
\[
    \sum_{i=1}^n|C_i|_{\op}
    \leq
    \dim_{\mathbb R}\operatorname{span}\{C_1,\ldots,C_n\}.
\]
\end{lemma}

\begin{proof}
For each $i$, choose a unit vector $v_i$ such that
$v_i^*C_iv_i=|C_i|_{\op}$, and define the entrywise nonnegative matrix
\[
    B_{ij}:=v_j^*C_iv_j.
\]
Every column sum is at most one, so its spectral radius satisfies
\[
    \rho(B)\leq\max_{j\in[n]}\sum_{i\in[n]}|B_{ij}|\leq1.
\]
Set $r:=\dim_{\mathbb R}\operatorname{span}\{C_1,\ldots,C_n\}$. 
The rows of $B$ are the images of the $C_i$ under the real linear map
$C\mapsto(v_j^*Cv_j)_{j=1}^n$, hence $\operatorname{rank}(B)\leq r$.
Counting eigenvalues with algebraic multiplicity, at most
$\operatorname{rank}(B)$ of them are nonzero.  Therefore
\[
    \sum_{i=1}^n|C_i|_{\op}
    =\tr(B)
    \leq\sum_{\lambda\in\operatorname{Sp}(B)}|\lambda|
    \leq\operatorname{rank}(B)
    \leq r.
\]
\end{proof}

\begin{proof}[Proof of Lemma~\ref{lem:sigma_profile_interaction}]
Fix any pairwise commuting positive-semidefinite family
$T_1,\ldots,T_n$ over which the infimum in
\eqref{eq:def_sigma_profile} is taken. Write
\[
 r:=\dim_{\mathbb R}\operatorname{span}\{T_1,\ldots,T_n\},
 \qquad
 \varepsilon:=\frac1n\sum_i|\Sigma_i-T_i|_{\op}.
\]
We first estimate the interaction within this commuting family.  Set
\[
 R:=\left(I_p+\frac1n\sum_iD_iT_i\right)^{-1},
 \qquad
 C:=I_p+\frac1n\sum_i|D_i|T_i.
\]
The supporting-phase inequality gives
\[
 \mathfrak C_K(R^{-1})
 =\mathfrak C_K(1)I_p
  +\frac1n\sum_i\mathfrak C_K(D_i)T_i
 \succeq c_KC.
\]
In particular, the inverse defining \(R\) exists.
Multiplying the preceding inequality on the left and right by
\(C^{-1/2}\), then applying the inverse norm bound, gives
\[
 |C^{1/2}RC^{1/2}|_{\op}\leq c_K^{-1}.
\]
Because the $T_i$ commute, the matrices $R$, $C$, and all $T_i$ are
simultaneously diagonalizable; in particular
$C^{1/2}RC^{1/2}=CR$. Hence $|CR|_{\op}\leq c_K^{-1}$.
Moreover, the matrices $|D_i|T_iC^{-1}/n$ are positive semidefinite,
sum to $I_p-C^{-1}\preceq I_p$, and have real span dimension at most
$r$.  The span bound in Lemma~\ref{lem:positive_profile_span_bound}
therefore gives
\begin{align}\label{eq:commuting_approximation_interaction_bound}
 \frac1n\sum_i|D_i|\,|T_iR|_{\op}
 &\leq\frac1{c_K}
       \sum_i\left|\frac{|D_i|}{n}T_iC^{-1}\right|_{\op}
 \leq\frac r{c_K}.
\end{align}

To compare this interaction with that of the original profile, set
$E=\frac1n\sum_iD_i(\Sigma_i-T_i)$.  Then
$|E|_{\op}=O(\varepsilon)$, and the definitions of $R$ and $G_D^Z$ give
\[
 G_D^Z=R+R(Z+I_p-E)G_D^Z
      =R\bigl[I_p+(Z+I_p-E)G_D^Z\bigr].
\]
Splitting $\Sigma_i=T_i+(\Sigma_i-T_i)$ therefore gives
\[
 \Sigma_iG_D^Z
 =T_iR\bigl[I_p+(Z+I_p-E)G_D^Z\bigr]+(\Sigma_i-T_i)G_D^Z.
\]
Using \eqref{eq:commuting_approximation_interaction_bound},
$|G_D^Z|_{\op}=O(1)$, and $|Z|_{\op}=O(1)$, we obtain
\[
 J(D)\leq1+O\bigl(r(1+\varepsilon)+\varepsilon\bigr)
       =O\bigl((1+r)(1+\varepsilon)\bigr).
\]
The implicit constant is independent of the chosen commuting family
\(T_1,\ldots,T_n\). Taking the infimum in
\eqref{eq:def_sigma_profile} proves the assertion.
\end{proof}

\paragraph{Stability and comparison.}
If $|Z|_{\op}=O(1)$, the bound $|\check\Delta|_\infty=O(1)$ and
Lemma~\ref{lem:sigma_profile_interaction} give
$J(\check\Delta)=O(\sigma_\Sigma)$.
Combining Proposition~\ref{pro:interaction_controlled_check_pivot}
with part~\textnormal{(i)} of
Proposition~\ref{pro:two_inverse_defect_transfers}
gives the following bias comparison.

\begin{proposition}[Hilbert--Schmidt norm comparison]
\label{pro:sigma_profile_hs_comparison}
Assume that $\sqrt p\,|\Sigma|_{\hs}^\infty=O(n)$, and let
$Z\in\mathcal K$ satisfy $|Z|_{\op}=O(1)$.
For every deterministic $A\in\mathcal M_p$,
\begin{align*}
 \left|\tr\left(A(\mathbb E[G^Z]-\widetilde G^Z)\right)\right|
 &=O\left(\sigma_\Sigma|A|_{\hs}\left[
 \frac{|\Sigma|_{\hs}^\infty}n
 +\min\left\{
 M_{\hs}^{(1)},
 \frac{(M_{\hs}^{(2)})^2\sqrt p}{n}
 \right\}\right]\right).
\end{align*}
\end{proposition}

\begin{proof}
Put
\[
 \varepsilon:=\frac{|\Sigma|_{\hs}^\infty}n
 +\min\left\{
  M_{\hs}^{(1)},
  \frac{(M_{\hs}^{(2)})^2\sqrt p}{n}
 \right\}.
\]
Lemma~\ref{lem:sigma_profile_interaction} and
Proposition~\ref{pro:interaction_controlled_check_pivot} give
\[
 |\mathbb E[G^Z]-\check G^Z|_{\hs}
 =O\left(\sigma_\Sigma \varepsilon\right),
 \qquad
 |r(\check\Delta)|_\infty
 =O\left(\frac{|\Sigma|_{\hs}^\infty}n
           \sigma_\Sigma \varepsilon\right).
\]
Since $|\Sigma|_*^\infty\leq\sqrt p\,|\Sigma|_{\hs}^\infty
=O(n)$, the Hilbert--Schmidt stability estimate
\eqref{eq:inverse_defect_stability_hs} in
Proposition~\ref{pro:two_inverse_defect_transfers} applies and yields
\[
 |\check G^Z-\widetilde G^Z|_{\hs}
 =O\left(\min\{|\Sigma|_{\hs}^\infty,\sqrt p\}
              |r(\check\Delta)|_\infty\right)
    =O\left(\frac{|\Sigma|_{\hs}^\infty
                   \min\{|\Sigma|_{\hs}^\infty,\sqrt p\}}n
                  \sigma_\Sigma \varepsilon\right)
 =O\left(\sigma_\Sigma \varepsilon\right).
\]
Here $\min\{|\Sigma|_{\hs}^\infty,\sqrt p\}\leq\sqrt p$, and
\eqref{eq:hs_profile_size_condition} gives
$\sqrt p\,|\Sigma|_{\hs}^\infty/n=O(1)$.
Combining the two resolvent comparisons gives
\[
 |\mathbb E[G^Z]-\widetilde G^Z|_{\hs}
 =O\left(\sigma_\Sigma \varepsilon\right).
\]
The trace bound follows by the Hilbert--Schmidt Cauchy--Schwarz
inequality.
\end{proof}

\begin{proof}[Proof of
Theorem~\ref{the:noncommuting_limited_second_moment_profiles}]
Fix $Z\in\mathcal K$ with $|Z|_{\op}=O(1)$.
Under \eqref{eq:hs_profile_size_condition},
Propositions~\ref{pro:sigma_profile_hs_comparison}
and~\ref{pro:resolvent_concentration}, combined by the triangle inequality
in $L^q$, give for every deterministic $A\in\mathcal M_p$ and every $q\geq1$,
\begin{align}
\label{eq:matrix_hs_total_comparison}
&\left\|\tr\!\left(A(G^Z-\widetilde G^Z)\right)\right\|_{L^q}
\nonumber\\
&\quad=O\!\left(|A|_{\hs}\left[
 c_{n,q}M_{\hs}^{(q)}
 +\sigma_\Sigma\left(
  \frac{|\Sigma|_{\hs}^\infty}{n}
  +\min\!\left\{M_{\hs}^{(1)},
       \frac{(M_{\hs}^{(2)})^2\sqrt p}{n}\right\}
 \right)\right]\right).
\end{align}

Now specialize to $Z=zI_p$, so that $\widetilde G^{zI_p}=\tilde G^z$.
The additional norm bound holds uniformly for $z\in K$ by compactness.
For an arbitrary compact subset of $\Omega$, apply the estimate on the
finitely many compact pieces described at the beginning of
Section~\ref{sec:probabilistic_results} and take the maximum of the
constants.
For \(q\geq2\), \(c_{n,q}=n^{-1/2}\); taking the second entry
in the minimum gives
\eqref{eq:noncommuting_limited_second_moment_refined}.
For \(q\geq1\), taking the first entry and using
\(c_{n,q}\leq1\), \(M_{\hs}^{(1)}\leq M_{\hs}^{(q)}\), and
\(\sigma_\Sigma\geq1\) gives
\eqref{eq:noncommuting_limited_second_moment_uniform}.
\end{proof}

\begin{proof}[Proof of Corollary~\ref{cor:mp_iid_finite_variance}]
Fix a sequence in $\Theta$ with $p\to\infty$ and a compact set
$K\subset\Omega$. The finite-cover argument in
Section~\ref{sec:probabilistic_results} makes the scalar comparison
uniform on $K$. Here $\Sigma_i=I_p$ for every $i$, so
$|\Sigma|_{\hs}^\infty=\sqrt p$ and $\sigma_\Sigma\leq2$.
For every $L>0$, apply
\eqref{eq:iid_finite_variance_quadratic_form_truncation} from
Proposition~\ref{pro:quadratic_form_finite_variance} with $\sigma=1$.
Taking the supremum over $|A|_{\hs}\leq1$ and dividing by $\sqrt p$ gives
\[
    \frac{M_{\hs}^{(1)}}{\sqrt p}
    \leq\frac{\sqrt2+L}{\sqrt p}
    +2\mathbb E\!\left[
      \xi_{p,n}^2\mathbf1_{\{|\xi_{p,n}|>L\}}
    \right].
\]
Taking the limsup along the sequence and then letting \(L\to\infty\)
gives \(M_{\hs}^{(1)}/\sqrt p\to0\) by uniform integrability
of the squared entries.
Since $p=O(n)$, the profile satisfies \eqref{eq:hs_profile_size_condition}.
The first-moment estimate of
Theorem~\ref{the:noncommuting_limited_second_moment_profiles}, with
$A=I_p/p$, therefore yields
\[
 \sup_{z\in K}\|g(z)-\tilde g(z)\|_{L^1}
 =O\left(\frac{M_{\hs}^{(1)}}{\sqrt p}+\frac1n\right)
 \longrightarrow0.
\]

For the identity profile, the matrix fixed-point equation gives
$\tilde G^z=\tilde g(z)I_p$ and
\[
 \tilde g(z)
 =\left(-z+\frac1{1+(p/n)\tilde g(z)}\right)^{-1}.
\]
Equivalently,
\[
 \frac pn z\tilde g(z)^2
 +\left(z+\frac pn-1\right)\tilde g(z)+1=0.
\]

The Stieltjes solution of this equation is the transform of
\(\mu_{\mathrm{MP},p/n}\). If \(p/n\to c\in[0,\infty)\),
continuity of the Marchenko--Pastur Stieltjes transform with respect
to the aspect ratio gives
\[
    \tilde g(z)=\tilde g_{\mathrm{MP},p/n}(z)
    \longrightarrow \tilde g_{\mathrm{MP},c}(z),
    \qquad \Im z>0,
\]
where
\[
    \tilde g_{\mathrm{MP},c}(z)
    :=\int_{\mathbb R_+}\frac{\mu_{\mathrm{MP},c}(d\lambda)}{\lambda-z},
    \qquad z\in\Omega.
\]
For $c=0$, the limiting transform is $(1-z)^{-1}$, corresponding to
$\delta_1=\mu_{\mathrm{MP},0}$.

The preceding $L^1$ comparison now gives
$g(z)\to \tilde g_{\mathrm{MP},c}(z)$ in probability for every
$z$ with $\Im z>0$.  From any subsequence, diagonal extraction gives a
further subsequence with almost-sure convergence simultaneously on a
fixed countable dense subset of the upper half-plane.
The bound
\[
 |g(z)-g(w)|\leq\frac{|z-w|}{\Im z\,\Im w},
 \qquad \Im z,\Im w>0,
\]
valid for every probability Stieltjes transform, extends convergence on
this event to the entire upper half-plane.
The Stieltjes-transform continuity theorem
\cite[Theorem~B.9]{BaiSilverstein2010} gives almost-sure weak convergence
along that further subsequence.  The subsequence characterization of
convergence in probability proves the claimed weak convergence in
probability.
\end{proof}

\subsection{Operator-norm route}
\label{sse:operator_profile_route}

To display the trace dependence in the estimates below, write
\begin{align}\label{eq:def_rho_profile}
    \rho_\Sigma:=1+\frac{|\Sigma|_*^\infty}n.
\end{align}
We first prove one-class comparisons for parameters in \(\mathcal K\)
with \(|Z|_{\op}=O(1)\).
In the classwise application, we normalize each background before applying
these comparisons. We replace the classes in a fixed order, using a
conditional average pivot followed by its expectation. The last step
compares the resulting deterministic pivot with $\widetilde\Delta^Z$.

\subsubsection{One-class case}

For auxiliary models with $N$ columns, the notation $Z\in\mathcal K$
and the uniformity convention are understood over their index set of
pairs $(p,N)$. The statements below hold uniformly for every $p,N\geq1$
and every individual parameter $Z_p\in\mathcal K_p$.
No upper bound on $|Z|_{\op}$ is assumed unless stated.
For arbitrary independent real columns $x_1,\ldots,x_N$ with finite second
moments $\Sigma_i=\mathbb E[x_ix_i^\top]$, write
\[
 G=\left(\frac1N\sum_i x_ix_i^\top-Z\right)^{-1},\qquad
 \Delta_i=\left(1+\frac{x_i^\top G_{-i}x_i}{N}\right)^{-1},
\]
where $G_{-i}$ omits column $i$ while retaining normalization $1/N$.
For an auxiliary family of $N$ columns, the modulus and trace factor
are computed from that family and its normalization $1/N$. Given $q\geq1$:
\begin{align*}
 M_{\op}^{(q)}&:=\max_{i\in[N]}\sup_{|B|_{\op}\leq1}
 \|x_i^\top Bx_i-\tr(\Sigma_iB)\|_{L^q},
 \qquad \rho_\Sigma:=1+\frac{|\Sigma|_*^\infty}N.
\end{align*}
The preceding concentration and conditional replacement arguments apply
with $N$ in place of $n$, using the moduli of this auxiliary family.

The conditional versions below are understood with respect to a sigma-field
$\mathscr F$ such that $Z$ is $\mathscr F$-measurable and the selected
columns are conditionally independent given $\mathscr F$. Their profiles
are then $\Sigma_i=\mathbb E[x_ix_i^\top\mid\mathscr F]$, and their
quadratic-form moduli are defined as above using conditional $L^q$ norms.
If all stated sectorial, profile, and moment hypotheses, including
any bound on \(|Z|_{\op}\), hold almost surely with the same
deterministic constants, all expectations and \(L^q\) norms in
the definitions and conclusions are taken conditionally on
\(\mathscr F\). The common sector and phase remain deterministic and
fixed under this conditioning.

\begin{lemma}[Expected pivots and partial pivot sums]
\label{lem:operator_expected_pivots}
Every convex average $d$ of the numbers $(\mathbb E[\Delta_i])_{i\in[N]}$ belongs to
$\mathcal C_Z\setminus\{0\}$ and satisfies
\[
 |d|=O(1),\qquad |d|^{-1}=O(\rho_\Sigma).
\]
Moreover, for every deterministic $I\subset[N]$ and $q\geq1$,
\begin{align}\label{eq:partial_schur_pivot_concentration}
 \left\|\sum_{i\in I}(\Delta_i-\mathbb E\Delta_i)\right\|_{L^q}
 =O\!\left(c_{N,q}M_{\op}^{(q)}\right).
\end{align}
\end{lemma}
\begin{proof}
By \eqref{eq:def_uniform_sector_constant} and the resolvent bounds
\eqref{eq:universal_sectorial_resolvent_bounds}, Jensen's inequality gives
\[
 \mathfrak C_K(\mathbb E\Delta_i)
 \geq c_K\mathbb E\left[|\Delta_i|\right]
 \geq\frac{c_K}
 {1+\tr(\Sigma_i)/\bigl(Nc_K\eta_K\bigr)}.
\]
Together with the uniform upper pivot bound from
Lemma~\ref{lem:distance_resolvent_and_ward}, this proves the first
assertion, including convex averages.

For the second assertion, we will apply
Theorem~\ref{the:efron_stein_Lq_simple} to $\sum_{i\in I}\Delta_i$,
viewed as a function of the $N$ independent columns. It is enough to show
that replacing any one column $x_j$ by an independent copy changes this
sum by $O(M_{\op}^{(q)}/N)$ in $L^q$, uniformly in $j\in[N]$.
The Schur identity \eqref{eq:schur_pivot_trace_identity}, with $n=N$, gives
\[
 \sum_{i\in I}(1-\Delta_i)
 =\frac1N\sum_{i\in I}x_i^\top Gx_i.
\]
For a fixed $j\in[N]$, Sherman--Morrison gives the exact equality
\[
 \frac1N\sum_{i\in I}x_i^\top Gx_i
 -\frac1N\sum_{i\in I\setminus\{j\}}x_i^\top G_{-j}x_i
 =\frac{x_j^\top B_jx_j/N}{1+x_j^\top G_{-j}x_j/N},
\]
where both sums retain the normalization $1/N$ and
\[
 B_j:=\mathbf1_{\{j\in I\}}G_{-j}
 -\frac1N\sum_{i\in I\setminus\{j\}}
   G_{-j}x_ix_i^\top G_{-j}.
\]
The matrix \(B_j\) depends only on the other columns.
The positive semidefinite partial sum
\(N^{-1}\sum_{i\in I\setminus\{j\}}x_ix_i^\top\)
is bounded above by the full leave-one-out sum in the Loewner order.
The supporting-phase bound and the two inverse identities in
Lemma~\ref{lem:distance_resolvent_and_ward}, followed by
Cauchy--Schwarz, therefore give
\[
 |B_j|_{\op}=O(1),\qquad
 |\tr(\Sigma_jB_j)|
 =O\!\left(\mathfrak C_{\bar K}(\tr(\Sigma_jG_{-j}))\right).
\]
The reciprocal bounds of Lemma~\ref{lem:distance_resolvent_and_ward}
therefore imply
\[
 \frac{|\tr(\Sigma_jB_j)|}
 {N|1+\tr(\Sigma_jG_{-j})/N|}=O(1).
\]
Conditionally on the other columns, $B_j$ and $G_{-j}$ are fixed matrices
with uniformly bounded operator norms, so the definition of $M_{\op}^{(q)}$
bounds the conditional $L^q$ norms of
$x_j^\top B_jx_j-\tr(\Sigma_jB_j)$ and
$x_j^\top G_{-j}x_j-\tr(\Sigma_jG_{-j})$ by $O(M_{\op}^{(q)})$.
Equation~\eqref{eq:self_normalized_decomposition_profile}, applied
conditionally on the other columns, now bounds an independent-copy
replacement of the sum in $L^q$ by $O(M_{\op}^{(q)}/N)$.
Applying Theorem~\ref{the:efron_stein_Lq_simple} to
$\sum_{i\in I}\Delta_i$ therefore yields
\[
 \left\|\sum_{i\in I}(\Delta_i-\mathbb E\Delta_i)\right\|_{L^q}
 =O\!\left(N^{1/(q\wedge2)}\frac{M_{\op}^{(q)}}N\right)
 =O\!\left(c_{N,q}M_{\op}^{(q)}\right).
\]
All constants are uniform in $I$ and require no upper bound on
$|Z|_{\op}$.
\end{proof}

Now assume $\Sigma_i=\Sigma$ for every $i$, and set
\[
 \widehat\delta:=\frac1N\sum_{i=1}^N\mathbb E[\Delta_i],
 \qquad \widehat G:=(\widehat\delta\Sigma-Z)^{-1}.
\]
Lemma~\ref{lem:operator_expected_pivots} gives
$\widehat\delta\in\mathcal C_Z\setminus\{0\}$, so
Lemma~\ref{lem:distance_resolvent_and_ward} makes $\widehat G$ well defined.
The average is needed even within one class: equality of second moments
does not imply equality of the column laws.
For this common profile, $\rho_\Sigma=1+\tr(\Sigma)/N$.

The following comparison gives a common estimate for the first- and
second-moment cases.

\begin{proposition}[One-class comparison in a bounded sectorial background]
\label{pro:one_class_sectorial_comparison}
Assume that $\Sigma_i=\Sigma$ for every $i\in[N]$, and let
$Z\in\mathcal K$ satisfy $|Z|_{\op}=O(1)$.
For every $q\in\{1,2\}$ such that $M_{\op}^{(q)}<\infty$,
\begin{equation}
\label{eq:common_sectorial_bias}
    |\mathbb EG-\widehat G|_*
    =
    O\!\left(
        \frac{\rho_\Sigma}{N^{q-1}}
        \bigl(
            \min\{1,|\Sigma|_{\op}\}
            +M_{\op}^{(q)}
        \bigr)^q
    \right).
\end{equation}
\end{proposition}

To control the profile in the second-moment estimate, we first prove
the following spectral-spread bound.

\begin{lemma}[Spectral spread from quadratic-form fluctuations]
\label{lem:quadratic_form_spectral_spread}
Let $x\in\mathbb R^p$ have finite second moments. Then
\begin{align*}
    |\Sigma|_{\hs} 
    \leq|\Sigma|_{\op}  + \sup_{|B|_{\op}\leq1}
    \left\|x^\top Bx-\tr(\Sigma B)\right\|_{L^2},
    \qquad \Sigma:=\mathbb E[xx^\top].
\end{align*}

\end{lemma}

\begin{proof}
    The result is immediate if the supremum is infinite, and it is trivial for
$p=1$.  Otherwise, let
$\lambda_1(\Sigma)\geq\cdots\geq\lambda_p(\Sigma)\geq0$ be the
eigenvalues of $\Sigma$, choose an associated orthonormal eigenbasis
$u_1,\ldots,u_p$, and write $\xi_j:=u_j^\top x$.  Thus
\begin{align*}
    \mathbb E[\xi_j\xi_k]
    =\lambda_j(\Sigma)\mathbf1_{\{j=k\}}.
\end{align*}

For $1\leq\ell\leq\lfloor p/2\rfloor$, put $a=2\ell-1$ and
$b=2\ell$.  Consider the two centered quadratic forms
\[
    Y_-:=\xi_a^2-\xi_b^2
       -(\lambda_a(\Sigma)-\lambda_b(\Sigma)),
    \qquad
    Y_+:=2\xi_a\xi_b,
\]
associated respectively with the test matrices
$B_-:=u_au_a^\top-u_bu_b^\top$ and
$B_+:=u_au_b^\top+u_bu_a^\top$.
Their squared $L^2$ norms satisfy
\begin{align*}
    \|Y_-\|_{L^2}^2+\|Y_+\|_{L^2}^2
    &=\mathbb E\left[(\xi_a^2+\xi_b^2)^2\right]
      -(\lambda_a(\Sigma)-\lambda_b(\Sigma))^2
    \\
    &\geq
      (\lambda_a(\Sigma)+\lambda_b(\Sigma))^2
      -(\lambda_a(\Sigma)-\lambda_b(\Sigma))^2
    =4\lambda_a(\Sigma)\lambda_b(\Sigma).
\end{align*}
Choose $B_\ell\in\{B_-,-B_-,B_+,-B_+\}$ whose centered quadratic
form $Y_\ell:=x^\top B_\ell x-\tr(\Sigma B_\ell)$
has the largest $L^2$ norm, with its sign chosen so that
\[
    \mathbb E\!\left[
        Y_\ell\sum_{k<\ell}Y_k
    \right]\geq0.
\]
Set $B:=\sum_{\ell=1}^{\lfloor p/2\rfloor}B_\ell$.
The matrices $B_\ell$ act on mutually orthogonal two-dimensional
subspaces and have operator norm one, so $|B|_{\op}=1$.
Our choice of signs also gives
\begin{align}
    \left\|x^\top Bx-\tr(\Sigma B)\right\|_{L^2}^2
    &=
    \sum_{\ell=1}^{\lfloor p/2\rfloor}\|Y_\ell\|_{L^2}^2
    +
    2\sum_{\ell=1}^{\lfloor p/2\rfloor}
        \mathbb E\!\left[
            Y_\ell\sum_{k<\ell}Y_k
        \right]
    \geq
    \sum_{\ell=1}^{\lfloor p/2\rfloor}\|Y_\ell\|_{L^2}^2\nonumber\\
    &\geq
    2\sum_{\ell=1}^{\lfloor p/2\rfloor}\lambda_{2\ell-1}(\Sigma)\lambda_{2\ell}(\Sigma)
    \geq
    2\sum_{\ell=1}^{\lfloor p/2\rfloor}\lambda_{2\ell}(\Sigma)^2 
    \geq \sum_{j=2}^p\lambda_j(\Sigma)^2.
    \label{eq:quadratic_form_spectral_spread}
\end{align}
The result follows from 
\begin{align*}
    |\Sigma|_{\hs} 
    \leq |\lambda_1(\Sigma)u_1u_1^\top|_{\hs}  + 
    |\Sigma-\lambda_1(\Sigma) u_1u_1^\top|_{\hs}
    \leq|\Sigma|_{\op}  + 
    \sqrt{\sum_{j=2}^p\lambda_j(\Sigma)^2}.
\end{align*}
\end{proof}
\begin{remark}[Second-moment size and quadratic-form fluctuations]
The centered quadratic-form modulus does not control the size of the
second-moment matrices.  Indeed, if \(x_i=\varepsilon_i v_i\), where
\(v_i\in\mathbb R^p\) is deterministic and \(\varepsilon_i\)
is a Rademacher variable, then
$\Sigma_i=v_iv_i^{\top}$ and
\[
    x_i^{\top}Bx_i-\tr(\Sigma_iB)=0
\]
almost surely for every $B$, whereas
$|\Sigma_i|_{\hs}=|v_i|^2$ may be arbitrarily large. This shows that the
largest eigenvalue cannot be included in the spectral sum in
\eqref{eq:quadratic_form_spectral_spread}; profile sizes must therefore
remain separate from the quadratic-form moduli.
\end{remark}
\begin{proof}[Proof of Proposition~\ref{pro:one_class_sectorial_comparison}]
If $\Sigma=0$, every column vanishes almost surely, and the conclusion
is immediate. Assume henceforth that $\Sigma\neq0$.
All estimates involving $M_{\op}^{(2)}$ below are used only when this
quantity is finite.
Fix a deterministic matrix $A$ with $|A|_{\op}\leq1$, and set
\[
    \delta:=\frac1N\sum_{i=1}^N\Delta_i,
    \qquad
    s_i:=x_i^\top\widehat GAG_{-i}x_i.
\]
Let $\mathbb E_i$ denote expectation conditional on all columns except
$x_i$. Since $\mathbb E_i s_i=\tr(\Sigma\widehat GAG_{-i})$ and
$\mathbb E[\delta-\widehat\delta]=0$, 
the resolvent identity and $Gx_i=\Delta_iG_{-i}x_i$ give
\begin{align}
    \tr\!\left(A(\mathbb EG-\widehat G)\right)
    &={}
    \widehat\delta\,
    \tr\!\left(\Sigma\widehat GA\,\mathbb EG\right)
    -\frac1N\sum_{i=1}^N\mathbb E[\Delta_i s_i]
    \notag\\
    &={}
    \mathbb E\!\left[
        \delta\,\tr(\Sigma\widehat GAG)
    \right]
    -\mathbb E\!\left[
        (\delta-\widehat\delta)
        \tr\!\left(\Sigma\widehat GA(G-\mathbb EG)\right)
    \right]
    \notag\\
    &\quad-
    \frac1N\sum_{i=1}^N
    \mathbb E\!\left[
        \Delta_i\tr(\Sigma\widehat GAG_{-i})
    \right]
    -\frac1N\sum_{i=1}^N
    \mathbb E\!\left[
        (\Delta_i-\mathbb E_i\Delta_i)
        (s_i-\mathbb E_i s_i)
    \right]
    \notag\\
    &=a+b+c,
\end{align}
where
\begin{itemize}
    \item $\displaystyle
        a:=-\mathbb E\!\left[
            (\delta-\widehat\delta)
            \tr\!\left(\Sigma\widehat GA(G-\mathbb EG)\right)
        \right]$;
    \item $\displaystyle
        b:=-\frac1N\sum_{i=1}^N
        \mathbb E\!\left[
            (\Delta_i-\mathbb E_i\Delta_i)
            (s_i-\mathbb E_i s_i)
        \right]$;
    \item $\displaystyle
        c:=\frac1N\sum_{i=1}^N
        \mathbb E\!\left[
            \Delta_i\,
            \tr\!\left(\Sigma\widehat GA(G-G_{-i})\right)
        \right]$.
\end{itemize}

We first bound $|\Sigma\widehat G|_{\op}$.
The identity $\widehat\delta\Sigma\widehat G=I_p+Z\widehat G$ 
together with $|\widehat\delta|^{-1}=O(\rho_\Sigma)$,
$|Z|_{\op}=O(1)$, and the uniform resolvent bounds, gives
\begin{equation}
    |\Sigma\widehat G|_{\op}
    =O\!\left(\min\{|\Sigma|_{\op},\rho_\Sigma\}\right).
\end{equation}
We now bound $a$, $b$, and $c$ in turn.
\begin{enumerate}
    \item \emph{Bound on $a$.}
    The uniform pivot bounds and
    Lemma~\ref{lem:operator_expected_pivots} give
    \[
        |\delta-\widehat\delta|=O(1),
        \qquad
        \|\delta-\widehat\delta\|_{L^2}
        =O\!\left(\frac{M_{\op}^{(2)}}{N\sqrt N}\right).
    \]
    For $q\in\{1,2\}$ with $M_{\op}^{(q)}<\infty$,
    Proposition~\ref{pro:resolvent_concentration} gives
    \[
        \left\|
            \tr\!\left(\Sigma\widehat GA(G-\mathbb EG)\right)
        \right\|_{L^q}
        =O\!\left(
            |\Sigma\widehat G|_{\op}
            N^{1/q-1}M_{\op}^{(q)}
        \right).
    \]
    Using the boundedness of $\delta-\widehat\delta$ for the first
    estimate and Cauchy--Schwarz for the second, we obtain
    \begin{align*}
        |a|
        &=O\!\left(
            \min\{|\Sigma|_{\op},\rho_\Sigma\}M_{\op}^{(1)}
        \right)
        \qquad \text{and}\qquad
        |a|
        =O\!\left(
            \frac{
                \min\{|\Sigma|_{\op},\rho_\Sigma\}
                (M_{\op}^{(2)})^2
            }{N^2}
        \right).
    \end{align*}

    \item \emph{Bound on $b$.}
    The uniform pivot bounds and the reciprocal Lipschitz estimate give
    \[
        |\Delta_i-\mathbb E_i\Delta_i|=O(1),
        \qquad
        \|\Delta_i-\mathbb E_i\Delta_i\|_{L^2}
        =O\!\left(\frac{M_{\op}^{(2)}}N\right).
    \]
    Since $\widehat GAG_{-i}$ is independent of $x_i$ and has uniformly
    bounded operator norm, conditional quadratic-form control gives
    \[
        \|s_i-\mathbb E_i s_i\|_{L^q}=O(M_{\op}^{(q)}),
        \qquad q\in\{1,2\},\quad M_{\op}^{(q)}<\infty.
    \]
    Therefore
    \[
        |b|=O(M_{\op}^{(1)})
        \qquad \text{and}\qquad
        |b|=O\!\left(\frac{(M_{\op}^{(2)})^2}{N}\right).
    \]

    \item \emph{Bound on $c$.}
    Set $T:=G+GZG=N^{-1}GXX^\top G$.
    The Sherman--Morrison identities give
    \[
        \sum_{i=1}^N\Delta_i(G-G_{-i})
        =-\frac1N\sum_{i=1}^N
            \Delta_i^2G_{-i}x_ix_i^\top G_{-i}
        =-T.
    \]
    Therefore $c=-\frac1N\tr\!\left(\Sigma\widehat GA\,\mathbb ET\right)$.
    The uniform resolvent bounds and $|Z|_{\op}=O(1)$ imply
    $|T|_{\op}=O(1)$. Since $\operatorname{rank}(T)\leq N$,
    $|T|_*=O(N)$, and hence
    \[
        |c|=O(|\Sigma\widehat G|_{\op})
        =O\!\left(\min\{|\Sigma|_{\op},\rho_\Sigma\}\right).
    \]

    For the second-moment bound, apply
    Lemma~\ref{lem:uncentered_quadratic_form_bound} conditionally on the
    other columns. Together with the uniform pivot and resolvent bounds,
    this gives
    \begin{align}
        \label{eq:one_class_c_profile_bound}
        |c|
        \leq\!\frac{O(1)}{N^2}\sum_{i=1}^N
            \mathbb E\left[
                \left|x_i^\top G_{-i}\Sigma\widehat GAG_{-i}x_i\right|
            \right]
        \!=O\!\left(\frac{|\Sigma|_{\hs}^2}{N}\right)
        \!=O\!\left(
            \frac{|\Sigma|_{\op}^2+(M_{\op}^{(2)})^2}{N}
        \right),
    \end{align}
    where the last bound follows from
    Lemma~\ref{lem:quadratic_form_spectral_spread}.

    To refine this bound, choose a unit eigenvector $v$ associated
    with the largest eigenvalue of $\Sigma$, set $P:=vv^\top$, and write
    $\Sigma=|\Sigma|_{\op}P+\Sigma_0$. Set
    \[
        G_0:=\left(
            \frac1N(I_p-P)XX^\top(I_p-P)-Z
        \right)^{-1},
        \qquad
        T_0:=G_0+G_0ZG_0.
    \]
    The common second-moment matrix of the projected columns $(I_p-P)x_i$ is $\Sigma_0$. Applying the
    same conditional quadratic-form estimate to these columns gives
    \[
        \left|\tr\!\left(\Sigma_0\widehat GA\,\mathbb ET_0\right)\right|
        =O\!\left(|\Sigma_0|_{\hs}^2\right)
        =O\!\left((M_{\op}^{(2)})^2\right),
    \]
    thanks to Lemma~\ref{lem:quadratic_form_spectral_spread}. Moreover,
    \[
        \operatorname{rank}\!\left(
            XX^\top-(I_p-P)XX^\top(I_p-P)
        \right)\leq2.
    \]
    The resolvent identity and the uniform operator-norm bounds
    therefore imply $|G-G_0|_*=O(1)$. Since
    \[
        T-T_0=(G-G_0)+(G-G_0)ZG+G_0Z(G-G_0),
    \]
    we also have $|T-T_0|_*=O(1)$.
    Finally, $P\Sigma\widehat G$ has rank at most one and
    $\Sigma_0=(I_p-P)\Sigma$. Splitting $\Sigma=P\Sigma+\Sigma_0$, and
    then replacing $T$ by $T_0+(T-T_0)$ in the $\Sigma_0$ term, gives
    \begin{align*}
        N|c|
        &\leq
        |P\Sigma\widehat G|_*\,|\mathbb ET|_{\op}
        +\left|\tr\!\left(\Sigma_0\widehat GA\,\mathbb ET_0\right)\right|
        +|\Sigma_0\widehat G|_{\op}\,\mathbb E|T-T_0|_*
        \\
        &\leq
        O(|\Sigma\widehat G|_{\op})
        +O\!\left((M_{\op}^{(2)})^2\right)
        +O(|\Sigma\widehat G|_{\op})
        =O\!\left(
            \min\{|\Sigma|_{\op},\rho_\Sigma\}
            +(M_{\op}^{(2)})^2
        \right).
    \end{align*}
    Combining this with \eqref{eq:one_class_c_profile_bound}, we obtain
    \begin{equation}
    \label{eq:one_class_c_refined_bound}
        |c|=O\!\left(
            \frac{
                \min\{|\Sigma|_{\op}^2,|\Sigma|_{\op},\rho_\Sigma\}
                +(M_{\op}^{(2)})^2
            }{N}
        \right).
    \end{equation}
\end{enumerate}

For $q=1$, the preceding estimates and $\rho_\Sigma\geq1$ give
\begin{align*}
    |a|+|b|+|c|
    &=O\!\left(
        \min\{|\Sigma|_{\op},\rho_\Sigma\}
        +\bigl(1+\min\{|\Sigma|_{\op},\rho_\Sigma\}\bigr)
            M_{\op}^{(1)}
    \right)
    \\
    &=O\!\left(
        \rho_\Sigma
        \bigl(\min\{1,|\Sigma|_{\op}\}+M_{\op}^{(1)}\bigr)
    \right).
\end{align*}

For $q=2$, use \eqref{eq:one_class_c_refined_bound} and
\[
    \min\{|\Sigma|_{\op}^2,|\Sigma|_{\op},\rho_\Sigma\}
    \leq\rho_\Sigma\min\{1,|\Sigma|_{\op}^2\},
    \qquad
    1+\frac{\min\{|\Sigma|_{\op},\rho_\Sigma\}}N
    \leq2\rho_\Sigma.
\]

It follows that
\[
    |a|+|b|+|c|
    =O\!\left(
        \frac{\rho_\Sigma}{N}
        \bigl(\min\{1,|\Sigma|_{\op}^2\}+(M_{\op}^{(2)})^2\bigr)
    \right),
\]
which implies \eqref{eq:common_sectorial_bias} for $q=2$.
All bounds are uniform over deterministic $A$ with $|A|_{\op}\leq1$,
so nuclear duality completes the proof. 
\end{proof}

\subsubsection{Multiple classes}

Return to the original \(n\)-column family and the classes
\(I_1,\ldots,I_k\) introduced in the main-results section.  Thus
\[
    \mathbb E[x_ix_i^\top]=\Sigma^{(\ell)}
    \quad\text{for }i\in I_{\ell},
    \qquad
    n_{\ell}:=|I_{\ell}|,
    \qquad
    \pi_{\ell}:=\frac{n_{\ell}}{n},
    \qquad
    S_{\ell}:=\frac1n\sum_{i\in I_{\ell}}x_ix_i^\top .
\]
Fix \(Z\in\mathcal K\) and a deterministic ordering
\(I_1,\ldots,I_k\) of the classes. Put
\[
    \mathcal F^{\ell}
    :=
    \sigma(x_i:i\in I_{\ell+1}\cup\cdots\cup I_k),
    \qquad
    Z^{(0)}:=Z,
    \qquad
    G^{(0)}:=G^Z.
\]
At each stage, we first average the Schur pivots conditionally on
the unreplaced classes. We then take the expectation of this
conditional average and retain the resulting deterministic
coefficient at subsequent stages.
Given \(Z^{(\ell-1)}\), define, for \(i\in I_{\ell}\),
\[
    G_{-i}^{(\ell-1)}
    :=
    \left(
        \sum_{b\geq \ell}S_b-\frac{x_ix_i^\top}{n}
        -Z^{(\ell-1)}
    \right)^{-1},
\]
and
\begin{align}
\label{eq:class_stage_pivots}
    \widehat\delta_{\mathcal F^{\ell}}^{(\ell)}
    &:=
    \frac1{n_{\ell}}\sum_{i\in I_{\ell}}
    \mathbb E\!\left[
        \left(
            1+\frac{x_i^\top G_{-i}^{(\ell-1)}x_i}{n}
        \right)^{-1}
        \middle|\mathcal F^{\ell}
    \right],
    &
    \widehat\delta^{(\ell)}
    &:=
    \mathbb E\!\left[
        \widehat\delta_{\mathcal F^{\ell}}^{(\ell)}
    \right].
\end{align}
Define the conditional and deterministic parameters by
\begin{align}
\label{eq:class_stage_parameters}
    Z_{\mathcal F^{\ell}}^{(\ell)}
    &:=
    Z-\sum_{b<\ell}\pi_b\widehat\delta^{(b)}\Sigma^{(b)}
    -\pi_{\ell}\widehat\delta_{\mathcal F^{\ell}}^{(\ell)}\Sigma^{(\ell)},
    \qquad
    Z^{(\ell)}
    =
    Z-\sum_{b\leq \ell}\pi_b\widehat\delta^{(b)}\Sigma^{(b)}.
\end{align}
The corresponding resolvents are
\begin{align}
\label{eq:class_stage_resolvents}
    G_{\mathcal F^{\ell}}^{(\ell)}
    &:=
    \left(
        \sum_{b>\ell}S_b-Z_{\mathcal F^{\ell}}^{(\ell)}
    \right)^{-1},
    \qquad
    G^{(\ell)}
    :=
    \left(
        \sum_{b>\ell}S_b-Z^{(\ell)}
    \right)^{-1}.
\end{align}

The parameters \(Z^{(\ell)}\) are deterministic.
Both \(G_{\mathcal F^{\ell}}^{(\ell)}\) and \(G^{(\ell)}\) are
\(\mathcal F^{\ell}\)-measurable and may depend on the unreplaced classes.
In particular: 
\[
    G^{(k)}=(-Z^{(k)})^{-1}=G_{\widehat\Delta}^Z
    \qquad \text{where}\quad \widehat\Delta_i:=\widehat\delta^{(a)}
    \qquad\text{for }i\in I_a.
\]

\begin{lemma}[Iterated pivots]
\label{lem:class_average_pivot_estimates}
For every \(\ell\in[k]\), \(Z^{(\ell)}\in \mathcal K_p\) and
\(Z_{\mathcal F^{\ell}}^{(\ell)}\in \mathcal K_p\).
All stage resolvents are well defined. Moreover, $\widehat\delta_{\mathcal F^{\ell}}^{(\ell)},
    \widehat\delta^{(\ell)}
    \in\mathcal C_Z\setminus\{0\}$, and
\[
    \left|
        \widehat\delta_{\mathcal F^{\ell}}^{(\ell)}
    \right|
    +
    \left|
        \widehat\delta^{(\ell)}
    \right|
    =
    O(1),
    \qquad
    \left|
        \widehat\delta_{\mathcal F^{\ell}}^{(\ell)}
    \right|^{-1}
    +
    \left|
        \widehat\delta^{(\ell)}
    \right|^{-1}
    =
    O(\rho_\Sigma).
\]
In addition,
\begin{align}
\label{eq:class_pivot_first_moment}
    n_{\ell}
    \left\|
        \widehat\delta_{\mathcal F^{\ell}}^{(\ell)}
        -
        \widehat\delta^{(\ell)}
    \right\|_{L^1}
    &=
    O\!\left(
        \pi_{\ell}\tau_\Sigma M_{\op}^{(1)}
    \right),
    \\
\label{eq:class_pivot_second_moment}
    n_{\ell}
    \left\|
        \widehat\delta_{\mathcal F^{\ell}}^{(\ell)}
        -
        \widehat\delta^{(\ell)}
    \right\|_{L^2}
    &=
    O\!\left(
        \pi_{\ell}
        \left[
            \frac{\tau_\Sigma M_{\op}^{(2)}}{\sqrt n}
            +
            \frac{(M_{\op}^{(2)})^2}{n}
        \right]
    \right).
\end{align}
If \(Z=Z^*\prec0\), all stage resolvents are positive definite and all
stage pivots are positive.
\end{lemma}

\begin{proof}
Starting from \(Z^{(0)}=Z\), proceed by induction with the fixed common
phase. If \(Y=Z-\sum_j d_jS_j\), where
\(d_j\in\mathcal C_Z\) and \(S_j\succeq0\), then
\(\mathcal C_Y\subseteq\mathcal C_Z\subseteq\mathcal C_K\) and
\begin{align*}
    \mathfrak C_K(-Y)
    &=\mathfrak C_K(-Z)
      +\sum_j\mathfrak C_K(d_j)S_j
    \succeq c_K
      \left(\eta_KI_p+\sum_j|d_j|S_j\right).
\end{align*}
Consequently, \(Y\in\mathcal K_p\). This applies to the retained pivot
terms and to the unreplaced sample blocks, for which \(d_j=1\).
At stage \(\ell\), condition on \(\mathcal F^{\ell}\) and apply
Lemmas~\ref{lem:distance_resolvent_and_ward}
and~\ref{lem:operator_expected_pivots} with \(N=n\),
zero columns outside \(I_{\ell}\), and parameter
\(Z^{(\ell-1)}-\sum_{b>\ell}S_b\).
This gives the pivot inclusions and bounds and ensures that the
stage resolvents are well defined.
The definitions \eqref{eq:class_stage_parameters} and convexity of
\(\mathcal C_Z\) give the sector inclusions and preserve
membership in \(\mathcal K_p\) for \(Z^{(\ell)}\) and
\(Z_{\mathcal F^{\ell}}^{(\ell)}\);
positivity is preserved when \(Z=Z^*\prec0\).

It remains to prove \eqref{eq:class_pivot_first_moment} and
\eqref{eq:class_pivot_second_moment}. Set
\[
    \Delta_{i,\ell}
    :=
    \left(
        1+\frac{x_i^\top G_{-i}^{(\ell-1)}x_i}{n}
    \right)^{-1},
    \qquad i\in I_{\ell}.
\]
Since
\[
    n_{\ell}\widehat\delta_{\mathcal F^{\ell}}^{(\ell)}
    =
    \mathbb E\!\left[
        \sum_{i\in I_{\ell}}\Delta_{i,\ell}
        \middle|\mathcal F^{\ell}
    \right],
    \qquad
    n_{\ell}\widehat\delta^{(\ell)}
    =
    \mathbb E\!\left[
        \sum_{i\in I_{\ell}}\Delta_{i,\ell}
    \right],
\]
apply \eqref{eq:partial_schur_pivot_concentration} with \(N=n\),
\(I=I_{\ell}\), parameter \(Z^{(\ell-1)}\), and zero columns in the
already replaced classes. Conditional Jensen then gives
\begin{align}
\label{eq:class_pivot_coarse_bound}
    n_{\ell}
    \left\|
        \widehat\delta_{\mathcal F^{\ell}}^{(\ell)}
        -
        \widehat\delta^{(\ell)}
    \right\|_{L^q}
    =
    O\!\left(
        c_{n,q}M_{\op}^{(q)}
    \right),
    \qquad q=1,2.
\end{align}

To obtain a bound involving \(|\Sigma^{(\ell)}|_{\op}\), define
\[
    \widehat\delta_\Sigma^{(\ell)}
    :=
    \frac1{n_{\ell}}\sum_{i\in I_{\ell}}
    \mathbb E\!\left[
        \left(
            1+\frac1n
            \tr(\Sigma^{(\ell)}G_{-i}^{(\ell-1)})
        \right)^{-1}
        \middle|\mathcal F^{\ell}
    \right].
\]
This quantity remains \(\mathcal F^{\ell}\)-measurable: only the
quadratic form inside the reciprocal has been replaced by its
conditional mean.
Apply both estimates of Lemma~\ref{lem:reciprocal_difference_estimate}
with \(s=1\) and \(t=x_i^\top G_{-i}^{(\ell-1)}x_i\), conditionally on all
columns other than \(x_i\). The conditional mean of \(t\) is
\(\tr(\Sigma^{(\ell)}G_{-i}^{(\ell-1)})\), and
Lemma~\ref{lem:distance_resolvent_and_ward} bounds the reciprocal
denominators. Conditional quadratic-form control from
Lemma~\ref{lem:complex_conditional_test_reduction}, after taking expectations
given \(\mathcal F^{\ell}\) and averaging over \(i\in I_{\ell}\), gives, uniformly in
the conditional background,
\[
    \left|
        \widehat\delta_{\mathcal F^{\ell}}^{(\ell)}
        -
        \widehat\delta_\Sigma^{(\ell)}
    \right|
    =
    O\!\left(
        \min\left\{
            \frac{M_{\op}^{(1)}}{n},
            \frac{(M_{\op}^{(2)})^2}{n^2}
        \right\}
    \right).
\]
Moreover, combine the reciprocal Lipschitz bound
\eqref{eq:reciprocal_lipschitz_bound} in the proof of
Lemma~\ref{lem:reciprocal_difference_estimate} with
Proposition~\ref{pro:resolvent_concentration} applied to
\(\tr(\Sigma^{(\ell)}G_{-i}^{(\ell-1)})\) in the same auxiliary \(n\)-column model.
Conditional Jensen then gives
\[
    \left\|
        \widehat\delta_\Sigma^{(\ell)}
        -
        \mathbb E[
            \widehat\delta_\Sigma^{(\ell)}
        ]
    \right\|_{L^q}
    =
    O\!\left(
        \frac{
            c_{n,q}
            |\Sigma^{(\ell)}|_{\op}
            M_{\op}^{(q)}
        }{n}
    \right),
    \qquad q=1,2.
\]
Centering
\(\widehat\delta_{\mathcal F^{\ell}}^{(\ell)}
-\widehat\delta_\Sigma^{(\ell)}\)
and combining the two preceding bounds by the triangle inequality,
with \(n_{\ell}=n\pi_{\ell}\), gives
\begin{align*}
    n_{\ell}
    \left\|
        \widehat\delta_{\mathcal F^{\ell}}^{(\ell)}
        -
        \widehat\delta^{(\ell)}
    \right\|_{L^1}
    &=
    O\!\left(
        \pi_{\ell}
        (1+|\Sigma^{(\ell)}|_{\op})
        M_{\op}^{(1)}
    \right),
    \\
    n_{\ell}
    \left\|
        \widehat\delta_{\mathcal F^{\ell}}^{(\ell)}
        -
        \widehat\delta^{(\ell)}
    \right\|_{L^2}
    &=
    O\!\left(
        \pi_{\ell}|\Sigma^{(\ell)}|_{\op}
        \frac{M_{\op}^{(2)}}{\sqrt n}
        +
        \pi_{\ell}\frac{(M_{\op}^{(2)})^2}{n}
    \right).
\end{align*}
Taking the better of these estimates and
\eqref{eq:class_pivot_coarse_bound}, and using
\[
    \min\{1,\pi_{\ell}(1+|\Sigma^{(\ell)}|_{\op})\}
    \leq
    2\pi_{\ell}\tau_\Sigma,
    \qquad
    \min\{1,\pi_{\ell}|\Sigma^{(\ell)}|_{\op}\}
    \leq
    2\pi_{\ell}\tau_\Sigma,
\]
proves \eqref{eq:class_pivot_first_moment} and
\eqref{eq:class_pivot_second_moment}.
\end{proof}

Throughout the three results below, set
\begin{equation}
\label{eq:classwise_comparison_error}
    \varepsilon
    :=\min\left\{
        \tau_\Sigma\bigl(1+M_{\op}^{(1)}\bigr),
        \frac{\tau_\Sigma M_{\op}^{(2)}}{\sqrt n}
        +\frac{\tau_\Sigma^2+(M_{\op}^{(2)})^2}{n}
    \right\}.
\end{equation}
The second entry is interpreted as $+\infty$ when
$M_{\op}^{(2)}=\infty$. In these three results, write
\[
    \widehat G:=G_{\widehat\Delta}^Z=G^{(k)}.
\]

\begin{proposition}
\label{pro:classwise_averaged_pivot_comparison}
Let $Z\in\mathcal K$ satisfy $|Z|_{\op}=O(1)$. Then
\begin{equation}
    |\mathbb E[G^Z]-\widehat G|_*
    =O\!\left(\rho_\Sigma\varepsilon\right).
\end{equation}
\end{proposition}
\begin{proof}
The proof relies on the telescoping decomposition:
\begin{align}\label{eq:def_Rr}
    \widehat G - \mathbb E[G^Z]
    = \sum_{\ell=1}^k \mathbb{E}[R_{\ell}]
    \qquad \text{with}\quad 
   R_{\ell} := G^{(\ell)}-\mathbb E[G^{(\ell-1)}\mid\mathcal F^{\ell}].
\end{align}
For each $\ell\in[k]$, $G^{(\ell)}$ and $G^{(\ell-1)}$ both share the $\mathcal F^{\ell}$-measurable parameter matrix
\begin{align*}
    B_{\ell}:=\sum_{b>\ell}S_b-Z^{(\ell-1)}.
\end{align*}
Indeed, the definitions of the
stage parameters \eqref{eq:class_stage_parameters} and
resolvents \eqref{eq:class_stage_resolvents} give\begin{align}
    G^{(\ell-1)}&=(S_{\ell}+B_{\ell})^{-1},
    \qquad 
    G^{(\ell)}
    =\bigl(\pi_{\ell}\widehat\delta^{(\ell)}\Sigma^{(\ell)}+B_{\ell}\bigr)^{-1}.
\end{align}
    The background $B_{\ell}$ need not be uniformly bounded. Set
\begin{align*}
    U_{\ell}:=\left(
        I_p+\sum_{b>\ell}S_b
        +\sum_{b<\ell}\pi_b|\widehat\delta^{(b)}|\Sigma^{(b)}
    \right)^{-1/2}.
\end{align*}
We will prove, uniformly in $\ell$, that
$\mathbb E\left\vert U_{\ell}^{-1}R_{\ell}U_{\ell}^{-1}\right\vert_*
    =O\!\left(\rho_\Sigma\pi_{\ell}\varepsilon\right)$,
and then sum over the classes.

Fix $\ell\in[k]$ and condition on $\mathcal F^{\ell}$. 
The matrix $U_{\ell}$ is real symmetric, $\mathcal F^{\ell}$-measurable, and
satisfies $|U_{\ell}|_{\op}\leq1$. First note that:
\[
    \sum _{b>\ell}U_{\ell}S_bU_{\ell} +\sum_{b<\ell}
    \pi_b|\widehat\delta^{(b)}|U_{\ell}\Sigma^{(b)}U_{\ell}
    =I_p-U_{\ell}^2\preceq I_p.
\]
Hence, by Cauchy--Schwarz,
\begin{equation}
\label{eq:class_background_phase_bound}
    \sup_{|\theta_b|\leq1}
    \left|
        \sum_{b>\ell}\theta_bU_{\ell}S_bU_{\ell}
        +\sum_{b<\ell}\theta_b\pi_b\widehat\delta^{(b)}
            U_{\ell}\Sigma^{(b)}U_{\ell}
    \right|_{\op}
    \leq1.
\end{equation}
In particular,
\[
    |U_{\ell}B_{\ell}U_{\ell}|_{\op}\leq1+|Z|_{\op}=O(1).
\]
Moreover, $-U_{\ell}B_{\ell}U_{\ell}\in\mathcal K_p$ and therefore satisfies the conditional one-class comparison with uniform deterministic constants. Indeed $\mathcal C_{-U_{\ell}B_{\ell}U_{\ell}}
    \subseteq\mathcal C_Z\subseteq\mathcal C_K$ and
\begin{align}
    \label{eq:use_of_eta_leq_1}
    \mathfrak C_K(U_{\ell}B_{\ell}U_{\ell})
    =U_{\ell}\mathfrak C_K(B_{\ell})U_{\ell}
    &\succeq c_KU_{\ell}\left(
        \sum_{b>\ell}S_b
        +\sum_{b<\ell}\pi_b|\widehat\delta^{(b)}|\Sigma^{(b)}
        +\eta_KI_p
      \right)U_{\ell} \nonumber\\
    &=c_K\bigl[I_p-(1-\eta_K)U_{\ell}^2\bigr]
      \succeq c_K\eta_KI_p.
\end{align}
The last inequality uses $0<\eta_K\leq1$, as defined in
\eqref{eq:def_family_parameter_bounds}.

Apply Proposition~\ref{pro:one_class_sectorial_comparison} conditionally to
\[
    y_i:=\sqrt{\pi_{\ell}}\,U_{\ell}x_i,
    \qquad i\in I_{\ell},
\]
with $N=n_{\ell}$ and parameter $-U_{\ell}B_{\ell}U_{\ell}$. With normalization $1/n_{\ell}$,
these columns give the sample matrix $U_{\ell}S_{\ell}U_{\ell}$ and the same
conditional average pivot $\widehat\delta_{\mathcal F^{\ell}}^{(\ell)}$.
The random and comparison resolvents are, respectively,
\[
    U_{\ell}^{-1}G^{(\ell-1)}U_{\ell}^{-1},
    \qquad
    U_{\ell}^{-1}G_{\mathcal F^{\ell}}^{(\ell)}U_{\ell}^{-1}.
\]
For these columns, the common conditional second moment,
quadratic-form moduli, and trace factor satisfy $|\pi_{\ell}U_{\ell}\Sigma^{(\ell)}U_{\ell}|_{\op}
    \leq\pi_{\ell}|\Sigma^{(\ell)}|_{\op},$ and
\begin{align*}
    M_{\op}^{(q)}(y_i:i\in I_{\ell}\mid\mathcal F^{\ell})
    &\leq\pi_{\ell}M_{\op}^{(q)},
    \qquad 
    1+\frac{\tr(\pi_{\ell}U_{\ell}\Sigma^{(\ell)}U_{\ell})}{n_{\ell}}
    \leq\rho_\Sigma.
\end{align*}
Using $\min\{1,\pi_{\ell}|\Sigma^{(\ell)}|_{\op}\}
    \leq2\pi_{\ell}\tau_\Sigma$ and $n_{\ell}=\pi_{\ell} n$, 
the unified estimate \eqref{eq:common_sectorial_bias} gives,
for $q\in\{1,2\}$ such that
$M_{\op}^{(q)}<\infty$,
\begin{align*}
    \left|
        U_{\ell}^{-1}\left(
            \mathbb E[G^{(\ell-1)}\mid\mathcal F^{\ell}]
            -G_{\mathcal F^{\ell}}^{(\ell)}
        \right)U_{\ell}^{-1}
    \right|_*
    &=O\!\left(
        \frac{\rho_\Sigma\pi_{\ell}}{n^{q-1}}
        \bigl(\tau_\Sigma+M_{\op}^{(q)}\bigr)^q
    \right).
\end{align*}
Since $\tau_\Sigma\geq1$, the $q=1$ bound is
$O(\rho_\Sigma\pi_{\ell}\tau_\Sigma(1+M_{\op}^{(1)}))$.
When $M_{\op}^{(2)}<\infty$, the $q=2$ bound is
$O(\rho_\Sigma\pi_{\ell}(\tau_\Sigma^2+(M_{\op}^{(2)})^2)/n)$.
Taking the better of the available estimates therefore gives
\begin{equation}
    \left|
        U_{\ell}^{-1}\left(
            \mathbb E[G^{(\ell-1)}\mid\mathcal F^{\ell}]
            -G_{\mathcal F^{\ell}}^{(\ell)}
        \right)U_{\ell}^{-1}
    \right|_*
    =O\!\left(\rho_\Sigma\pi_{\ell}\varepsilon\right)
\end{equation}
almost surely, with uniform deterministic constants.

The resolvent identity gives
\begin{equation*}
    G_{\mathcal F^{\ell}}^{(\ell)}-G^{(\ell)}
    =\pi_{\ell}\bigl(
        \widehat\delta^{(\ell)}-\widehat\delta_{\mathcal F^{\ell}}^{(\ell)}
    \bigr)G_{\mathcal F^{\ell}}^{(\ell)}\Sigma^{(\ell)}G^{(\ell)}.
\end{equation*}
The uniform bounds from Lemma~\ref{lem:distance_resolvent_and_ward},
applied with the common phase to the normalized resolvents, imply,
pointwise,
\[
    \left|
        U_{\ell}^{-1}\bigl(G_{\mathcal F^{\ell}}^{(\ell)}-G^{(\ell)}\bigr)U_{\ell}^{-1}
    \right|_*
    \leq C(\rho_\Sigma-1)n_{\ell}
    \left|\widehat\delta_{\mathcal F^{\ell}}^{(\ell)}
        -\widehat\delta^{(\ell)}\right|.
\]
Taking expectations and applying
\eqref{eq:class_pivot_first_moment} or
\eqref{eq:class_pivot_second_moment} from
Lemma~\ref{lem:class_average_pivot_estimates}, with
$\|\cdot\|_{L^1}\leq\|\cdot\|_{L^2}$ for the second-moment bound, yields finally
\begin{align}
    \label{eq:normalized_stage_replacement}
    \mathbb E\left|
        U_{\ell}^{-1}R_{\ell}U_{\ell}^{-1}
    \right|_*
    =O\!\left(\rho_\Sigma\pi_{\ell}\varepsilon\right).
\end{align}
One can then conclude with \eqref{eq:def_Rr} since $|U_{\ell}|_{\op}\leq1$ and $\sum_{\ell}\pi_{\ell}=1$.
\end{proof}

For stability, we need to bound the weighted inverse
defect $\sum_i|\widehat\Delta_i r_i(\widehat\Delta)|$ of Proposition~\ref{pro:two_inverse_defect_transfers}, rather than only the resolvent difference
or an unweighted bound on $r(\widehat\Delta)$.

\begin{lemma}[Weighted inverse defect]
\label{lem:classwise_weighted_inverse_defect}
Let $Z\in\mathcal K$ satisfy $|Z|_{\op}=O(1)$. Then
\begin{equation}
    \sum_{i=1}^n|\widehat\Delta_i r_i(\widehat\Delta)|
    =O\!\left(\rho_\Sigma\varepsilon\right).
\end{equation}
\end{lemma}

\begin{proof}
For $i\in I_a$, the inverse defect is constant on the class; write
$r^{(a)}(\widehat\Delta):=r_i(\widehat\Delta)$.
Starting from the definition given in~\eqref{eq:def_inverse_defect},
we have
\begin{align}\label{eq:first_identity_r_D_class}
    \sum_{i=1}^n|\widehat\Delta_i r_i(\widehat\Delta)|
    &=\sum_{a=1}^k n_a
        |\widehat\delta^{(a)}r^{(a)}(\widehat\Delta)|
    =\sum_{a=1}^k
        \left|
            n_a(1-\widehat\delta^{(a)})
            -\pi_a\widehat\delta^{(a)}
                \tr(\Sigma^{(a)}\widehat G)
        \right|.
\end{align}
The Schur identity \eqref{eq:schur_pivot_trace_identity}, applied
at stage $a$ with the pivots \eqref{eq:class_stage_pivots}, gives
\[
    n_a(1-\widehat\delta^{(a)})
    =\mathbb E[\tr(S_aG^{(a-1)})].
\]
Moreover, with the definition of $R_{\ell}$ in \eqref{eq:def_Rr}:
\[
    \widehat G
    =\mathbb E[G^{(a)}] - \sum_{\ell>a}\mathbb E[G^{(\ell-1)}-G^{(\ell)}]
    =\mathbb E[G^{(a)}] + \sum_{\ell>a}\mathbb E[R_{\ell}].
\]
Finally, taking traces in the identities
\[
    G^{(\ell-1)}(S_{\ell}+B_{\ell})=I_p
    =G^{(\ell)}
        \bigl(\pi_{\ell}\widehat\delta^{(\ell)}\Sigma^{(\ell)}+B_{\ell}\bigr),
\]
and using the $\mathcal F^{\ell}$-measurability of $B_{\ell}$ and $G^{(\ell)}$,
we obtain:
\begin{align*}
    \mathbb E[\tr(S_{\ell}G^{(\ell-1)})]
        -\pi_{\ell}\widehat\delta^{(\ell)}
            \tr\!\left(\Sigma^{(\ell)}\mathbb E[G^{(\ell)}]\right)
    = \mathbb E[\tr(B_{\ell}R_{\ell})].
\end{align*}
Therefore,~\eqref{eq:first_identity_r_D_class} gives
\begin{align*}
    \sum_{i=1}^n|\widehat\Delta_i r_i(\widehat\Delta)|
    &=\sum_{a=1}^k
        \left|
            \mathbb E[\tr(B_aR_a)]-\pi_a\widehat\delta^{(a)}
                \sum_{\ell>a}\tr(\Sigma^{(a)}\mathbb E[R_{\ell}])
        \right|\\
    &\leq
        \sum_{\ell=1}^k\left|\mathbb E[\tr(B_{\ell}R_{\ell})]\right|
        +\sum_{\ell=1}^k\sum_{a<\ell}
            \pi_a|\widehat\delta^{(a)}|
            |\tr(\Sigma^{(a)}\mathbb E[R_{\ell}])|\\
    &\leq \sum_{\ell=1}^k
        \mathbb E\!\left[
            \left|U_{\ell}^{-1}R_{\ell}U_{\ell}^{-1}\right|_*
            \left(
                |U_{\ell}B_{\ell}U_{\ell}|_{\op}
                +\left|\sum_{a<\ell}
                    \theta_{a,\ell}\pi_a\widehat\delta^{(a)}
                    U_{\ell}\Sigma^{(a)}U_{\ell}
                \right|_{\op}
            \right)
        \right],
\end{align*}
where, for each $\ell$, we choose deterministic
unit phases $\theta_{a,\ell}$, $a<\ell$, such that
\[
    \theta_{a,\ell}\widehat\delta^{(a)}
        \tr(\Sigma^{(a)}\mathbb E[R_{\ell}])\in\mathbb R_+.
\]
By~\eqref{eq:class_background_phase_bound},
\[
    \left|
        \sum_{a<\ell}
            \theta_{a,\ell}\pi_a\widehat\delta^{(a)}
            U_{\ell}\Sigma^{(a)}U_{\ell}
    \right|_{\op}\leq1.
\]
The result follows from~\eqref{eq:normalized_stage_replacement}, $|U_{\ell}B_{\ell}U_{\ell}|_{\op}=O(1)$, and $\sum_{\ell}\pi_{\ell}=1$.
\end{proof}

\begin{proposition}[Classwise comparison with the fixed-point resolvent]
\label{pro:classwise_operator_profile_comparison}
Let $Z\in\mathcal K$ satisfy $|Z|_{\op}=O(1)$. Then
\begin{equation}
\label{eq:tau_operator_bias}
    |\mathbb E[G^Z]-\widetilde G^Z|_*
    =O\!\left(\rho_\Sigma^{3/2}\varepsilon\right).
\end{equation}
If $Z=Z^*\prec0$, one can replace $\rho_\Sigma^{3/2}$ by $\rho_\Sigma$ in this bound.
\end{proposition}

\begin{proof}
Lemma~\ref{lem:class_average_pivot_estimates} gives
$\widehat\Delta\in\mathcal D_n(\mathcal C_Z\setminus\{0\})$, $|\widehat\Delta|_\infty=O(1)$ and $\max_i|\widehat\Delta_i|^{-1}=O(\rho_\Sigma)$.
The fixed-point equation \eqref{eq:def_matrix_parameter_resolvent}
and Proposition~\ref{pro:matrix_parameter_fixed_point}
also give $|\widetilde\Delta|_\infty=O(1)$ and
\[
    |\widetilde\Delta_i|^{-1}
    =\left|1+\frac1n\tr(\Sigma_i\widetilde G^Z)\right|
    \leq1+\frac{\tr(\Sigma_i)}n|\widetilde G^Z|_{\op}
    =O(\rho_\Sigma).
\]
Therefore \eqref{eq:sectorial_weighted_defects_bounded_pivots} in
Proposition~\ref{pro:two_inverse_defect_transfers}, followed by
Lemma~\ref{lem:classwise_weighted_inverse_defect}, yields
\[
    |\widehat G-\widetilde G^Z|_*
    =O\!\left(
        \rho_\Sigma^{1/2}\sum_{i=1}^n|\widehat\Delta_i r_i(\widehat\Delta)|
    \right)
    =O\!\left(\rho_\Sigma^{3/2}\varepsilon\right).
\]
Combining this with
Proposition~\ref{pro:classwise_averaged_pivot_comparison} and
$\rho_\Sigma\geq1$ proves \eqref{eq:tau_operator_bias}.

If $Z=Z^*\prec0$, the pivots are positive by
Lemma~\ref{lem:class_average_pivot_estimates} and one can use 
Corollary~\ref{cor:negative_hermitian_inverse_defect_transfer}.
\end{proof}

\begin{proof}[Proof of Theorem~\ref{the:operator_profile_comparison}]
For every deterministic \(A\in\mathcal M_p\), every \(q\geq1\), and
every \(Z\in\mathcal K\),
\begin{align*}
    \left\|
        \tr\!\left(
            A(G^Z-\widetilde G^Z)
        \right)
    \right\|_{L^q}
    &\leq
    \left\|
        \tr\!\left(
            A(G^Z-\mathbb E[G^Z])
        \right)
    \right\|_{L^q}
    +
    |A|_{\op}
    |\mathbb E[G^Z]-\widetilde G^Z|_*.
\end{align*}
Assume now that \(|Z|_{\op}=O(1)\).
Proposition~\ref{pro:resolvent_concentration} and
\eqref{eq:tau_operator_bias} then give
\begin{align}
\label{eq:matrix_operator_total_comparison}
&\left\|
    \tr\!\left(
        A(G^Z-\widetilde G^Z)
    \right)
\right\|_{L^q}
\nonumber\\
&\quad=
O\!\Biggl(
    |A|_{\op}
    \Biggl[
        c_{n,q}M_{\op}^{(q)}
        +
        \rho_\Sigma^{3/2}
        \min\left\{
            \tau_\Sigma\bigl(1+M_{\op}^{(1)}\bigr),
            \,
            \frac{\tau_\Sigma M_{\op}^{(2)}}{\sqrt n}
            +
            \frac{
                \tau_\Sigma^2+(M_{\op}^{(2)})^2
            }{n}
        \right\}
    \Biggr]
\Biggr).
\end{align}
The minimum uses the same $+\infty$ convention as
\eqref{eq:classwise_comparison_error} when
$M_{\op}^{(2)}=\infty$.
Specialize \eqref{eq:matrix_operator_total_comparison} to
\(Z=zI_p\). The additional norm bound holds uniformly for $z\in K$
by compactness. 
When \(z<0\), Proposition~\ref{pro:classwise_operator_profile_comparison} allows us to replace the factor
\(\rho_\Sigma^{3/2}\) by \(\rho_\Sigma\).
\end{proof}

\subsection{Resolvent sandwiches}
\label{sec:linear_response_superoperator}

Wagner, Couillet, Debbah, and Slock
\cite{WagnerCouilletDebbahSlock2012} derive two-resolvent equivalents
by differentiating a deterministic equivalent with respect to the
parameter of an additive matrix perturbation.
Hoydis, ten Brink, and Debbah \cite{HOY13} state a normalized-trace
equivalent for uniformly bounded positive semidefinite insertions.
In a random-feature model, Louart, Liao, and Couillet
\cite{LOU17b} obtain operator-norm error bounds for the expected
sandwich by a direct resolvent expansion.
Xiao and Slock \cite{XiaoSlock2026} formulate the differentiation
argument in a matrix-Dyson-equation setting using a one-resolvent
approximation uniform under perturbation.

Here Cauchy's integral formula and Minkowski's inequality transfer
the preceding one-resolvent \(L^q\) estimates to resolvent sandwiches,
with explicit dependence on the dimension and second-moment profile.

We retain the compact set $K$ and the constant $\eta_K$ from
\eqref{eq:def_family_parameter_bounds}. For $z\in K$, we have
$\eta_z\geq\eta_K>0$, and compactness gives
$\sup_{z\in K}|z|=O(1)$.
The finite-cover reduction following
\eqref{eq:def_family_parameter_bounds} extends the estimates to
arbitrary compact subsets of $\Omega$.

For a deterministic matrix \(A\in\mathcal M_p\) and \(t\in\mathbb C\),
define \(Z_t\) and, whenever the inverse exists, \(G^{Z_t}\) by
\[
    Z_t:=zI_p-tA,
    \qquad
    G^{Z_t}
    =\left(\frac1nXX^\top+tA-zI_p\right)^{-1}.
\]
Note that
\begin{align}\label{eq:linear_response_path_separation}
    d\bigl(\mathcal W(Z_t),\mathbb R_+\bigr)
    \geq\eta_z-|t|\,|A|_{\op}.
\end{align}
Now choose $r:=\frac{\eta_K}{4}$. Uniformly for $z\in K$, $|t|<2r$, and $|A|_{\op}\leq1$,
\eqref{eq:linear_response_path_separation} gives
\[
    d\bigl(\mathcal W(Z_t),\mathbb R_+\bigr)\geq\eta_K/2,
    \qquad
    |Z_t|_{\op}\leq\sup_{z\in K}|z|+2r=O(1).
\]
Set $K^+:=K+\{w\in\mathbb C:|w|\leq2r\}$. Then
$d(\operatorname{conv}K^+,\mathbb R_+)\geq\eta_K/2$ and
$\mathcal W(Z_t)\subseteq K^+$, so the preceding estimates apply with
$K^+$ in place of $K$. Proposition~\ref{pro:matrix_parameter_fixed_point} gives
\begin{align}\label{eq:linear_response_deterministic_deformation}
 \widetilde G^{Z_t}
 =\left(
 tA-zI_p+\frac1n\sum_{i=1}^n
 \frac{\Sigma_i}
 {1+\frac1n\tr(\Sigma_i\widetilde G^{Z_t})}
 \right)^{-1}.
\end{align}
Proposition~\ref{pro:matrix_parameter_analyticity} shows that
$t\mapsto\widetilde G^{Z_t}$ is holomorphic on $|t|<2r$.

\begin{proposition}[Deterministic equivalent of the resolvent sandwich]
\label{pro:deterministic_linear_response_equation}
For every $z\in\Omega$ and every deterministic matrix $A\in\mathcal M_p$,
the linear fixed-point equation
\eqref{eq:deterministic_linear_response_equation} has a unique solution, and
\begin{align}
    \widetilde{\mathcal S}^{\,z}[A]
    =
    -\left.
    \frac{\mathrm d}{\mathrm dt}
    \widetilde G^{Z_t}
    \right|_{t=0}.
\end{align}
The solution depends linearly on $A$. Moreover, $\widetilde{\mathcal S}^{\,\bar z}[\overline A]
    =\overline{\widetilde{\mathcal S}^{\,z}[A]}$.
\end{proposition}

\begin{proof}
If $A=0$, the path $Z_t=zI_p$ is constant. Otherwise apply the preceding
disk construction to $A/|A|_{\op}$ and the scalar parameter
$s=t|A|_{\op}$. This gives holomorphic dependence in a neighborhood of
$t=0$ for every $A\in\mathcal M_p$. Put
\[
 T:=-\left.\partial_t\widetilde G^{Z_t}\right|_{t=0}.
\]
Differentiating \eqref{eq:linear_response_deterministic_deformation}
gives
\[
 T=\tilde G^zA\tilde G^z
 +\frac1{n^2}\sum_i(\widetilde\Delta_i^{zI_p})^2
   \tr(\Sigma_iT)\tilde G^z\Sigma_i\tilde G^z.
\]
Since
$\widetilde\Delta_i^{zI_p}=(1+\frac1n\tr(\Sigma_i\tilde G^z))^{-1}$,
this is \eqref{eq:deterministic_linear_response_equation}.

For uniqueness, let $U$ solve the homogeneous equation and set
\[
 s_i:=(\widetilde\Delta_i^{zI_p})^2\tr(\Sigma_iU).
\]
Taking the weighted traces in that equation gives
\[
 s=\Psi^{zI_p}(\widetilde\Delta^{zI_p},
                 \widetilde\Delta^{zI_p})s.
\]
Proposition~\ref{pro:fixed_point_secant_contraction} implies $s=0$.
The homogeneous equation then gives $U=0$, proving uniqueness and
$T=\widetilde{\mathcal S}^{\,z}[A]$.
Linearity in $A$ and the stated conjugation identity follow immediately
from uniqueness of the solution.
\end{proof}

\begin{proof}[Proof of
Corollary~\ref{cor:resolvent_sandwich_comparison}]
The case $A=0$ is immediate. We first prove the estimates for
$|A|_{\op}\leq1$.
For a deterministic test matrix \(B\), set
\[
 F(t):=\tr\!\left(B[G^{Z_t}-\widetilde G^{Z_t}]\right).
\]
The inverse defining $G^{Z_t}$ and
Proposition~\ref{pro:matrix_parameter_analyticity} show that $F$ is
holomorphic pathwise on $|t|<2r$. On the integration circle, the resolvent
bounds in Lemma~\ref{lem:distance_resolvent_and_ward} dominate $|F(t)|$ by a deterministic multiple of
$|B|_*$, so Cauchy's integral formula and Minkowski's inequality apply.  They give
\[
 \|F'(0)\|_{L^q}
 \leq\frac1{2\pi r}\int_0^{2\pi}
       \|F(re^{\ii\theta})\|_{L^q}\,d\theta
 \leq\frac1r\sup_{|t|=r}\|F(t)\|_{L^q}.
\]
Differentiating the random resolvent and applying
Proposition~\ref{pro:deterministic_linear_response_equation} gives
\[
 -\left.\partial_tG^{Z_t}\right|_{t=0}=G^zAG^z,
 \qquad
 -\left.\partial_t\widetilde G^{Z_t}\right|_{t=0}
 =\widetilde{\mathcal S}^{\,z}[A],
\]
so we obtain the common transfer bound
\begin{equation}
 \left\|\tr\!\left(
 B[\mathcal S^z[A]-\widetilde{\mathcal S}^{\,z}[A]]
 \right)\right\|_{L^q}
    \leq\frac1r\sup_{|t|=r}
 \left\|\tr\!\left(B[G^{Z_t}-\widetilde G^{Z_t}]\right)\right\|_{L^q}.
\end{equation}

Apply \eqref{eq:matrix_operator_total_comparison} and
\eqref{eq:matrix_hs_total_comparison} with test matrix $B$, uniformly on
the contour. The same simplifications used in the
proofs of the one-resolvent theorems give the two stated bounds, since
$1/r=4/\eta_K=O(1)$. The general operator norm comparison is used even when
$z<0$, because the contour is not a family of negative Hermitian
parameters.

For general $A\neq0$, apply the preceding estimates to
$A/|A|_{\op}$. Both $\mathcal S^z[A]$ and
$\widetilde{\mathcal S}^{\,z}[A]$ depend linearly on $A$, so multiplication
by $|A|_{\op}$ gives the stated bounds.

\end{proof}

\clearpage
\appendix

\section{Notation index}
\label{sec:specific_notations_appendix}

This appendix indexes the principal notation used in the paper.

{\small
\setlength{\tabcolsep}{4pt}
\renewcommand{\arraystretch}{1.15}
\begin{longtable}{@{}
>{\raggedright\arraybackslash}p{0.25\textwidth}
>{\raggedright\arraybackslash}p{0.46\textwidth}
>{\raggedright\arraybackslash}p{0.19\textwidth}
@{}}
\hline
\textbf{Symbol} & \textbf{Meaning} & \textbf{Definition} \\
\hline
\endhead
\hline
\endfoot
$M_{\hs}^{(q)},M_{\op}^{(q)}$
& Quadratic-form moduli for the Hilbert--Schmidt and operator test balls
& \eqref{eq:def_M_hs_op_main} \\
$I_a,n_a,\pi_a$
& Exact second-moment classes, their sizes, and their proportions
& Section~\ref{sssec:main_operator_route} \\
$\tau_\Sigma$
& Operator-size-adjusted exact-class coefficient
& \eqref{eq:def_tau_profile} \\
$|\Sigma|_\nu^\infty$
& Maximum of $|\Sigma_i|_\nu$, for $\nu\in\{\op,\hs,*\}$
& \eqref{eq:def_profile_hs_infty} \\
$\rho_\Sigma$
& Trace factor $1+\frac{|\Sigma|_*^\infty}{n}$
& \eqref{eq:def_rho_profile} \\
$\sigma_\Sigma$
& Commuting-approximation coefficient of the finite profile
& \eqref{eq:def_sigma_profile} \\
$c_{n,q}$
& $n^{1/(q\wedge2)-1}$
& \eqref{eq:def_c_nq} \\
$\Omega$
& Scalar parameter domain $\mathbb C\setminus[0,\infty)$
& Section~\ref{sec:deterministic_study} \\
$\mathcal W(Z)$
& Numerical range $\{u^*Zu:|u|=1\}$
& Section~\ref{sec:deterministic_study} \\
$\eta_Z,\eta_z$
& Distances $d(\mathcal W(Z),\mathbb R_+)$
  and $d(z,\mathbb R_+)$, respectively
& Section~\ref{sec:deterministic_study} \\
$\mathcal C_Z$
& $\cone(\mathbb R_+\cup-\mathcal W(Z))$, closed for admissible $Z$
& \eqref{eq:def_matrix_resolvent_sector} \\
$\theta_Z^-,\theta_Z^+$
& Boundary arguments of $\mathcal C_Z$
& Lemma~\ref{lem:matrix_parameter_sector_aperture} \\
$\theta_Z,\omega_Z,c_Z$
& Bisecting angle and phase; cosine of half the aperture
& \eqref{eq:matrix_parameter_angular_phase} \\
$\mathfrak C_Z(w),\mathfrak C_{\bar Z}(w)$
& Sector parts of $w$ relative to the sector and its conjugate
& \eqref{eq:def_weighted_real_part} \\
$K$
& Compact scalar set determining the common sector; in the probabilistic proof setup, $d(\operatorname{conv}K,\mathbb R_+)>0$
& \eqref{eq:def_family_parameter_bounds} \\
$\eta_K$
& Convex-hull separation capped at one: $1\wedge d(\operatorname{conv}K,\mathbb R_+)$
& \eqref{eq:def_family_parameter_bounds} \\
$\mathcal C_K,\overline{\mathcal C_K}$
& Common closed sector $\operatorname{cone}(\mathbb R_+\cup(-K))$ and its conjugate
& \eqref{eq:def_common_probabilistic_sector} \\
$\omega_K,c_K$
& Bisecting unit phase and supporting constant of $\mathcal C_K$
& \eqref{eq:def_uniform_sector_constant} \\
$\mathfrak C_K(w),\mathfrak C_{\bar K}(w)$
& Sector parts of $w$ relative to the common sector and its conjugate
& \eqref{eq:def_common_probabilistic_real_parts} \\
$\mathcal K_p,\mathcal K$
& Matrix class in dimension $p$ and corresponding parameter families over $\Theta$; no upper-norm restriction
& \eqref{eq:universal_sectorial_base_class} \\
$G^Z,G_{-i}^Z$
& Random resolvent and its $i$-th leave-one-out version
& Section~\ref{sse:universal_probabilistic_toolbox} \\
$G_D^Z,\Phi^Z(D)$
& Deterministic resolvent at a diagonal profile and its fixed-point map
& \eqref{eq:def_matrix_parameter_resolvent} \\
$\Psi^Z(D,D')$
& Deterministic secant kernel
& Lemma~\ref{lem:weighted_secant_estimate} \\
$\widetilde\Delta^Z,\widetilde G^Z$
& Fixed point of $\Phi^Z$ and its deterministic resolvent
& Proposition~\ref{pro:matrix_parameter_fixed_point} \\
$\Delta_i^Z$
& Random Schur pivot
& \eqref{eq:def_matrix_schur_pivot} \\
$\check\Delta$
& Diagonal profile with entries $(\mathbb E[\Delta_i^{-1}])^{-1}$
& Section~\ref{sse:universal_profile_pivots} \\
$\widehat\delta_{\mathcal F^{\ell}}^{(\ell)},\widehat\delta^{(\ell)}$
& Conditional average Schur pivot for class $\ell$ and its expectation
& \eqref{eq:class_stage_pivots} \\
$Z_{\mathcal F^{\ell}}^{(\ell)},Z^{(\ell)}$
& Class-replacement parameters using the conditional and expected pivots
& \eqref{eq:class_stage_parameters} \\
$G_{\mathcal F^{\ell}}^{(\ell)},G^{(\ell)}$
& Class-replacement resolvents using the conditional and expected pivots
& \eqref{eq:class_stage_resolvents} \\
$r(D)$
& Inverse defect of a diagonal pivot
& \eqref{eq:def_inverse_defect} \\
$J(D)$
& One plus the average weighted profile interaction
& \eqref{eq:def_matrix_check_pivot_interaction} \\
$\tilde\Gamma^z$
& Companion diagonal profile for the scalar parameter $z$,
  $\tilde\Gamma^z=-z^{-1}\widetilde\Delta^{zI_p}$
& Corollary~\ref{cor:tilde_Gamma_analytic} \\
$\mathcal S^z[A],\widetilde{\mathcal S}^{\,z}[A]$
& Random resolvent sandwich and its deterministic equivalent
& Section~\ref{sssec:main_resolvent_sandwiches} \\
\end{longtable}}

\section{Auxiliary results}
\label{app:auxiliary_results}

This appendix collects the moment inequalities, quadratic-form estimates,
and Stieltjes-transform facts used in the proofs and applications.

For the replacement-moment inequality, fix $q\geq1$.
Let $U_1,\ldots,U_n$ be independent
random variables, and let
\[
    F=f(U_1,\ldots,U_n)\in L^q
\]
be real- or complex-valued.  For each $i\in[n]$, let $U_i'$ be an
independent copy of $U_i$, independent of everything else, and define
\[
    F^{(i)}
    :=f(U_1,\ldots,U_{i-1},U_i',U_{i+1},\ldots,U_n).
\]
Thus $F^{(i)}$ has the same law as $F$ and belongs to $L^q$.

\begin{theorem}[Martingale replacement-moment inequality]
\label{the:efron_stein_Lq_simple}
If $1\leq q\leq2$, then
\[
    \|F-\mathbb E[F]\|_{L^q}
    = O\!\left(
      \left(\sum_{i=1}^n\|F-F^{(i)}\|_{L^q}^q\right)^{1/q}
    \right)
    = O\!\left(
      n^{1/q}\sup_{i\in[n]}\|F-F^{(i)}\|_{L^q}
    \right).
\]
If $q\geq2$, then
\[
    \|F-\mathbb E[F]\|_{L^q}
    = O\!\left(
      \left(\sum_{i=1}^n\|F-F^{(i)}\|_{L^q}^2\right)^{1/2}
    \right)
    = O\!\left(
      \sqrt n\sup_{i\in[n]}\|F-F^{(i)}\|_{L^q}
    \right).
\]
\end{theorem}

\begin{proof}
We first assume that $F$ is real-valued. Set
\[
    \mathcal F_i:=\sigma(U_1,\ldots,U_i),
    \qquad
    M_i:=\mathbb E[F\mid\mathcal F_i],
    \qquad
    D_i:=M_i-M_{i-1},
    \qquad
    V_i:=F-F^{(i)},
\]
with $\mathcal F_0$ trivial.  Then $(D_i)_{i=1}^n$ is a martingale
difference sequence and
\[
    F-\mathbb E[F]=\sum_{i=1}^nD_i.
\]
Since $U_i'$ has the same law as $U_i$, and since
$U_i',U_{i+1},\ldots,U_n$ are independent of $\mathcal F_i$,
\[
    \mathbb E[F^{(i)}\mid\mathcal F_i]
    =\mathbb E[F\mid\mathcal F_{i-1}]
    =M_{i-1}.
\]
Therefore
\[
    D_i
    =\mathbb E[F-F^{(i)}\mid\mathcal F_i]
    =\mathbb E[V_i\mid\mathcal F_i],
\]
and conditional Jensen's inequality gives, for every $q\geq1$,
\[
    \|D_i\|_{L^q}\leq\|V_i\|_{L^q}.
\]

For $1<q\leq2$, set $S_0:=0$ and $S_i:=\sum_{j=1}^iD_j$.
The scalar inequality underlying the von Bahr--Esseen estimate
\cite{vonBahr-Esseen-65} is
\[
    |x+y|^q
    \leq |x|^q+q|x|^{q-1}\operatorname{sgn}(x)y+C_q|y|^q,
    \qquad x,y\in\mathbb R.
\]
Apply it with $x=S_{i-1}$ and $y=D_i$. The linear term is integrable by
H\"older's inequality and has conditional mean zero given
$\mathcal F_{i-1}$. Thus
\[
    \mathbb E[|S_i|^q\mid\mathcal F_{i-1}]
    \leq |S_{i-1}|^q
       +C_q\mathbb E[|D_i|^q\mid\mathcal F_{i-1}].
\]
Taking expectations and iterating from $S_0=0$ gives
\[
    \mathbb E\left[\left|\sum_{i=1}^nD_i\right|^q\right]
    \leq C_q\sum_{i=1}^n\mathbb E[|D_i|^q].
\]
For $q=1$, the corresponding bound follows directly from the triangle
inequality. Hence, for $1\leq q\leq2$,
\[
    \|F-\mathbb E[F]\|_{L^q}
    = O\!\left(
      \left(\sum_{i=1}^n\|V_i\|_{L^q}^q\right)^{1/q}
    \right).
\]

For $q\geq2$, Burkholder's square-function inequality
\cite{Burkholder1988Sharp} gives
\[
    \left\|\sum_{i=1}^nD_i\right\|_{L^q}
    \leq(q-1)
    \left\|\left(\sum_{i=1}^nD_i^2\right)^{1/2}\right\|_{L^q}.
\]
Moreover, Minkowski's inequality in $L^{q/2}$ gives
\[
    \left\|\left(\sum_{i=1}^nD_i^2\right)^{1/2}\right\|_{L^q}^2
    =\left\|\sum_{i=1}^nD_i^2\right\|_{L^{q/2}}
    \leq\sum_{i=1}^n\|D_i\|_{L^q}^2.
\]
It follows that
\[
    \|F-\mathbb E[F]\|_{L^q}
    \leq(q-1)
    \left(\sum_{i=1}^n\|V_i\|_{L^q}^2\right)^{1/2}
    \leq(q-1)\sqrt n\sup_{i\in[n]}\|V_i\|_{L^q}.
\]
If $F$ is complex-valued, apply the real result separately to
$\Re(F)$ and $\Im(F)$.  Since the corresponding replacement differences
have $L^q$ norm at most $\|V_i\|_{L^q}$, the triangle inequality
introduces at most a factor of two.
\end{proof}

\begin{lemma}[Reciprocal difference estimate]
\label{lem:reciprocal_difference_estimate}
Let $s$ and $t$ be integrable complex random variables, with $s,t\in L^2$
for the second estimate. If
\[
    \left|1+\frac{\mathbb E[t]}n\right|^{-1}
    =O(1),
    \qquad
    \left|1+\frac tn\right|^{-1}
    =O(1)
\]
almost surely, then
\begin{align*}
    \left|
      \mathbb E\left[
        s\left(
          \frac1{1+\mathbb E[t]/n}
          -
          \frac1{1+t/n}
        \right)
      \right]
    \right|
    &=
    O\left(
       \|s-\mathbb E[s]\|_{L^1}
       +\frac{|\mathbb E[s]|}{n}
        \|t-\mathbb E[t]\|_{L^1}
    \right),
    \\
    \left|
      \mathbb E\left[
        s\left(
          \frac1{1+\mathbb E[t]/n}
          -
          \frac1{1+t/n}
        \right)
      \right]
    \right|
    &=
    O\left(
       \frac{\|s-\mathbb E[s]\|_{L^2}
             \|t-\mathbb E[t]\|_{L^2}}{n}
       +\frac{|\mathbb E[s]|
             \|t-\mathbb E[t]\|_{L^2}^2}{n^2}
    \right).
\end{align*}
\end{lemma}

\begin{proof}
Put
\[
    d:=
    \frac1{1+\mathbb E[t]/n}
    -
    \frac1{1+t/n}.
\]
The exact expansions
\begin{align}
    d
    &=
    \frac{t-\mathbb E[t]}
    {n(1+\mathbb E[t]/n)(1+t/n)}
    \label{eq:reciprocal_first_order_identity}
    \\
    &=
    \frac{t-\mathbb E[t]}
    {n(1+\mathbb E[t]/n)^2}
    -
    \frac{(t-\mathbb E[t])^2}
    {n^2(1+\mathbb E[t]/n)^2(1+t/n)}
    \label{eq:reciprocal_second_order_identity}
\end{align}
and the denominator bounds give
\begin{equation}\label{eq:reciprocal_lipschitz_bound}
    |d|=O(1),
    \qquad
    |d|=O\left(\frac{|t-\mathbb E[t]|}{n}\right).
\end{equation}
Taking expectations in \eqref{eq:reciprocal_first_order_identity} gives
\[
    |\mathbb E[d]|
    =O\left(\frac{\|t-\mathbb E[t]\|_{L^1}}{n}\right).
\]
If $t\in L^2$, the linear term in
\eqref{eq:reciprocal_second_order_identity} has mean zero, so
\[
    |\mathbb E[d]|
    =O\left(\frac{\|t-\mathbb E[t]\|_{L^2}^2}{n^2}\right).
\]
Finally,
\[
    \mathbb E[sd]
    =\mathbb E\left[(s-\mathbb E[s])d\right]
     +\mathbb E[s]\,\mathbb E[d].
\]
The first conclusion follows from $|d|=O(1)$ and the first bound on
$|\mathbb E[d]|$. For the second, apply Cauchy--Schwarz to the centered
term and use the second bounds on $|d|$ and $|\mathbb E[d]|$.
\end{proof}

For \(B\in\mathcal M_p\), write
\(B_{\mathrm s}:=(B+B^\top)/2\).

\begin{lemma}[Complex and conditional test reduction]
\label{lem:complex_conditional_test_reduction}
Let $x\in\mathbb R^p$ have finite second moments, and set $\Sigma:=\mathbb E[xx^\top]$. Then  \begin{align*}
    x^{\top}Bx-\tr(\Sigma B)
    &=x^{\top}B_{\mathrm s}x-\tr(\Sigma B_{\mathrm s}),
    \qquad 
    |B_{\mathrm s}|_{\op}&\leq|B|_{\op},
    \qquad
    |B_{\mathrm s}|_{\hs}&\leq|B|_{\hs}.
\end{align*}
\end{lemma}
Consequently, each modulus in \eqref{eq:def_M_hs_op_main} may be applied
to a complex, nonsymmetric test after symmetrization.  If a bound has only
been established for real symmetric tests, applying it to the real and
imaginary parts gives the complex version with at most twice the numerical
constant.  Finally, if \(B\) is independent of \(x\) and lies almost surely
in the relevant deterministic test ball, the same \(L^q\) bound
follows by conditioning on \(B\) and then taking expectations.

\begin{proof}
The skew-transpose part of $B$ has zero quadratic form against a real
vector, and its trace against the real symmetric matrix $\Sigma$ is zero.
The two norm inequalities follow from the triangle inequality and the
invariance of the operator and Hilbert--Schmidt norms under transpose.
Write $B_{\mathrm s}=C+\ii D$ with $C,D$ real symmetric.  The real and
imaginary parts of the centered quadratic form are the centered quadratic
forms associated with $C$ and $D$, respectively; moreover
$|C|_{\op},|D|_{\op}\leq|B_{\mathrm s}|_{\op}$, and the same inequalities
hold with the Hilbert--Schmidt norm in place of the operator norm.  
\end{proof}

\begin{lemma}[Uncentered quadratic-form bound]
\label{lem:uncentered_quadratic_form_bound}
For every $i\in[n]$ and every deterministic $B\in\mathcal M_p$,
\begin{align}
    \mathbb E\left[|x_i^\top Bx_i|\right]
    &\leq|\Sigma_i|_{\hs}|B|_{\hs},
    \qquad \mathbb E\left[|x_i^\top Bx_i|\right]
    \leq|B|_{\op}\tr(\Sigma_i).
\end{align}
Consequently,
\[
    M_{\hs}^{(1)}\leq2|\Sigma|_{\hs}^\infty,
    \qquad M_{\op}^{(1)}\leq2|\Sigma|_*^\infty.
\]
\end{lemma}

\begin{proof}
Polar decomposition and Cauchy--Schwarz give the pointwise inequality
\[
    |x_i^*Bx_i|
    \leq
    \bigl(x_i^*(BB^*)^{1/2}x_i\bigr)^{1/2}
    \bigl(x_i^*(B^*B)^{1/2}x_i\bigr)^{1/2}.
\]
Since $x_i$ is real, $x_i^*Bx_i=x_i^\top Bx_i$. Taking expectations
and applying Cauchy--Schwarz, followed by Hilbert--Schmidt duality, gives
\[
    \mathbb E\left[|x_i^\top Bx_i|\right]
    \leq\sqrt{\tr\!\left(\Sigma_i(BB^*)^{1/2}\right)
               \tr\!\left(\Sigma_i(B^*B)^{1/2}\right)}
    \leq|\Sigma_i|_{\hs}|B|_{\hs}.
\]
Here both positive square-root matrices have Hilbert--Schmidt norm
$|B|_{\hs}$.
For the operator norm, the pointwise bound
$|x_i^\top Bx_i|\leq|B|_{\op}|x_i|^2$ gives
\[
    \mathbb E\left[|x_i^\top Bx_i|\right]
    \leq|B|_{\op}\mathbb E[|x_i|^2]
    =|B|_{\op}\tr(\Sigma_i).
\]
Since $\tr(\Sigma_iB)=\mathbb E[x_i^\top Bx_i]$, the triangle inequality gives
\[
    \left\|x_i^\top Bx_i-\tr(\Sigma_iB)\right\|_{L^1}
    \leq2\mathbb E\left[|x_i^\top Bx_i|\right].
\]
Taking the supremum over $i\in[n]$ and over the Hilbert--Schmidt or
operator unit ball proves the two modulus bounds.
\end{proof}

\begin{theorem}[Von Bahr--Esseen inequality {\cite{vonBahr-Esseen-65}}]
\label{the:von_bahr_esseen}
Let $1\leq q\leq2$, and let $Y_1,\ldots,Y_m$ be independent centered
real random variables with finite $q$th absolute moments. Then
\[
    \mathbb E\left|\sum_{j=1}^mY_j\right|^q
    \leq 2\sum_{j=1}^m\mathbb E|Y_j|^q.
\]
\end{theorem}

\begin{proposition}[Quadratic forms with independent coordinates]
\label{pro:quadratic_form_independent_coordinates_L1}
Let $x\in\mathbb R^p$ have centered independent coordinates, and put
$\Sigma:=\mathbb E[xx^\top]$.
Let \(r\geq1\). If
\[
    K_{2r}:=\sup_{i\in[p]}\|x_i\|_{L^{2r}}<\infty,
\]
then, for every deterministic $A\in\mathcal M_p$,
\begin{align*}
    \alpha&:=\max\left\{\frac1r-\frac12,0\right\},
    &
    \left\|x^{\top}Ax-\tr\bigl(A\Sigma\bigr)\right\|_{L^1}
    &= O\!\left(K_{2r}^2p^\alpha|A|_{\hs}\right).
\end{align*}
For $r=1$, one has the explicit bound
\[
    \left\|x^\top Ax
      -\tr\bigl(A\Sigma\bigr)\right\|_{L^1}
    \leq2K_2^2\sqrt p\,|A|_{\hs}.
\]
\end{proposition}

\begin{proof}[Proof of Proposition~\ref{pro:quadratic_form_independent_coordinates_L1}]
For $r=1$, the matrix $\Sigma$ is diagonal and has
Hilbert--Schmidt norm at most $K_2^2\sqrt p$. Thus
Lemma~\ref{lem:uncentered_quadratic_form_bound} and the triangle inequality give the explicit bound.

Now let $r>1$.
By Lemma~\ref{lem:complex_conditional_test_reduction}, we may replace $A$
by its complex symmetric part.  Write $A=D+R$, where
$D=\diag(d_1,\ldots,d_p)$ and $R$ has zero diagonal.  Independence and
centering give the off-diagonal identity
\begin{align}
 \mathbb E\left[|x^\top Rx|^2\right]
 &=4\sum_{1\leq i<j\leq p}|R_{ij}|^2
     \mathbb E[x_i^2]\mathbb E[x_j^2].
 \label{eq:independent_coordinate_off_diagonal_L2}
\end{align}
This identity uses only second moments; the absolute squares are retained
because $R$ may be complex.

Put $K_2:=\sup_{i\in[p]}\|x_i\|_{L^2}\leq K_{2r}$.
For the diagonal part, apply
Theorem~\ref{the:von_bahr_esseen} with \(q:=r\wedge2\).
Then \(1<q\leq2\) and \(2q\leq2r\), so
\[
    \left\|x_i^2-\mathbb E[x_i^2]\right\|_{L^q}
    \leq2\|x_i\|_{L^{2q}}^2
    \leq2K_{2r}^2.
\]
The centered squares are independent real random variables. Since the
coefficients $d_i$ may be complex, apply the theorem separately to
$\Re(d_i)(x_i^2-\mathbb E[x_i^2])$ and
$\Im(d_i)(x_i^2-\mathbb E[x_i^2])$, then combine the bounds by
the triangle inequality in $L^q$. Using
$\|\cdot\|_{L^1}\leq\|\cdot\|_{L^q}$ gives
\[
 \left\|x^\top Dx-\tr(D\Sigma)\right\|_{L^1}
 = O\!\left(K_{2r}^2\left(\sum_{i\in[p]}|d_i|^q\right)^{1/q}\right)
    = O\!\left(K_{2r}^2p^{\alpha}|D|_{\hs}\right).
\]
The last step follows from H\"older's inequality:
\[
    \left(\sum_{i\in[p]}|d_i|^q\right)^{1/q}
    \leq p^{1/q-1/2}\left(\sum_{i\in[p]}|d_i|^2\right)^{1/2}
    =p^\alpha|D|_{\hs},
\]
since $1/q-1/2=\alpha$.
Equation~\eqref{eq:independent_coordinate_off_diagonal_L2} gives
\[
 \mathbb E\left[|x^\top Rx|^2\right]
 \leq2K_2^4|R|_{\hs}^2.
\]
Consequently,
\[
    \|x^\top Rx\|_{L^1}
    \leq\sqrt2K_2^2|R|_{\hs}.
\]
Since $p^\alpha\geq1$, and
$|D|_{\hs}+|R|_{\hs}\leq\sqrt2|A|_{\hs}$, the two estimates prove the
proposition.
\end{proof}

\begin{proposition}[Quadratic forms under finite variance]
\label{pro:quadratic_form_finite_variance}
Let the coordinates of $x\in\mathbb R^p$ be independent copies of a
centered real random variable $\xi$ with $\sigma^2:=\mathbb E[\xi^2]<\infty$,
where $\sigma\geq0$.
For every deterministic $A\in\mathcal M_p$ and every $L>0$,
\begin{equation}
    \left\|x^\top Ax-\sigma^2\tr(A)\right\|_{L^1}
    \leq
    |A|_{\hs}\left(
        \sqrt2\,\sigma^2+L\sigma
        +2\sqrt p\,
          \mathbb E\!\left[\xi^2\mathbf1_{\{|\xi|>L\}}\right]
    \right).
    \label{eq:iid_finite_variance_quadratic_form_truncation}
\end{equation}
Consequently, if all entries of $X$ are independent copies of $\xi$, and $\xi$ has a fixed law with finite variance,
$M_{\hs}^{(1)}=o(\sqrt p)$ as $p\to\infty$.
\end{proposition}

\begin{proof}
Use the symmetrization and decomposition
$A=D+R$ from the proof of
Proposition~\ref{pro:quadratic_form_independent_coordinates_L1}, with
$D=\diag(d_1,\ldots,d_p)$ and $R$ having zero diagonal.
For the diagonal part, decompose
\begin{align*}
    x_j^2-\sigma^2
    ={}&
    \left(
      x_j^2\mathbf1_{\{|x_j|\leq L\}}
      -\mathbb E\left[\xi^2\mathbf1_{\{|\xi|\leq L\}}\right]
    \right)
    +\left(
      x_j^2\mathbf1_{\{|x_j|>L\}}
      -\mathbb E\left[\xi^2\mathbf1_{\{|\xi|>L\}}\right]
    \right).
\end{align*}
The sum of the bounded centered terms has $L^2$ norm at most
\begin{align*}
    \left(
      \sum_{j=1}^p|d_j|^2
      \mathbb E\left[\xi^4\mathbf1_{\{|\xi|\leq L\}}\right]
    \right)^{1/2}
    \leq L\sigma\,|D|_{\hs},
\end{align*}
whereas the tail sum has $L^1$ norm at most
\begin{align*}
    2\sum_{j=1}^p|d_j|
      \mathbb E\left[\xi^2\mathbf1_{\{|\xi|>L\}}\right]
    \leq
    2\sqrt p\,|D|_{\hs}
      \mathbb E\left[\xi^2\mathbf1_{\{|\xi|>L\}}\right].
\end{align*}
Equation~\eqref{eq:independent_coordinate_off_diagonal_L2} also gives
\[
    \|x^\top Rx\|_{L^1}\leq\sqrt2\,\sigma^2|R|_{\hs}.
\]
Since $|D|_{\hs},|R|_{\hs}\leq|A|_{\hs}$, combining these three bounds
proves \eqref{eq:iid_finite_variance_quadratic_form_truncation}.
Taking the supremum over \(|A|_{\hs}\leq1\) and dividing by
\(\sqrt p\) gives
\[
    \frac{M_{\hs}^{(1)}}{\sqrt p}
    \leq
    \frac{\sqrt2\,\sigma^2+L\sigma}{\sqrt p}
    +2\mathbb E\!\left[\xi^2\mathbf1_{\{|\xi|>L\}}\right].
\]
For each fixed \(L\), taking the upper limit as \(p\to\infty\) yields
\[
    \limsup_{p\to\infty}\frac{M_{\hs}^{(1)}}{\sqrt p}
    \leq
    2\mathbb E\!\left[\xi^2\mathbf1_{\{|\xi|>L\}}\right].
\]
Letting \(L\to\infty\) proves the conclusion, since \(\xi^2\)
is integrable.
\end{proof}

Recall that $B_{\mathrm s}=(B+B^{\top})/2$ denotes the symmetric part of
$B\in\mathcal M_p$.

\begin{lemma}[Geometry of the two test balls]
\label{lem:quadratic_form_test_ball_geometry}
For every real symmetric positive semidefinite matrix
$\Sigma\in\mathcal M_p$,
\begin{align*}
    \sup_{\substack{B\in\mathcal M_p\\|B|_{\hs}\leq1}}
    |\Sigma^{1/2}B_{\mathrm s}\Sigma^{1/2}|_{\hs}
    &=|\Sigma|_{\op},
    &
    \sup_{\substack{B\in\mathcal M_p\\|B|_{\op}\leq1}}
    |\Sigma^{1/2}B_{\mathrm s}\Sigma^{1/2}|_{\hs}
    &=|\Sigma|_{\hs}.
\end{align*}
\end{lemma}

\begin{proof}
Lemma~\ref{lem:complex_conditional_test_reduction} gives
$|B_{\mathrm s}|_{\hs}\leq |B|_{\hs}$.  Also,
$|AXB|_{\hs}\leq |A|_{\op}|X|_{\hs}|B|_{\op}$.  Hence
\[
    |\Sigma^{1/2}B_{\mathrm s}\Sigma^{1/2}|_{\hs}
    \leq |\Sigma|_{\op}|B_{\mathrm s}|_{\hs}
    \leq |\Sigma|_{\op}|B|_{\hs}.
\]
For the other test ball, work in a real orthonormal eigenbasis of
\(\Sigma\), with eigenvalues \((\lambda_j)_{j\in[p]}\). Then
\[
    |\Sigma^{1/2}B_{\mathrm s}\Sigma^{1/2}|_{\hs}^2
    =\sum_{j,k=1}^p
      \lambda_j\lambda_k|(B_{\mathrm s})_{jk}|^2
    \leq\frac12\sum_{j,k=1}^p
      (\lambda_j^2+\lambda_k^2)|(B_{\mathrm s})_{jk}|^2
    \leq |\Sigma|_{\hs}^2|B|_{\op}^2.
\]
Every row and column of $B_{\mathrm s}$ has Euclidean norm at most
$|B_{\mathrm s}|_{\op}\leq|B|_{\op}$.  This proves the last inequality.
For the reverse inequalities, take respectively
$B=vv^{\top}$, where $v$ is a unit eigenvector corresponding to
$|\Sigma|_{\op}$, and $B=I_p$.
\end{proof}

\begin{lemma}[Independent-coordinate quadratic-form moduli]
\label{lem:profile_quadratic_form_moduli}
Let $x\in\mathbb R^p$ have second-moment matrix
$\Sigma:=\mathbb E[xx^\top]$, and assume that $x=\Sigma^{1/2}y$,
where $y$ has independent, centered, standardized coordinates.
If
\[
    K_4:=\sup_{j\in[p]}\|y_j\|_{L^4}<\infty,
\]
then
\begin{align*}
    \sup_{\substack{B\in\mathcal M_p\\|B|_{\hs}\leq1}}
    \left\|x^\top Bx-\tr(\Sigma B)\right\|_{L^2}
    &=O\!\left(K_4^2|\Sigma|_{\op}\right),
    \\
    \sup_{\substack{B\in\mathcal M_p\\|B|_{\op}\leq1}}
    \left\|x^\top Bx-\tr(\Sigma B)\right\|_{L^2}
    &=O\!\left(K_4^2|\Sigma|_{\hs}\right).
\end{align*}
\end{lemma}

\begin{proof}
Quadratic forms depend only on $B_{\mathrm s}$ by
Lemma~\ref{lem:complex_conditional_test_reduction}. For a fixed test matrix
$B$, put
\[
    C:=\Sigma^{1/2}B_{\mathrm s}\Sigma^{1/2}.
\]
Then
\[
    x^\top Bx-\tr(\Sigma B)
    =y^\top Cy-\tr(C).
\]
Writing $C=(c_{jk})_{j,k\in[p]}$, independence, centering, and
standardization give
\[
 \left\|y^\top Cy-\tr(C)\right\|_{L^2}^2
    ={}
 \sum_{j=1}^p|c_{jj}|^2\Var(y_j^2)
 +4\sum_{1\leq j<k\leq p}|c_{jk}|^2
 \leq 2K_4^4|C|_{\hs}^2.
\]
Taking square roots and applying
Lemma~\ref{lem:quadratic_form_test_ball_geometry} to the two test
balls proves the result.
\end{proof}

For columns $x_i=\Sigma_i^{1/2}y_i$ satisfying these hypotheses with
$\sup_{i\in[n],\,j\in[p]}\|(y_i)_j\|_{L^4}\leq K_4<\infty$,
taking the maximum over $i\in[n]$ gives
\begin{align*}
    M_{\hs}^{(2)}
    &=O\!\left(K_4^2|\Sigma|_{\op}^\infty\right),
    &
    M_{\op}^{(2)}
    &=O\!\left(K_4^2|\Sigma|_{\hs}^\infty\right).
\end{align*}

\begin{lemma}[Gaussian quadratic-form moduli]
\label{lem:gaussian_quadratic_form_moduli}
Let $x\sim\mathcal N(0,\Sigma)$, where $\Sigma\in\mathcal M_p$ is real
symmetric and positive semidefinite.  For every fixed $q\geq2$,
\begin{align}
    \sqrt2|\Sigma|_{\op}
    &\leq
    \sup_{\substack{B\in\mathcal M_p\\|B|_{\hs}\leq1}}
    \left\|x^{\top}Bx-\tr(\Sigma B)\right\|_{L^q}
    = O(|\Sigma|_{\op}),
    \label{eq:gaussian_hs_test_ball_exact_scale}
    \\
    \sqrt2|\Sigma|_{\hs}
    &\leq
    \sup_{\substack{B\in\mathcal M_p\\|B|_{\op}\leq1}}
    \left\|x^{\top}Bx-\tr(\Sigma B)\right\|_{L^q}
    = O(|\Sigma|_{\hs}).
    \label{eq:gaussian_operator_test_ball_exact_scale}
\end{align}
For \(q=2\), both lower bounds are equalities.
\end{lemma}

\begin{proof}
Write $x=\Sigma^{1/2}g$, where $g\sim\mathcal N(0,I_p)$, and let
$B_{\mathrm s}=(B+B^{\top})/2$.  Wick's identity gives, also for complex
test matrices,
\begin{align*}
    \left\|x^{\top}Bx-\tr(\Sigma B)\right\|_{L^2}
    =\sqrt2\,
      |\Sigma^{1/2}B_{\mathrm s}\Sigma^{1/2}|_{\hs}.
\end{align*}
For $q\geq2$, monotonicity of the $L^q$ norms and Gaussian
hypercontractivity for second-order chaos yield
\[
    \sqrt2\,
      |\Sigma^{1/2}B_{\mathrm s}\Sigma^{1/2}|_{\hs}
    \leq
    \left\|x^{\top}Bx-\tr(\Sigma B)\right\|_{L^q}
    =
    O\!\left(|\Sigma^{1/2}B_{\mathrm s}\Sigma^{1/2}|_{\hs}\right).
\]
Taking the two test-ball suprema and applying
Lemma~\ref{lem:quadratic_form_test_ball_geometry} proves
\eqref{eq:gaussian_hs_test_ball_exact_scale} and
\eqref{eq:gaussian_operator_test_ball_exact_scale}.  At $q=2$, Wick's
identity and the two exact geometric suprema give the asserted identities.
\end{proof}

For Gaussian columns $x_i\sim\mathcal N(0,\Sigma_i)$, taking the maximum
over $i\in[n]$ gives
$\sqrt2|\Sigma|_{\op}^\infty\leq M_{\hs}^{(q)}
=O(|\Sigma|_{\op}^\infty)$ and
$\sqrt2|\Sigma|_{\hs}^\infty\leq M_{\op}^{(q)}
=O(|\Sigma|_{\hs}^\infty)$ for every fixed $q\geq2$.

The deterministic construction uses the following characterization;
see \cite{BOL97,kammoun2016no}.

\begin{theorem}[Classical Stieltjes-transform characterization]
\label{the:condition_stieltjes_transform}
Let $f$ be holomorphic on the upper half-plane, with $\Im f(z)>0$ there. If
\[
    \lim_{y\to+\infty}iy\,f(iy)=-1,
\]
then $f$ is the Stieltjes transform of a probability measure $\mu$ on
$\mathbb R$:
\[
    f(z)=\int_{\mathbb R}\frac{\mu(d\lambda)}{\lambda-z}.
\]
Moreover, for \(a<b\) such that
\(\mu(\{a\})=\mu(\{b\})=0\),
\[
    \mu([a,b])
    =\lim_{y\to0^+}\frac1\pi
      \int_a^b\Im(f(x+iy))\,dx.
\]
If, in addition, $\Im(zf(z))\geq0$ for every $z$ with $\Im z>0$, then
$\mu(\mathbb R\setminus\mathbb R_+)=0$, and $f$ admits a holomorphic
continuation to $\mathbb C\setminus\mathbb R_+$.
\end{theorem}

\bibliographystyle{alpha}
\bibliography{biblio} 

@article{HachemLoubatonNajimVallet2013,
  author  = {Hachem, Walid and Loubaton, Philippe and
             Najim, Jamal and Vallet, Pascal},
  title   = {On bilinear forms based on the resolvent of large random matrices},
  journal = {Annales de l'Institut Henri Poincar{\'e},
             Probabilit{\'e}s et Statistiques},
  year    = {2013},
  volume  = {49},
  number  = {1},
  pages   = {36--63},
  doi     = {10.1214/11-AIHP450},
  url     = {https://numdam.org/articles/10.1214/11-AIHP450/}
}

@article{KnowlesYin2017,
  author  = {Knowles, Antti and Yin, Jun},
  title   = {Anisotropic local laws for random matrices},
  journal = {Probability Theory and Related Fields},
  year    = {2017},
  volume  = {169},
  number  = {1--2},
  pages   = {257--352},
  doi     = {10.1007/s00440-016-0730-4}
}

@article{WagnerCouilletDebbahSlock2012,
  author  = {Wagner, Sebastian and Couillet, Romain and Debbah, M{\'e}rouane and Slock, Dirk T. M.},
  title   = {Large System Analysis of Linear Precoding in Correlated {MISO} Broadcast Channels under Limited Feedback},
  journal = {IEEE Transactions on Information Theory},
  volume  = {58},
  number  = {7},
  pages   = {4509--4537},
  year    = {2012},
  doi     = {10.1109/TIT.2012.2191700}
}

@misc{XiaoSlock2026,
  author        = {Xiao, Fangqing and Slock, Dirk T. M.},
  title         = {Gaussian-Process Dynamics of Diagonal Expectation Propagation under Variance-Profile {Gaussian} Measurements},
  year          = {2026},
  eprint        = {2606.04531},
  archivePrefix = {arXiv},
  primaryClass  = {eess.SP},
  note          = {arXiv:\allowbreak 2606.04531}
}

@article{Yin2020DifferentCorrelations,
  author  = {Yin, Yanqing},
  title   = {On the Singular Value Distribution of Large-Dimensional Data Matrices Whose Columns Have Different Correlations},
  journal = {Statistics},
  volume  = {54},
  number  = {2},
  pages   = {353--374},
  year    = {2020},
  doi     = {10.1080/02331888.2020.1720021}
}

@article{BenaychGeorgesCouillet2016,
  author  = {Benaych-Georges, Florent and Couillet, Romain},
  title   = {Spectral Analysis of the {Gram} Matrix of Mixture Models},
  journal = {ESAIM: Probability and Statistics},
  volume  = {20},
  pages   = {217--237},
  year    = {2016},
  doi     = {10.1051/ps/2016007}
}

@article{MeiWangYao2023,
  author  = {Mei, Tianxing and Wang, Chen and Yao, Jianfeng},
  title   = {On Singular Values of Data Matrices with General Independent Columns},
  journal = {The Annals of Statistics},
  volume  = {51},
  number  = {2},
  pages   = {624--645},
  year    = {2023},
  doi     = {10.1214/23-AOS2263}
}

@misc{Chouard2022,
  author        = {Chouard, Cl{\'e}ment},
  title         = {Quantitative Deterministic Equivalent of Sample Covariance Matrices with a General Dependence Structure},
  year          = {2022},
  eprint        = {2211.13044},
  archivePrefix = {arXiv},
  primaryClass  = {math.PR},
  note          = {arXiv:2211.13044}
}

@misc{FanMaPaquetteWang2026,
  author        = {Fan, Zhou and Ma, Renyuan and Paquette, Elliot and Wang, Zhichao},
  title         = {Anisotropic Local Law for Non-Separable Sample Covariance Matrices},
  year          = {2026},
  eprint        = {2602.17960},
  archivePrefix = {arXiv},
  primaryClass  = {math.PR},
  note          = {arXiv:2602.17960}
}

@misc{MaMisiakiewicz2026,
  author        = {Ma, Renyuan and Misiakiewicz, Theodor},
  title         = {The anisotropic local law for sample covariance matrices under quadratic-form concentration},
  year          = {2026},
  eprint        = {2609.09440},
  archivePrefix = {arXiv},
  primaryClass  = {math.PR},
  doi           = {10.48550/arXiv.2609.09440},
  note          = {arXiv:2609.09440}
}

@article{ZhuangZhangXuSong2026,
  author  = {Zhuang, Zeyan and Zhang, Xin and Xu, Dongfang and Song, Shenghui},
  title   = {No Eigenvalues Outside the Limiting Support of Generally Correlated and Noncentral Sample Covariance Matrices},
  journal = {IEEE Transactions on Information Theory},
  year    = {2026},
  doi     = {10.1109/TIT.2026.3703020}
}

@article{AltErdosKruger2017,
  author  = {Alt, Johannes and Erd{\H{o}}s, L{\'a}szl{\'o} and Kr{\"u}ger, Torben},
  title   = {Local Law for Random {Gram} Matrices},
  journal = {Electronic Journal of Probability},
  volume  = {22},
  number  = {25},
  pages   = {1--41},
  year    = {2017},
  doi     = {10.1214/17-EJP42}
}

@book{AjankiErdosKruger2019,
  author    = {Ajanki, Oskari H. and Erd{\H{o}}s, L{\'a}szl{\'o} and Kr{\"u}ger, Torben},
  title     = {Quadratic Vector Equations on Complex Upper Half-Plane},
  series    = {Memoirs of the American Mathematical Society},
  volume    = {261},
  publisher = {American Mathematical Society},
  year      = {2019},
  doi       = {10.1090/memo/1261}
}

@article{MarchenkoPastur1967,
  author  = {Mar{\v c}enko, V. A. and Pastur, L. A.},
  title   = {Distribution of eigenvalues for some sets of random matrices},
  journal = {Mathematics of the USSR-Sbornik},
  volume  = {1},
  number  = {4},
  pages   = {457--483},
  year    = {1967},
  doi     = {10.1070/SM1967v001n04ABEH001994}
}

@book{BaiSilverstein2010,
  author    = {Bai, Zhidong D. and Silverstein, Jack W.},
  title     = {Spectral Analysis of Large Dimensional Random Matrices},
  edition   = {2nd},
  publisher = {Springer},
  address   = {New York},
  year      = {2010}
}

@book{PasturShcherbina2011,
  author    = {Pastur, Leonid and Shcherbina, Mariya},
  title     = {Eigenvalue Distribution of Large Random Matrices},
  publisher = {American Mathematical Society},
  series    = {Mathematical Surveys and Monographs},
  volume    = {171},
  year      = {2011}
}

@article{DozierSilverstein2007,
  author  = {Dozier, R. Brent and Silverstein, Jack W.},
  title   = {On the empirical distribution of eigenvalues of large dimensional information-plus-noise-type matrices},
  journal = {Journal of Multivariate Analysis},
  volume  = {98},
  number  = {4},
  pages   = {678--694},
  year    = {2007},
  doi     = {10.1016/j.jmva.2006.09.006}
}

@article{HachemLoubatonNajim2007,
  author  = {Hachem, Walid and Loubaton, Philippe and Najim, Jamal},
  title   = {Deterministic equivalents for certain functionals of large random matrices},
  journal = {The Annals of Applied Probability},
  volume  = {17},
  number  = {3},
  pages   = {875--930},
  year    = {2007},
  doi     = {10.1214/105051606000000925}
}

@article{BaiZhou2008,
  author  = {Bai, Zhidong D. and Zhou, Wang},
  title   = {Large sample covariance matrices without independence structures in columns},
  journal = {Statistica Sinica},
  volume  = {18},
  number  = {2},
  pages   = {425--442},
  year    = {2008}
}

@article{Yaskov2016ECP,
  author  = {Yaskov, Pavel},
  title   = {Necessary and sufficient conditions for the {Marchenko--Pastur} theorem},
  journal = {Electronic Communications in Probability},
  volume  = {21},
  pages   = {Paper No. 73, 8},
  year    = {2016},
  doi     = {10.1214/16-ECP4748}
}

@article{DembczakLytova2022,
  author  = {Dembczak-Ko{\l}odziejczyk, Alicja and Lytova, Anna},
  title   = {On the empirical spectral distribution for certain models related to sample covariance matrices with different correlations},
  journal = {Random Matrices: Theory and Applications},
  volume  = {11},
  number  = {3},
  pages   = {2250030},
  year    = {2022},
  doi     = {10.1142/S2010326322500307}
}

@article{vonBahr-Esseen-65,
  author  = {von Bahr, Bengt and Esseen, Carl-Gustav},
  title   = {Inequalities for the $r$th absolute moment of a sum of random variables, $1 \le r \le 2$},
  journal = {Annals of Mathematical Statistics},
  year    = {1965},
  volume  = {36},
  number  = {1},
  pages   = {299--303},
  doi     = {10.1214/aoms/1177700291}
}

@incollection{Burkholder1988Sharp,
  author  = {Burkholder, Donald L.},
  title   = {Sharp inequalities for martingales and stochastic integrals},
  booktitle = {Colloque Paul L\'evy sur les processus stochastiques},
  series  = {Ast\'erisque},
  number  = {157--158},
  pages   = {75--94},
  year    = {1988},
  publisher = {Soci\'et\'e math\'ematique de France}
}

@article{kammoun2016no,
  title={No eigenvalues outside the limiting support of generally correlated {Gaussian} matrices},
  author={Kammoun, Abla and Alouini, Mohamed-Slim},
  journal={IEEE Transactions on Information Theory},
  volume={62},
  number={7},
  pages={4312--4326},
  year={2016},
  doi={10.1109/TIT.2016.2561998},
  publisher={IEEE}
}

@Article{BOL97,
author={Bolotnikov, Vladimir},
journal={Linear Algebra and its Applications},
title={On a General Moment Problem on the Half Axis},
year={1997},
volume={255},
pages={57--112}
}

@article{YaskovUniversality2014,
  title={The universality principle for spectral distributions of sample covariance matrices},
  author={Yaskov, Pavel},
  journal={arXiv preprint arXiv:1410.5190},
  year={2014}
}

@ARTICLE{LOU17b,
  title={A Random Matrix Approach to Neural Networks},
  author={Louart, Cosme and Liao, Zhenyu and Couillet, Romain},
  journal={Annals of Applied Probability},
  pages={1190--1248},
  volume={28},
  number={2},
  year={2018},
  doi={10.1214/17-AAP1328}
}

@ARTICLE{HOY13,
  title={Massive {MIMO} in the {UL}/{DL} of Cellular Networks: How Many Antennas Do We Need?},
  author={Hoydis, Jakob and ten Brink, Stephan and Debbah, M{\'e}rouane},
  journal={IEEE Journal on Selected Areas in Communications},
  volume={31},
  number={2},
  pages={160--171},
  year={2013},
  publisher={IEEE},
  doi={10.1109/JSAC.2013.130205}
}

@ARTICLE{ELK09,
  title={Concentration of measure and spectra of random matrices: applications to correlation matrices, elliptical distributions and beyond},
  author={El Karoui, Noureddine},
  journal={The Annals of Applied Probability},
  volume={19},
  number={6},
  pages={2362--2405},
  year={2009},
  doi={10.1214/08-AAP548},
  publisher={Institute of Mathematical Statistics}
}

@article{ChengMontanari2024,
  author        = {Cheng, Chen and Montanari, Andrea},
  title         = {Dimension free ridge regression},
  journal       = {The Annals of Statistics},
  year          = {2024},
  volume        = {52},
  number        = {6},
  pages         = {2879--2912},
  doi           = {10.1214/24-AOS2449},
  eprint        = {2210.08571},
  archivePrefix = {arXiv},
  primaryClass  = {math.ST}
}

@misc{MisiakiewiczSaeed2024,
  author        = {Misiakiewicz, Theodor and Saeed, Basil},
  title         = {A non-asymptotic theory of {Kernel Ridge Regression}:
                   deterministic equivalents, test error, and {GCV} estimator},
  year          = {2024},
  eprint        = {2403.08938},
  archivePrefix = {arXiv},
  primaryClass  = {stat.ML},
  note          = {arXiv:\allowbreak 2403.08938}
}

@article{DingWang2025,
  author        = {Ding, Xiucai and Wang, Zhenggang},
  title         = {Global and Local {CLTs} for Linear Spectral Statistics of
                   General Sample Covariance Matrices When the Dimension
                   Is Much Larger Than the Sample Size With Applications},
  journal       = {IEEE Transactions on Information Theory},
  year          = {2025},
  volume        = {71},
  number        = {8},
  pages         = {6248--6296},
  doi           = {10.1109/TIT.2025.3576191},
  eprint        = {2308.08646},
  archivePrefix = {arXiv},
  primaryClass  = {math.ST}
}
 
\end{document}